\documentclass[preprint,3p,12pt]{elsarticle}
\usepackage{mathrsfs}
\usepackage{amsmath}
\usepackage{stmaryrd}
\usepackage{bbding}
\usepackage{dcolumn}
\usepackage{graphicx}
\usepackage{amsfonts}
\usepackage{amssymb}
\usepackage{psfrag}
\usepackage{wrapfig}
\usepackage{subfigure}
\usepackage{makeidx}
\usepackage{bm}
\usepackage{epsf}
\usepackage{epsfig}
\usepackage{setspace}
\usepackage{graphicx}
\usepackage{epstopdf}
\usepackage{psfrag}
\usepackage{subfigure}
\usepackage{color}

\begin{document}
\title{A Few Benchmark Test Cases for Higher-order Euler Solvers}

\author[IAPCM]{Liang Pan}
\ead{panliangjlu@sina.com}

\author[IAPCM]{Jiequan Li\corref{cor1}}
\ead{li\_jiequan@iapcm.ac.cn}

\author[HKUST1,HKUST2]{Kun Xu}
\ead{makxu@ust.hk}

\address[IAPCM]{Institute of Applied Physics and Computational Mathematics, Beijing, China}
\address[HKUST1]{Department of Mathematics, Hong Kong University of Science and Technology, Clear Water Bay, Kowloon, Hong Kong}
\address[HKUST2]{Department of Mechanical and Aerospace Engineering, Hong Kong University of Science and Technology, Clear Water Bay, Kowloon, Hong Kong}
\cortext[cor1]{Corresponding author}

\begin{abstract}
There have been great efforts on the development of higher-order
numerical schemes for compressible Euler equations. The traditional
tests mostly targeting on the strong shock interactions alone may
not be adequate to test the performance of higher-order schemes.
This study will introduce a few test cases with a wide range of wave
structures for testing higher-order schemes. As  reference
solutions, all test cases will be calculated by our recently
developed two-stage fourth-order gas-kinetic scheme (GKS). All
examples are selected so that the numerical settings are very simple
and any high order accurate scheme can be straightly used for these
test cases, and compare their performance with the GKS solutions.
The examples include highly oscillatory solutions and the large
density ratio problem in one dimensional case; hurricane-like
solutions, interactions of planar contact discontinuities (the
composite of entropy wave and vortex sheets) sheets with large Mach
number asymptotic, interaction of planar rarefaction waves with
transition from continuous flows to the presence of shocks, and
other types of interactions of two-dimensional planar waves. The
numerical results from the fourth-order gas-kinetic scheme provide
reference solutions only. These benchmark test cases will help CFD
developers to validate and further develop their schemes to a higher
level of accuracy and robustness.
\end{abstract}

\maketitle

\section{Introduction}
In past decades, there have been tremendous efforts on designing
high-order accurate numerical schemes for compressible fluid flows
and great success has been achieved. High-order accurate numerical
schemes were pioneered by Lax and Wendroff \cite{Lax-Wendroff}, and
extended into the version of high resolution methods by Kolgan
\cite{Kolgan}, Boris \cite{Boris}, van Leer \cite{Van-Leer}, Harten
\cite{Harten} et al, and other higher order versions, such as
Essentially Non-Oscillatory (ENO) \cite{ENO1, ENO2}, Weighted
Essentially Non-Oscillatory (WENO) \cite{WENO1, WENO2},
Discontinuous Galerkin (DG) \cite{DG1, DG2, DG3} methods etc.

In the past decades, the evaluation of the performance of numerical
scheme was mostly based on the test cases with strong shocks for
capturing sharp shock transition, such as the blast wave
interaction, the forward step-facing flows, and the double Mach
reflection \cite{Woodward-Colella}. Now it is not a problem at all
for shock capturing scheme to get stable sharp shock transition.
However, with the further development of higher order numerical
methods and practical demands (such as turbulent flow simulations),
more challenging test problems for capturing multiple wave structure
are expected to be used. For testing higher-order schemes, the
setting of these cases should be sufficiently simple and easy for
coding, and avoid the possible pollution from the boundary condition
and curvilinear meshes. To introduce a few tests which can be
truthfully used to evaluate the performance of higher-order scheme
is the motivation for the current paper. Our selected examples
include the following: one-dimensional cases, two-dimensional
Riemann problems, and the conservation law with source terms. For
the one-dimensional problems, the first case is a highly oscillatory
shock-turbulence interaction problem, which is the extension of
Shu-Osher problem by Titarev and Toro \cite{Titarev-Toro} with much
more severe oscillations, and the second one is a large density
ratio problem with a very strong rarefaction wave in the solution
\cite{Tang-Liu}, which is used to test how a numerical scheme
capture strong waves. For the two-dimensional cases, four groups are
tested. (i) Hurricane-like solutions \cite{Zhang-Zheng,
2d-riemann3}, which are highly nontrivial two-dimensional
time-dependent solutions with one-point vacuum in the center and
rotational velocity field. It is proposed to test the preservation
of positivity and symmetry of the numerical scheme. (ii) The
interaction of planar contact discontinuities for different Mach
numbers. The multidimensional contact discontinuities are the
composite of entropy waves and vortex sheets. The simulation of such
cases have difficulties due to the strong shear effects.  Since the
large Mach number limits for these cases have explicit solutions
\cite{Sheng-Zhang, 2d-riemann3}, they are proposed here in order to
check the ability of the current scheme for capturing wave
structures of various scales and the asymptotic property. (iii)
Interaction of planar rarefaction waves with the transition from
continuous fluid flows to the presence of shocks. (iv) Further
interaction of planar shocks showing the Mach reflection phenomenon.
These two-dimensional problems  fall into the category of
two-dimensional Riemann problems proposed in \cite{2d-riemann0}. The
two-dimensional Riemann problems reveal almost all substantial wave
patterns of shock reflections, spiral formations, vortex-shock
interactions and so on, through the simple classification of initial
data. The rich wave configurations conjectured in \cite{2d-riemann0}
have been confirmed numerically by several subsequent works
\cite{2d-riemann2, 2d-riemann3, 2d-riemann4, 2d-riemann5}. Since the
formulation of these problems are extremely simple, there is no need
of complicated numerical boundary treatment and they are suitable as
benchmark tests. The case for the conservation law with source term
is also proposed.

In order to provide reference solutions for all these test cases. A
gas-kinetic scheme will be used to calculate the solutions in this
paper. Recently, based on the time-dependent flux function of the
generalized Riemann problem (GRP) solver \cite{GRP1,GRP2,GRP3}, a
two-stage fourth-order time-accurate discretization was developed
for Lax-Wendroff type flow solvers, particularly applied for the
hyperbolic conservation laws \cite{GRP-high}. The reason for the
success of a two-stage L-W type time stepping method in achieving a
fourth-order temporal accuracy is solely due to the use of both flux
function and its temporal derivative. In terms of the gas evolution
model, the gas-kinetic scheme provides a temporal accurate flux
function as well, even though it depends on time through a much more
complicated relaxation process from the kinetic to the hydrodynamic
scale physics than the time-dependent flux function of GRP. Based on
this time-stepping method and the second-order gas-kinetic solver
\cite{GKS-Xu1,GKS-Xu2}, a fourth-order gas-kinetic scheme was
constructed for the Euler and Navier-Stokes equations
\cite{GKS-high3}. In comparison with the formal one-stage
time-stepping third-order gas-kinetic solver \cite{GKS-high1,
GKS-high2}, the fourth-order scheme not only reduces the complexity
of the flux function, but also improves the accuracy of the scheme,
even though the third-order and fourth-order schemes take similar
computation cost. The robustness of the fourth-order gas-kinetic
scheme is as good as the second-order one. Numerical tests show that
the fourth-order scheme not only has the expected order of accuracy
for the smooth flows, but also has favorable shock capturing
property for the discontinuous solutions.

This paper is organized as follows. In Section 2, we will briefly
review the fourth-order gas-kinetic scheme. In Section 3, we select
several groups of problems to show the performance of the scheme.
The final conclusion is made in the last section.

\section{The review of two-stage fourth-order gas-kinetic scheme}

In this section, we will briefly review our recently developed
two-stage fourth-order gas-kinetic scheme. This scheme is developed
in the framework of finite volume scheme, and it contains three
standard ingredients: spatial data reconstruction, two-stage time
stepping discretization, and second-order gas-kinetic flux function.

\subsection{Spatial reconstruction}
The spatial reconstruction for the gas-kinetic scheme contains two
parts, i.e. initial data reconstruction and reconstruction for
equilibrium. \vspace{0.2cm}

In this paper, the fifth-order WENO method \cite{WENO2} is used for
the initial data reconstruction. Assume that $W$ are the macroscopic
flow variables that need to be reconstructed. $W_i$ are the cell
averaged values, and $W_i^r, W_i^l$ are the two values obtained by
the reconstruction at two ends of the $i$-th cell. The fifth-order
WENO reconstruction is given as follows
\begin{align*}
W_i^r=\sum_{k=0}^2w_kw_{k}^r,~~ W_i^l=\sum_{k=0}^2\widetilde{w}_k
w_{k}^l,
\end{align*}
where all quantities involved are taken as
\begin{align*}
w_{0}^r= \frac{1}{3}W_{i}+\frac{5}{6}W_{i+1}-\frac{1}{6}W_{i+2},&~~
w_{0}^l=\frac{11}{6}W_{i}-\frac{7}{6}W_{i+1}+\frac{1}{3}W_{i+2},\\
w_{1}^r=-\frac{1}{6}W_{i-1}+\frac{5}{6}W_{i}+\frac{1}{3}W_{i+1},&~~
w_{1}^l=\frac{1}{3}W_{i-1}+\frac{5}{6}W_{i}-\frac{1}{6}W_{i+1},\\
w_{2}^r=\frac{1}{3}W_{i-2}-\frac{7}{6}W_{i-1}+\frac{11}{6}W_{i},&~~
w_{2}^l=-\frac{1}{6}W_{i-2}+\frac{5}{6}W_{i-1}+\frac{1}{3}W_{i},
\end{align*}
and $w_k, \widetilde{w}_k, k=0,1,2$ are the nonlinear weights. The
most widely used is the WENO-JS non-linear weights \cite{WENO2},
which can be written as follows
\begin{align*}
\displaystyle w_k^{JS}=\frac{\alpha_k^{JS}}{\sum_{p=0}^2
\alpha_p^{JS}},~\ \
\alpha_k^{JS}=\frac{d_k}{(\epsilon+\beta_k)^2},~\ \
\widetilde{w}_k^{JS}=\frac{\widetilde{\alpha}_s^{JS}}{\sum_{p=0}^2
\widetilde{\alpha}_p^{JS}},\ \
\widetilde{\alpha}_k^{JS}=\frac{d_k}{(\epsilon+\beta_k)^2},
\end{align*}
where
\begin{align*}
\displaystyle d_0=\widetilde{d}_2=\frac{3}{10},\ \
d_1=\widetilde{d}_1=\frac{3}{5},\ \
d_2=\widetilde{d}_0=\frac{1}{10},\ \ \epsilon=10^{-6},
\end{align*}
and $\beta_k$ is the smooth indicator, and the basic idea for its
construction can be found in \cite{WENO1}. In order to achieve a
better performance of the WENO scheme near smooth extrema, WENO-Z
\cite{WENO-Z} and WENO-Z+ \cite{WENO-Z+} reconstruction were
developed. The only difference is the nonlinear weights. The
nonlinear weights for the WENO-Z method is written as
\begin{align*}
w_k^{Z}=\frac{\alpha_k^{Z}}{\sum_{k=0}^2 \alpha_k^{Z}},~\ \ \ \ \
\alpha_k^{Z}=d_k\Big[1+\Big(\frac{\delta}{\epsilon+\beta_k}\Big)^2],
\end{align*}
and the nonlinear weights for the WENO-Z+ method is
\begin{align*}
w_k^{Z+}=\frac{\alpha_k^{Z+}}{\sum_{k=0}^2 \alpha_k^{Z+}},~\ \ \ \
\alpha_k^{Z+}=d_k\Big[1+\Big(\frac{\delta+\varepsilon}{\epsilon+\beta_k}\Big)^2+\lambda\Big(\frac{\epsilon+\beta_k}{\delta+\varepsilon}\Big)\Big],
\end{align*}
where $\beta_k$ is the same local smoothness indicator as in
\cite{WENO2}, $\delta=|\beta_0-\beta_2|$ is used for the fifth-order
reconstruction, and $\lambda$ is a parameter for fine-tuning the
size of the weight of less smooth stencils. In the numerical tests,
without special statement, WENO-JS method will be used for initial
data reconstruction. \vspace{0.2cm}

After the initial date reconstruction, the reconstruction of
equilibrium part is presented. For the cell interface $x_{i+1/2}$,
the reconstructed variables at both sides of the cell interface are
denoted as $W_l, W_r$. According to the compatibility condition,
which will be given later, the macroscopic variables at the cell
interface is obtained and denoted as $W_0$. The conservative
variables around the cell interface can be expanded as
\begin{align*}
\overline{W}(x)=W_{0}+S_1(x-x_*)+\frac{1}{2}S_2(x-x_*)^2+\frac{1}{6}S_3(x-x_*)^3+\frac{1}{24}S_4(x-x_*)^4.
\end{align*}
With the following conditions,
\begin{align*}
\int_{I_{i+k}} \overline{W}(x)=W_{i+k}, k=-1,...,2,
\end{align*}
the derivatives are given by
\begin{align*}
\overline{W}_x=S_1=\big[-\frac{1}{12}(W_{i+2}-W_{i-1})+\frac{5}{4}(W_{i+1}-W_{i})\big]/\Delta
x.
\end{align*}

\subsection{Two-stage fourth-order temporal discretization}
The two-stage fourth-order temporal discretization was developed for
Lax-Wendroff flow solvers, and was originally applied for the
generalized Riemann problem solver (GRP) \cite{GRP-high} for
hyperbolic equations. In \cite{Seal2014, Seal2016},  multi-stage
multi-derivative time stepping methods  were proposed and developed
as well under different framework of flux evaluation. Consider the
following time-dependent equation
\begin{align}\label{pde}
\frac{\partial \textbf{w}}{\partial t}=\mathcal {L}(\textbf{w}),
\end{align}
with the initial condition at $t_n$, i.e.,
\begin{align}\label{pde2}
\textbf{w}(t=t_n)=\textbf{w}^n,
\end{align}
where $\mathcal {L}$ is an operator for spatial derivative of flux.
The time derivatives are obtained using the Cauchy-Kovalevskaya
method,
\begin{align*}
\frac{\partial \textbf{w}^n}{\partial t}=\mathcal{L}(\textbf{w}^n),~
\ \ \ \frac{\partial }{\partial t}\mathcal
{L}(\textbf{w}^n)=\frac{\partial }{\partial \textbf{w}}\mathcal
{L}(\textbf{w}^n)\mathcal {L}(\textbf{w}^n).
\end{align*}
Introducing an intermediate state at $t^*=t_n+\Delta t/2$,
\begin{align}\label{step1}
\textbf{w}^*=\textbf{w}^n+\frac{1}{2}\Delta t\mathcal
{L}(\textbf{w}^n)+\frac{1}{8}\Delta t^2\frac{\partial}{\partial
t}\mathcal{L}(\textbf{w}^n),
\end{align}
the corresponding time  derivatives are obtained as well for the
intermediate stage state,
\begin{align*}
\frac{\partial \textbf{w}^*}{\partial t}=\mathcal{L}(\textbf{w}^*),~
\frac{\partial }{\partial t}\mathcal
{L}(\textbf{w}^*)=\frac{\partial }{\partial \textbf{w}}\mathcal
{L}(\textbf{w}^*)\cdot \mathcal {L}(\textbf{w}^*).
\end{align*}
Then, a fourth-order temporal accuracy solution of $\textbf{w}(t)$
at $t=t_n +\Delta t$ can be provided by the following equation
\begin{align}\label{step2}
\textbf{w}^{n+1}=\textbf{w}^n+\Delta t\mathcal
{L}(\textbf{w}^n)+\frac{1}{6}\Delta t^2\big(\frac{\partial}{\partial
t}\mathcal{L}(\textbf{w}^n)+2\frac{\partial}{\partial
t}\mathcal{L}(\textbf{w}^*)\big).
\end{align}
The details of this analysis can be found in \cite{GRP-high}. Thus,
a fourth-order temporal accuracy can be achieved by the two-stage
discretization Eq.\eqref{step1} and Eq.\eqref{step2}.

Consider the following conservation laws
\begin{align*}
\frac{\partial \textbf{w}}{\partial t}+\frac{\partial\textbf{
f}(\textbf{w})}{\partial x}=0.
\end{align*}
The semi-discrete form of a finite volume scheme can be written as
\begin{align*}
\frac{\partial \textbf{w}_i}{\partial t}=\mathcal
{L}_i(\textbf{w})=-\frac{1}{\Delta
x}(\textbf{f}_{i+1/2}-\textbf{f}_{i-1/2}),
\end{align*}
where $\textbf{w}_i$ are the cell averaged conservative variables,
$\textbf{f}_{i+1/2}$ are the fluxes at the cell interface
$x=x_{i+1/2}$, and $\Delta x$ is the cell size. With the temporal
derivatives of the flux, the two-stage fourth-order scheme can be
developed \cite{GKS-high3,GRP-high}. Similarly, for the conservation
laws with source terms
\begin{align*}
\frac{\partial \textbf{w}}{\partial
t}+\frac{\partial\textbf{f}(\textbf{w})}{\partial x}=S(\textbf{w}),
\end{align*}
the corresponding operator can be denoted as
\begin{align}\label{source}
\mathcal {L}_i (\textbf{w})=-\frac{1}{\Delta
x}(\textbf{f}_{i+1/2}-\textbf{f}_{i-1/2})+S(\textbf{w}_i).
\end{align}
The two-stage fourth-order temporal discretization can be directly
extended for conservation laws with source terms.

\subsection{Second-order gas-kinetic flux solver}
The two-dimensional BGK equation \cite{BGK-1,BGK-2} can be written
as
\begin{equation}\label{bgk}
f_t+\textbf{u}\cdot\nabla f=\frac{g-f}{\tau},
\end{equation}
where $f$ is the gas distribution function, $g$ is the corresponding
equilibrium state, and $\tau$ is the collision time. The collision
term satisfies the compatibility condition
\begin{equation}\label{compatibility}
\int \frac{g-f}{\tau}\psi d\Xi=0,
\end{equation}
where $\psi=(1,u,v,\displaystyle \frac{1}{2}(u^2+v^2+\xi^2))$,
$d\Xi=dudvd\xi^1...d\xi^{K}$, $K$ is number of internal freedom,
i.e. $K=(4-2\gamma)/(\gamma-1)$ for two-dimensional flows, and
$\gamma$ is the specific heat ratio. \vspace{0.2cm}

To update the flow variables in the finite volume framework, the
integral solution of BGK equation Eq.\eqref{bgk} is used to
construct the gas distribution function at a cell interface, which
can be written as
\begin{equation}\label{integral1}
f(x_{i+1/2},t,u,v,\xi)=\frac{1}{\tau}\int_0^t g(x',y',t',u,v,\xi)e^{-(t-t')/\tau}dt'\\
+e^{-t/\tau}f_0(-ut,y-vt,u,v,\xi),
\end{equation}
where $x_{i+1/2}=0$ is the location of the cell interface,
$x_{i+1/2}=x'+u(t-t')$ and $y=y'+v(t-t')$ are the trajectory of
particles, $f_0$ is the initial gas distribution function, and $g$
is the corresponding equilibrium state. The time dependent integral
solution at the cell interface $x_{i+1/2}$ can be expressed as
\begin{align}\label{flux}
f(x_{i+1/2},t,u,v,\xi)=&(1-e^{-t/\tau})g_0+((t+\tau)e^{-t/\tau}-\tau)(\overline{a}_1u+\overline{a}_2v)g_0\nonumber\\
+&(t-\tau+\tau e^{-t/\tau}){\bar{A}} g_0\nonumber\\
+&e^{-t/\tau}g_r[1-(\tau+t)(a_{1r}u+a_{2r}v)-\tau A_r)]H(u)\nonumber\\
+&e^{-t/\tau}g_l[1-(\tau+t)(a_{1l}u+a_{2l}v)-\tau A_l)](1-H(u)).
\end{align}
Based on the spatial reconstruction of macroscopic flow variables,
which is presented before, the conservative variables $W_l$ and
$W_r$ on the left and right hand sides of a cell interface, and the
corresponding equilibrium states $g_l$ and $g_r$, can be determined.
The conservative variables $W_{0}$ and the equilibrium state $g_{0}$
at the cell interface can be determined according to the
compatibility condition Eq.\eqref{compatibility} as follows
\begin{align*}
\int\psi g_{0}d\Xi=W_0=\int_{u>0}\psi g_{l}d\Xi+\int_{u<0}\psi
g_{r}d\Xi.
\end{align*}
The coefficients related to the spatial derivatives and time
derivative $a_{1k}, a_{2k}, A_k, k=l, r$ and $\overline{a}_1,
\overline{a}_2, \overline{A}$ in gas distribution function
Eq.\eqref{flux} can be determined according to the spatial
derivatives and compatibility condition. More details of the
gas-kinetic scheme can be found in \cite{GKS-Xu1}. \vspace{0.2cm}

As mentioned in the section before, in order to utilize the
two-stage temporal discretization, the temporal derivatives of the
flux function need to be determined. While in order to obtain the
temporal derivatives at $t_n$ and $t_*=t_n + \Delta t/2$ with the
correct physics, the flux function should be approximated as a
linear function of time within the time interval. Let's first
introduce the following notation,
\begin{align*}
\mathbb{F}_{i+1/2}(W^n,\delta)
=\int_{t_n}^{t_n+\delta}F_{i+1/2}(W^n,t)dt&=\int_{t_n}^{t_n+\delta}\int
u \psi f(x_{i+1/2},t,u, v,\xi)dud\xi dt.
\end{align*}
In the time interval $[t_n, t_n+\Delta t]$, the flux is expanded as
the following linear form
\begin{align}\label{expansion}
F_{i+1/2}(W^n,t)=F_{i+1/2}^n+ \partial_t F_{j+1/2}^n(t-t_n).
\end{align}
The coefficients $F_{j+1/2}^n$ and $\partial_tF_{j+1/2}^n$ can be
fully determined as follows
\begin{align*}
F_{i+1/2}(W^n,t_n)\Delta t&+\frac{1}{2}\partial_t
F_{i+1/2}(W^n,t_n)\Delta t^2 =\mathbb{F}_{i+1/2}(W^n,\Delta t) , \\
\frac{1}{2}F_{i+1/2}(W^n,t_n)\Delta t&+\frac{1}{8}\partial_t
F_{i+1/2}(W^n,t_n)\Delta t^2 =\mathbb{F}_{i+1/2}(W^n,\Delta t/2).
\end{align*}
Similarly, $\displaystyle F_{i+1/2}(W^*,t_*), \partial_t
F_{i+1/2}(W^*,t_*)$ for the intermediate stage $t_n+\Delta t/2$ can
be constructed as well.

Thus, we have completed the review of three building block of
fourth-order gas-kinetic scheme, i.e. spatial data reconstruction,
two-stage temporal discretization, and second-order gas-kinetic flux
solver. More details about the implementation of the fourth-order
gas-kinetic scheme can be found in \cite{GKS-high3}.

\section{Numerical experiments}
In this section, we mainly focus on the numerical solutions of the
compressible Euler equations. Based on the Chapman-Enskog expansion,
the macroscopic equations can be derived from the BGK equation
Eq.\eqref{bgk}. If the zeroth order truncation is taken, i.e. $f=g$,
the following two-dimensional Euler equations can be obtained
\begin{align}
\label{eq:Euler} \frac{\partial}{\partial{t}}
\begin{pmatrix}
   \rho \\
   \rho U \\
   \rho V \\
   \rho E \\
 \end{pmatrix}+
\frac{\partial}{\partial{x}}
\begin{pmatrix}
   \rho U\\
   \rho U^2+p \\
   \rho UV \\
    U(\rho E+p) \\
 \end{pmatrix}+
\frac{\partial}{\partial{y}}
\begin{pmatrix}
   \rho V\\
   \rho UV \\
   \rho V^2+p \\
   V(\rho E+p) \\
 \end{pmatrix}=0,
\end{align}
where $\rho
E=\displaystyle\frac{1}{2}\rho(U^2+V^2)+\frac{p}{\gamma-1}$. In this
section, a number of benchmark tests will be selected to test the
performance of the fourth-order gas-kinetic scheme.

\begin{figure}[!h]
\centering
\includegraphics[width=0.475\textwidth]{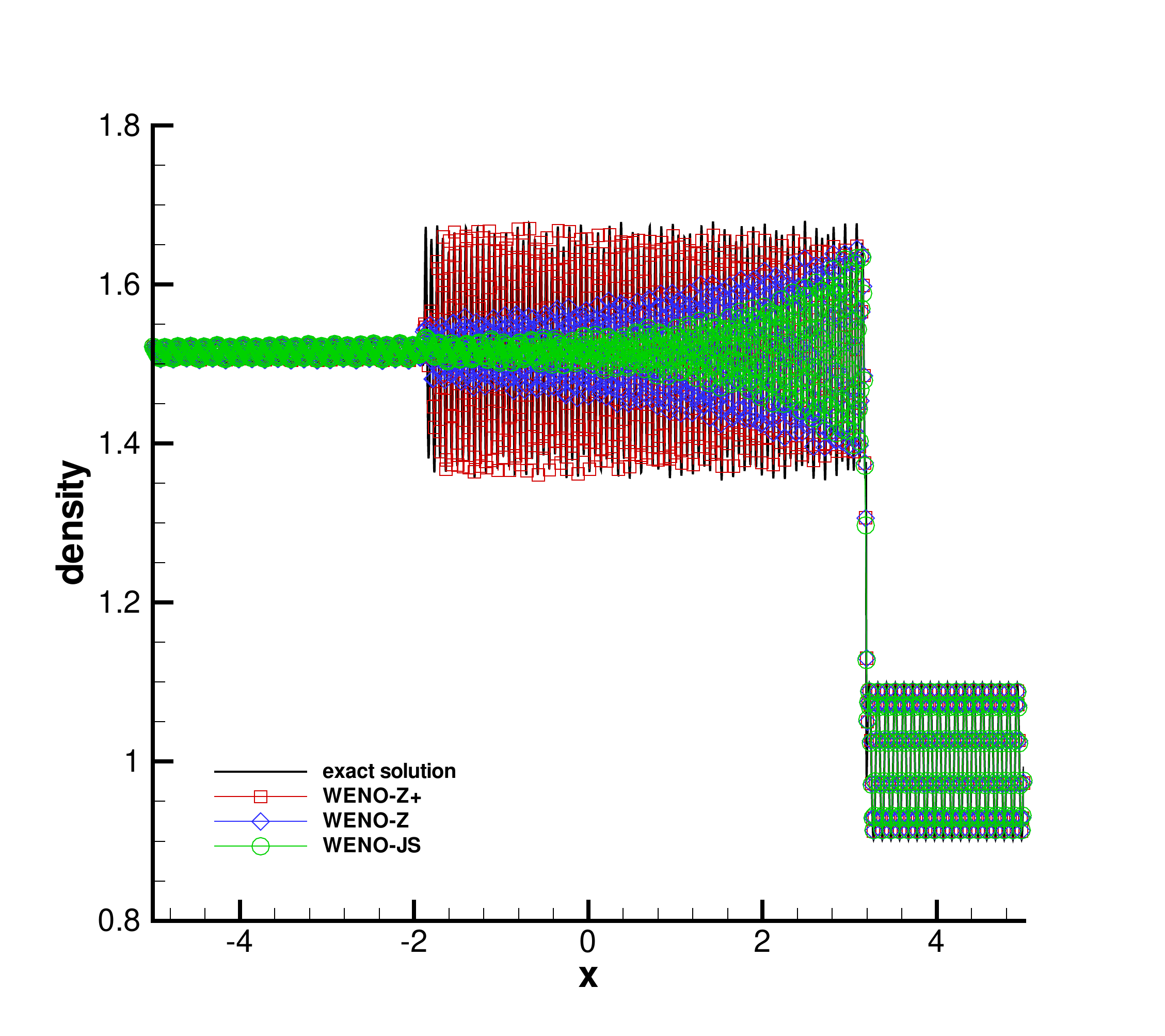}
\includegraphics[width=0.475\textwidth]{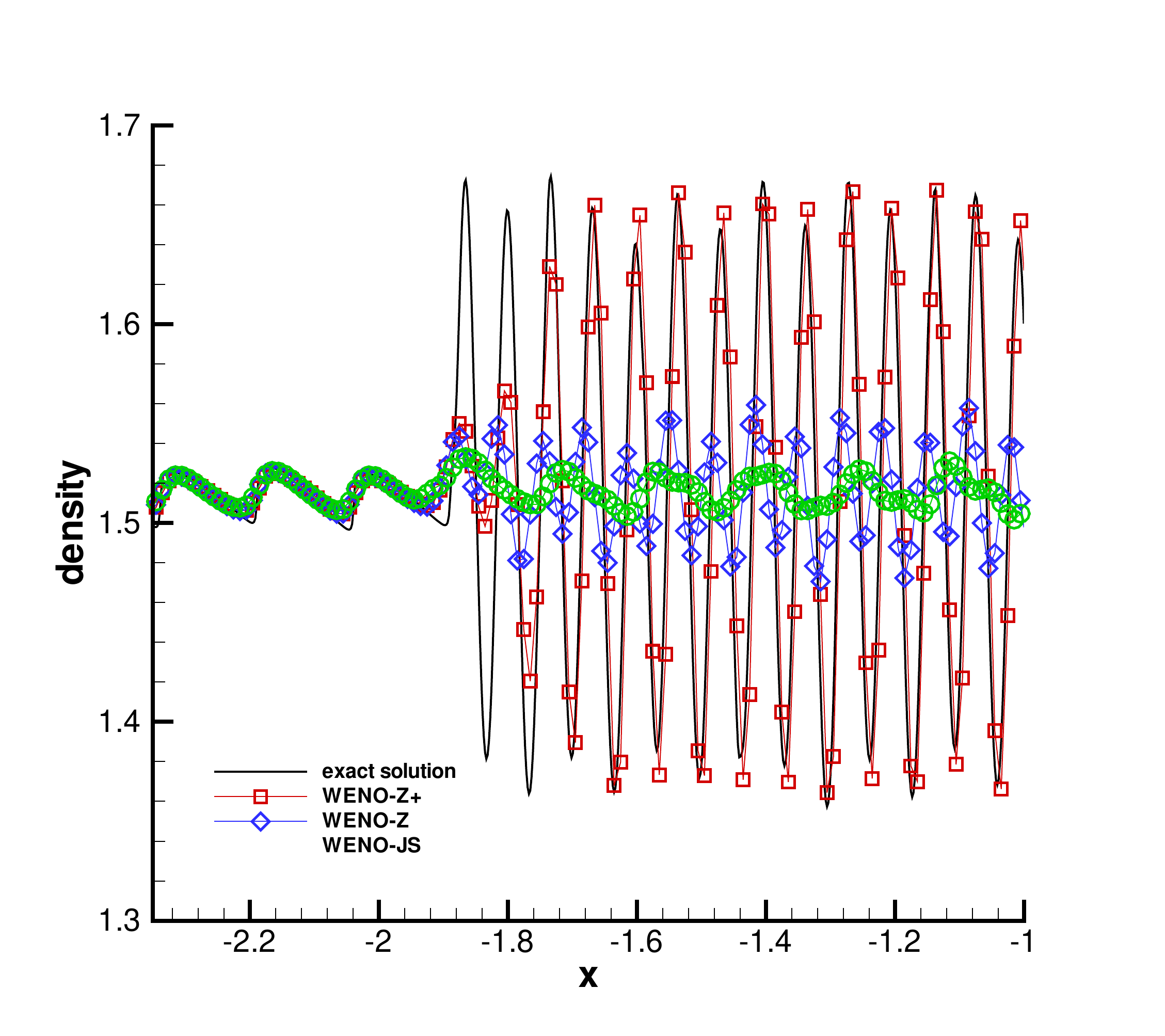}
\caption{\label{toro-1} The density distribution and local
enlargement of Titarev-Toro problem with WENO-JS, WENO-Z and WENO-Z+
weights at $t=5$, where $\lambda=\Delta x^{3/4}$ in WENO-Z+
methods.} \centering
\includegraphics[width=0.475\textwidth]{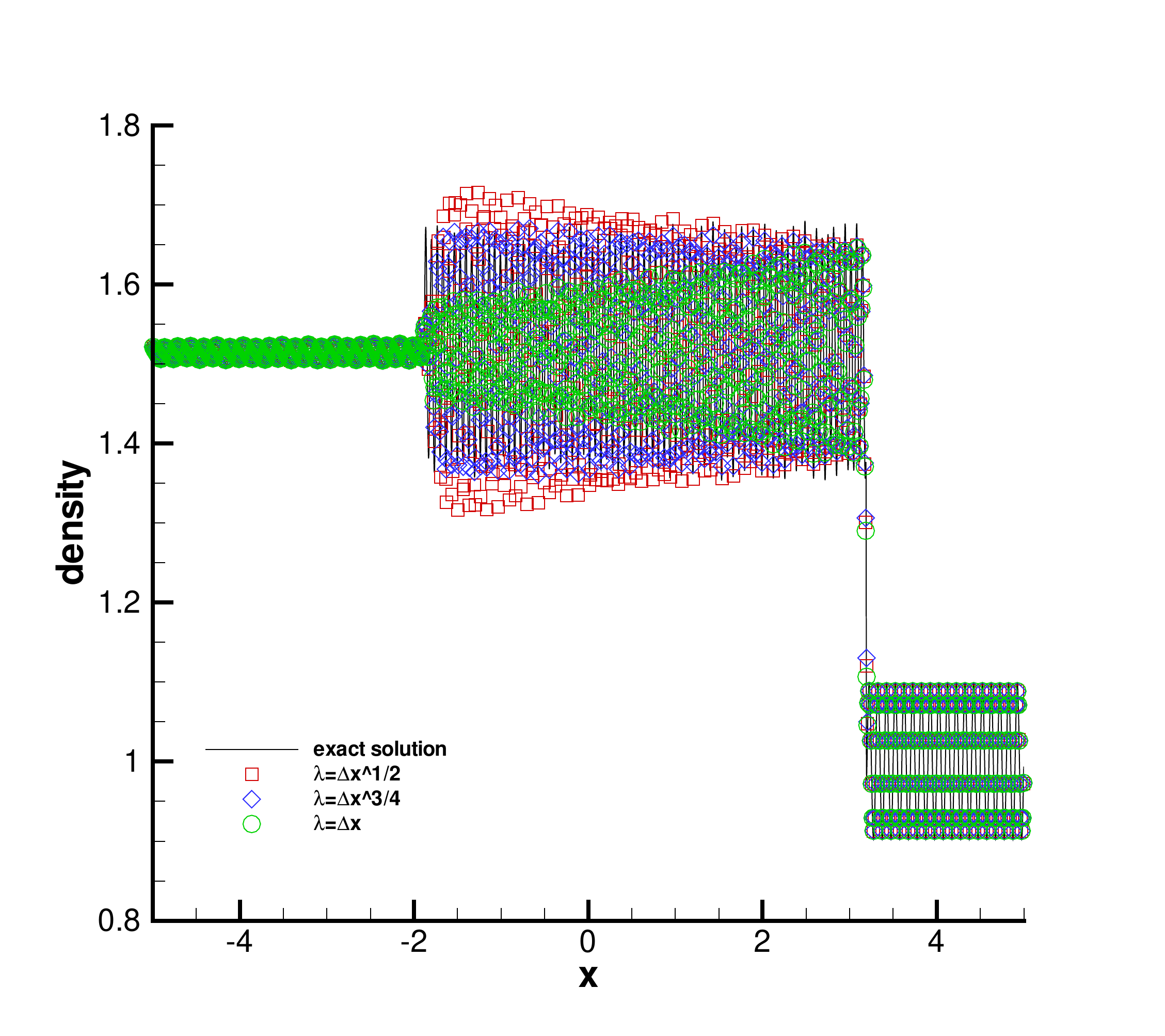}
\includegraphics[width=0.475\textwidth]{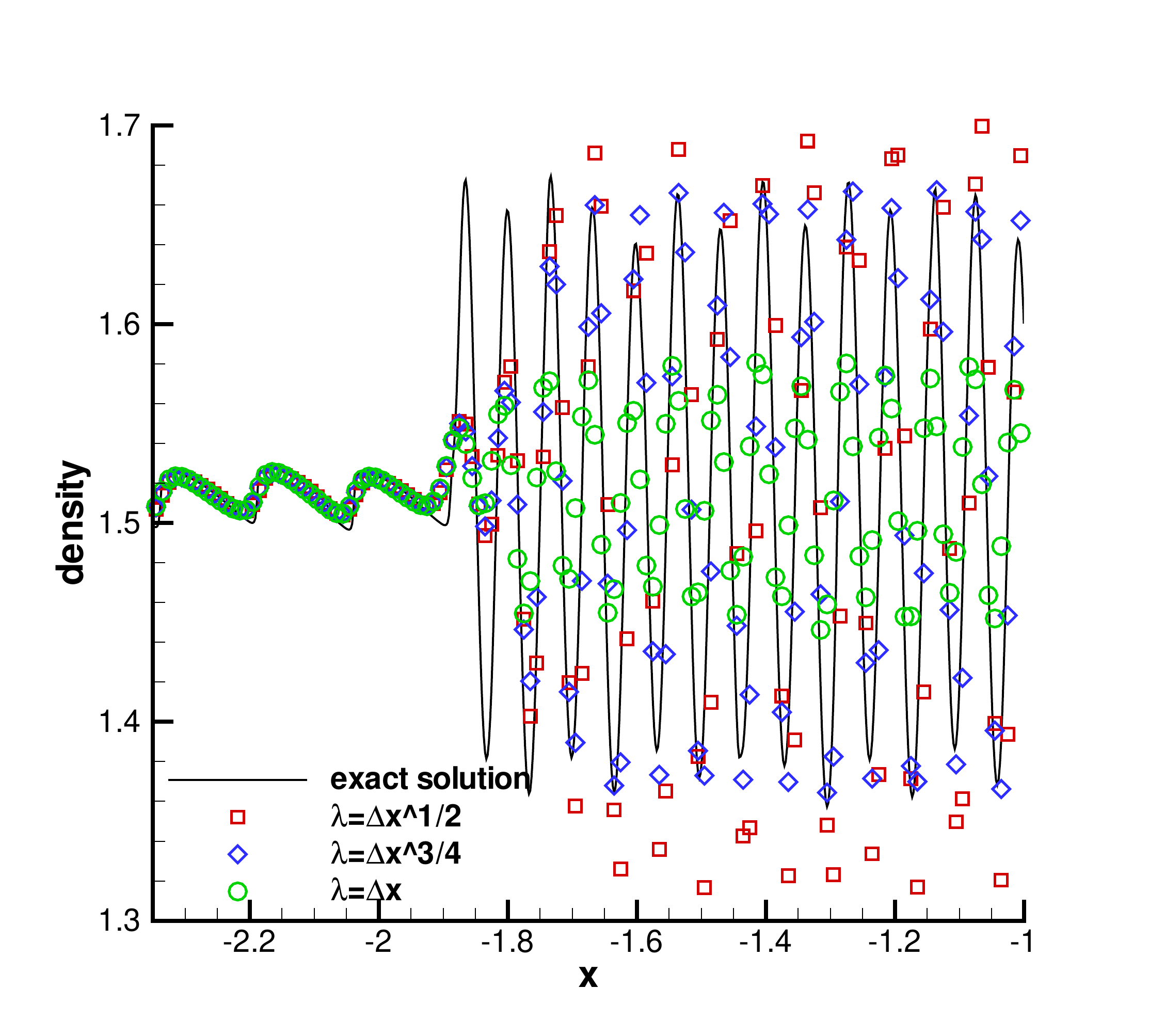}
\caption{\label{toro-2} The density distribution and local
enlargement of Titarev-Toro problem at $t=5$ for WENO-Z+ weights
with $\lambda=\Delta x^{1/2}, \Delta x^{3/4}$ and $\Delta x$.}
\end{figure}

For the inviscid flow, the collision time $\tau$ takes
\begin{align*}
\tau=\epsilon\Delta t+C \displaystyle|\frac{p_l-p_r}{p_l+p_r}|\Delta
t,
\end{align*}
where $\varepsilon=0.05$, $C=1$, $p_l$ and $p_r$ denote the pressure
on the left and right sides of the cell interface. The reason for
including artificial dissipation through the additional term in the
particle collision time is to enlarge the kinetic scale physics in
the discontinuous region for the construction of a numerical shock
structure through the particle free transport and inadequate
particle collision for keeping the non-equilibrium property in the
shock region. In all simulations, the Courant number $CFL=0.4$.
Without special statement, the specific heat ratio $\gamma=1.4$.

\subsection{One-dimensional problems}
The first one-dimensional problem is the extension of the Shu-Osher
problem given by Titarev and Toro \cite{Titarev-Toro} to test a
severely oscillatory wave interacting with shock. It aims to test
the ability of high-order numerical scheme to capture the extremely
high frequency waves. The initial condition for this cases is given
as follows,
\begin{equation*}
(\rho,U,p)=\left\{\begin{array}{ll}
(1.515695,0.523346,1.805),  \ \ \ \ & -5< x \leq -4.5,\\
(1 + 0.1\sin (20\pi x), 0, 1),  &  -4.5 <x <5.
\end{array} \right.
\end{equation*}
In the computation, $1000$ uniform cells are used, and the density
distributions at $t =5$ are presented. WENO-JS, WENO-Z and WENO-Z+
methods are used to test the performance of different nonlinear
weights in the WENO reconstruction. The density distribution and
local enlargement of Titarev-Toro problem with different weights are
given in Fig.\ref{toro-1}, where $\lambda=\Delta x^{3/4}$ in WENO-Z+
methods. Similar to the results in \cite{WENO-Z+}, with a proper
choice of the parameter $\lambda$, WENO-Z+ method performs much
better than WENO-JS and WENO-Z methods. In order to check the role
of $\lambda$ in WENO-Z+ method, the density distributions are
presented in Fig.\ref{toro-2}, for the cases with $\lambda=\Delta
x^{1/2}, \Delta x^{3/4}$ and $\Delta x$. It shows that the numerical
performance is very sensitive withe choice of $\lambda$, which is
also observed in \cite{WENO-Z+}. However, the current scheme tends
to provide better results than those in \cite{WENO-Z+} in all cases
with $\Delta x^{3/4}, \lambda=\Delta x^{1/2}$, and $\Delta x$ due to
the proper dynamics of the gas-kinetic flux solver.

\begin{figure}[!h]
\centering
\includegraphics[width=0.475\textwidth]{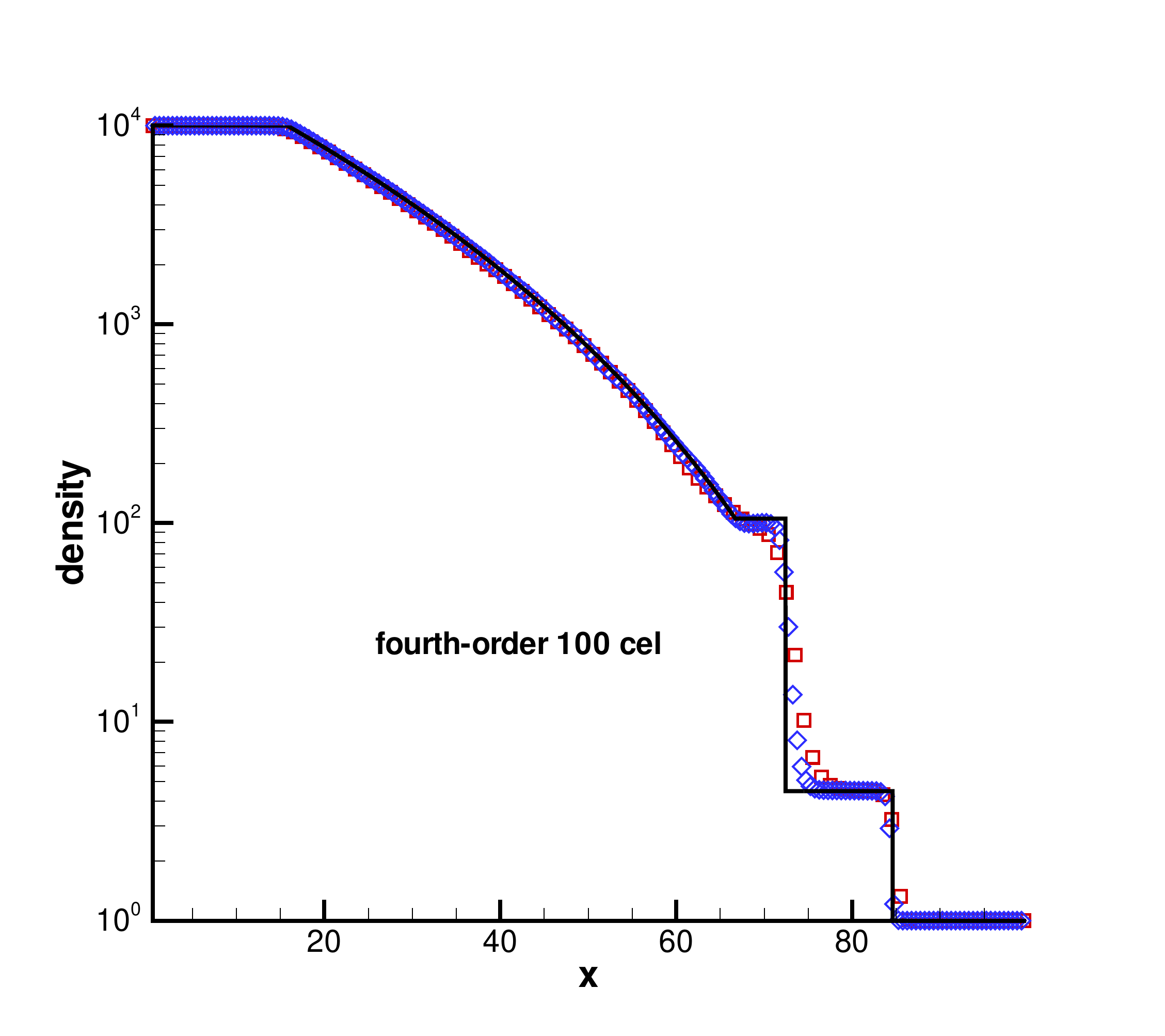}
\includegraphics[width=0.475\textwidth]{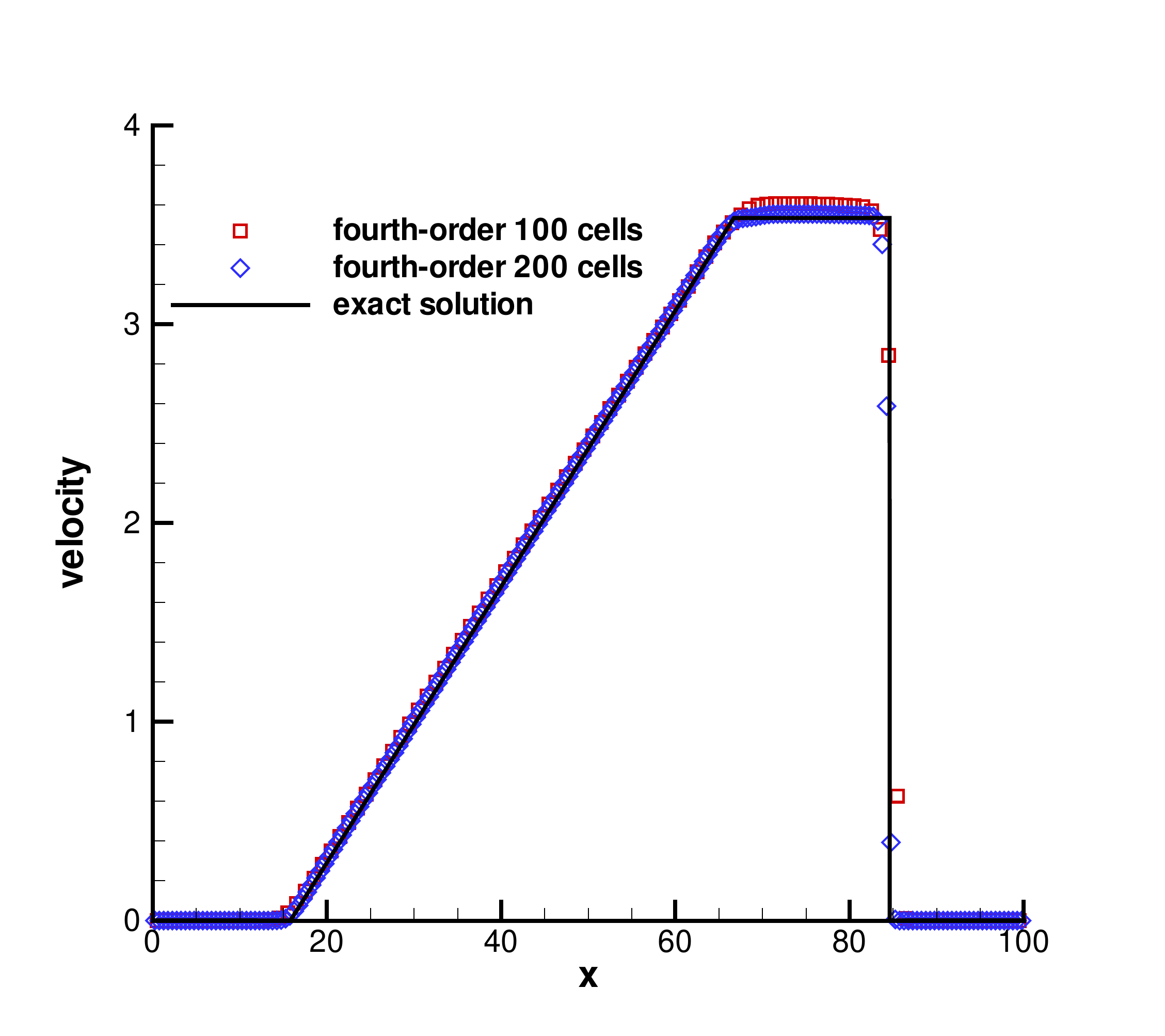}
\caption{\label{large-density-1}Large density ratio problem with
density variation from $10000$ to $1$. } \centering
\includegraphics[width=0.475\textwidth]{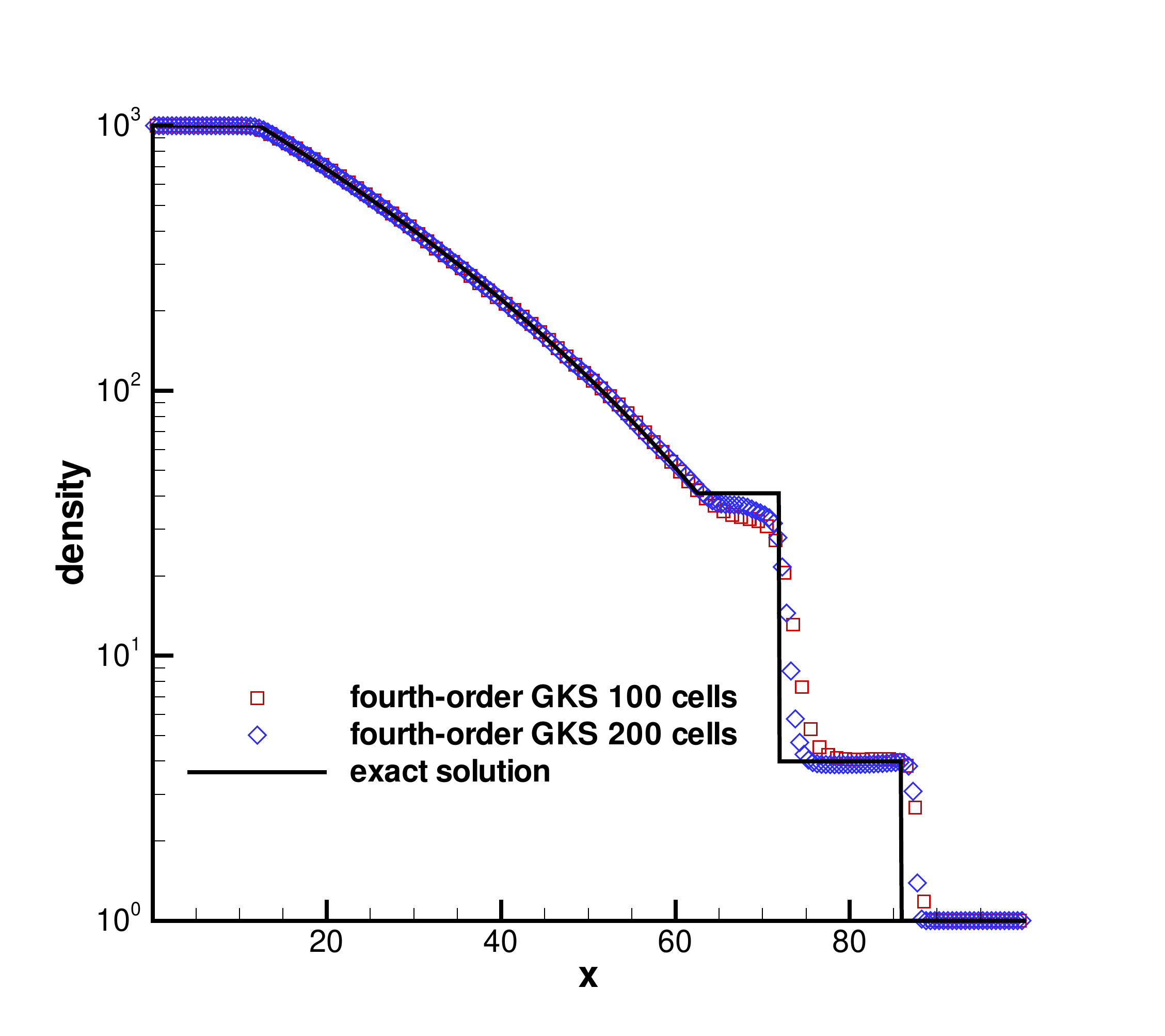}
\includegraphics[width=0.475\textwidth]{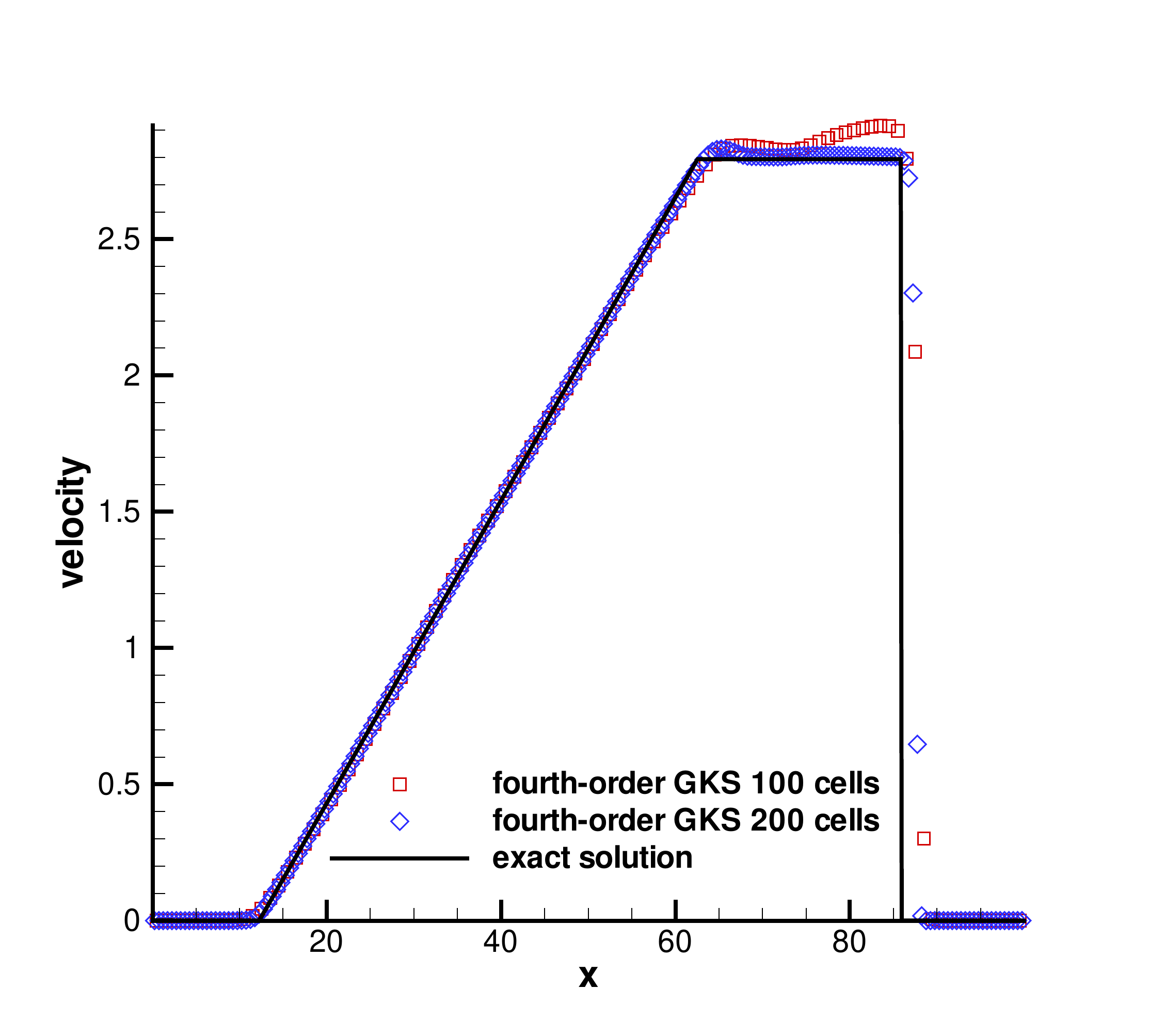}
\caption{\label{large-density-2}Large density ratio problem with
density variation from $1000$ to $1$. }
\end{figure}

The second example is the large density ratio problem with a very
strong rarefaction wave \cite{Tang-Liu}. It is proposed to test the
ability of a scheme for capturing strong waves. The initial data
\begin{equation*}
(\rho, U,p)=\left\{
\begin{array}{ll}
(10000,0,10000),  & 0\leq x<0.3, \\
(1,0,1), & 0.3<x<1.0.
\end{array}
\right.
\end{equation*}
It is a Sod-type problem, and the solution consists of a very strong
rarefaction wave with the density variation from $10000$ to $100$, a
contact discontinuity jumping from $100$ to $8$, and a shock wave
jumping from $1$ to $8$. If uniform grid points are used, the jump
of neighboring density values inside the rarefaction wave is about
$100$ times the density jump in the shock layer. Hence the strong
rarefaction wave is in a highly non-equilibrium state. It was
checked in \cite{Tang-Liu} that many high order numerical methods
have to use very fine grid points to capture the rarefaction wave,
the location of the shock and the contact discontinuity well. In
this case, we start with the exact solution at time $t=1.2$ as an
initial condition, and the output time for the solution is $t=12$,
with $100$ and $200$ uniform mesh points. The density and velocity
distributions are shown in Fig.\ref{large-density-1}. If we start
with the initial data at $t=0$, the numerical solution does not
display so well with the same grid points although it is still
acceptable. However, as the rarefaction wave is not so strong, the
scheme works well, which is verified through another example with
the density and pressure jumps from $10000$ to $1000$. The density
and velocity distributions are shown in Fig. \ref{large-density-2}.
The GRP simulation works better for such a case because the exact
GRP solver is used for the Euler equations \cite{GRP-high},
especially at the initial time for the solution starting from the
singular point. This example shows that the flux is extremely
important in the capturing of strong waves, refuting some misleading
statement that ``the construction of numerical flux function is not
so essential in the computation of fluid flows".

\subsection{Two-dimensional Riemann problems}
In this subsection, five groups of two-dimensional benchmark
problems are presented to illustrate the performance of the current
scheme.

\begin{figure}[!htb]
\centering
\includegraphics[width=0.475\textwidth]{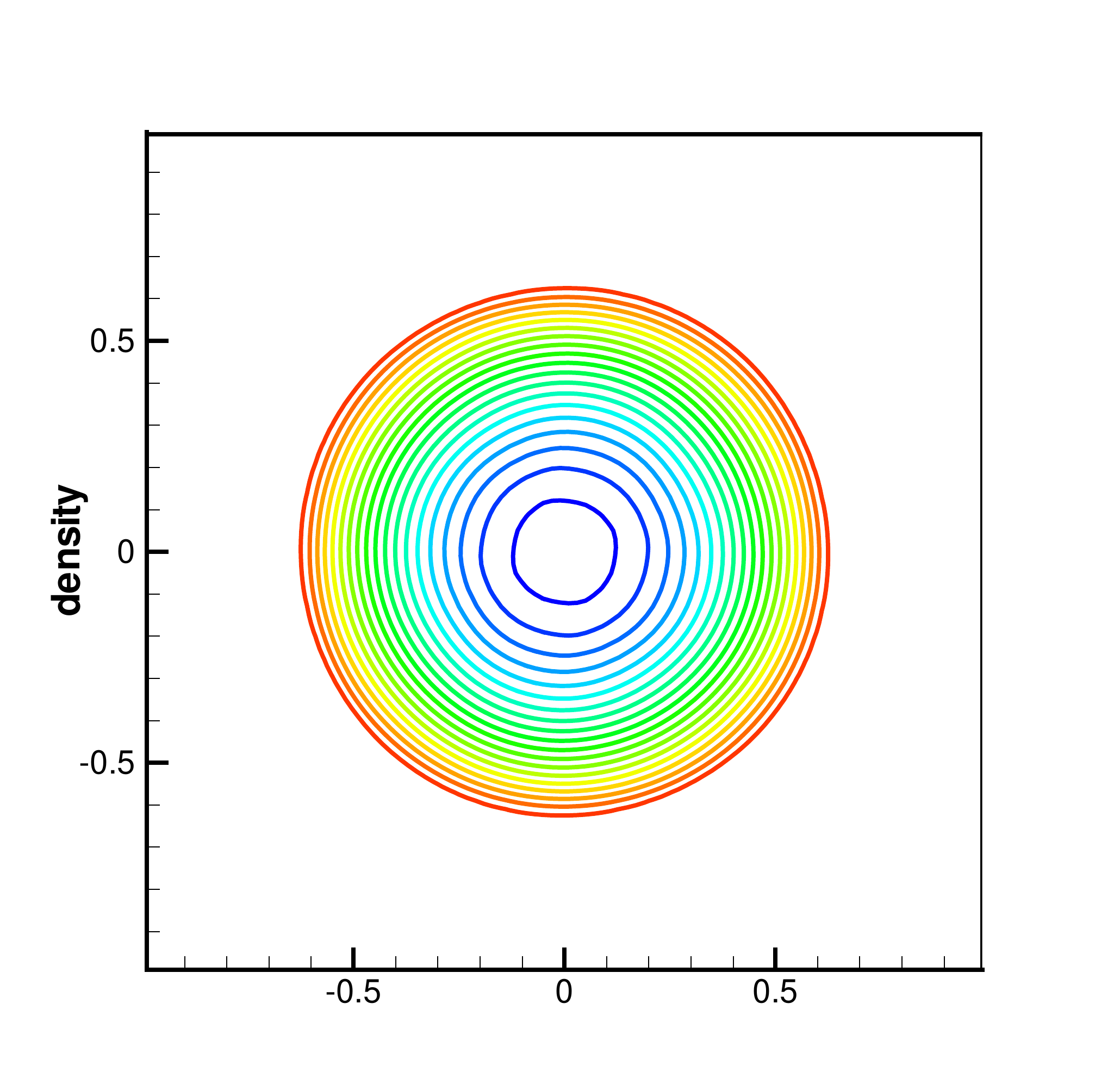}
\includegraphics[width=0.475\textwidth]{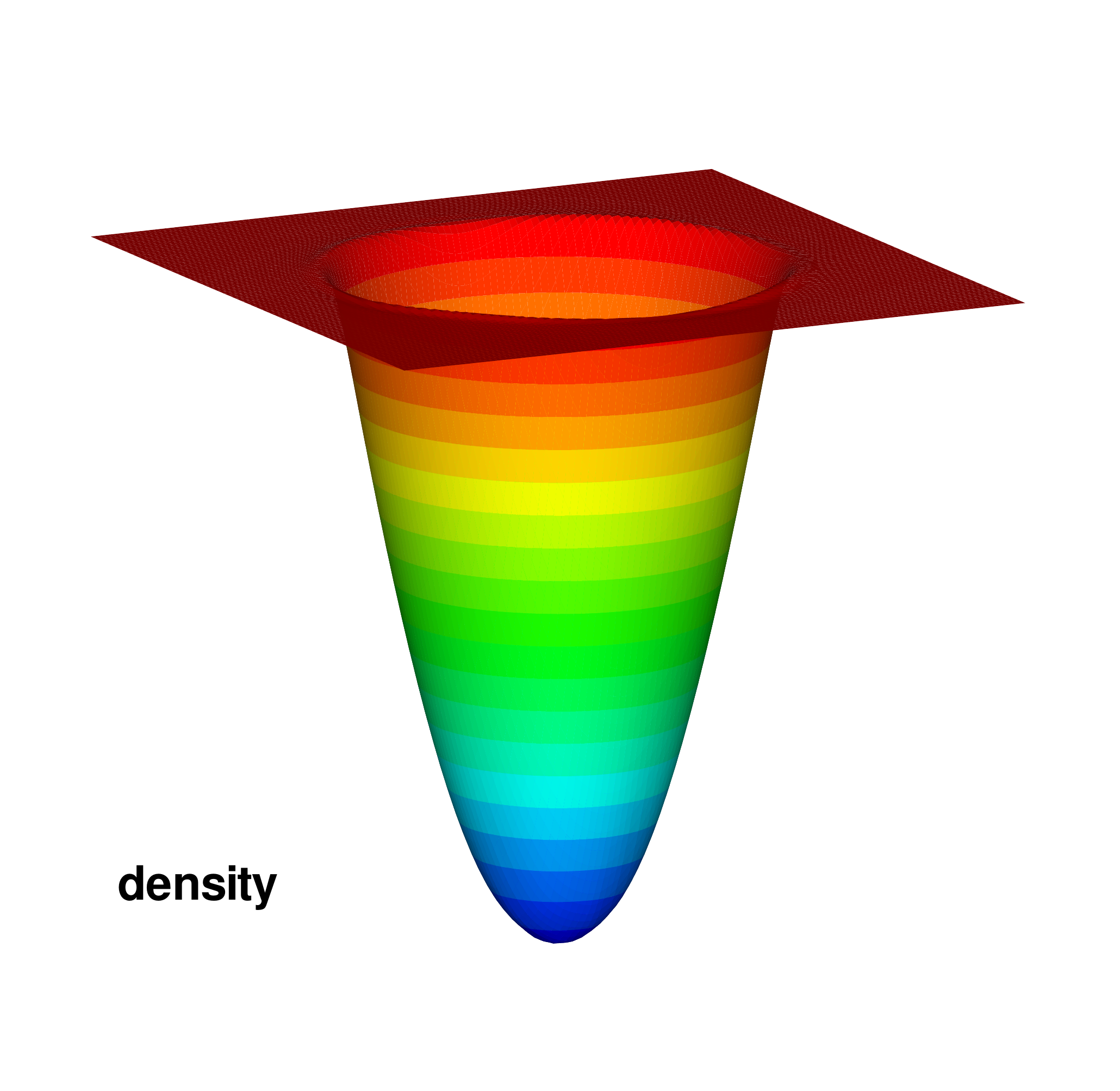}\\
\includegraphics[width=0.475\textwidth]{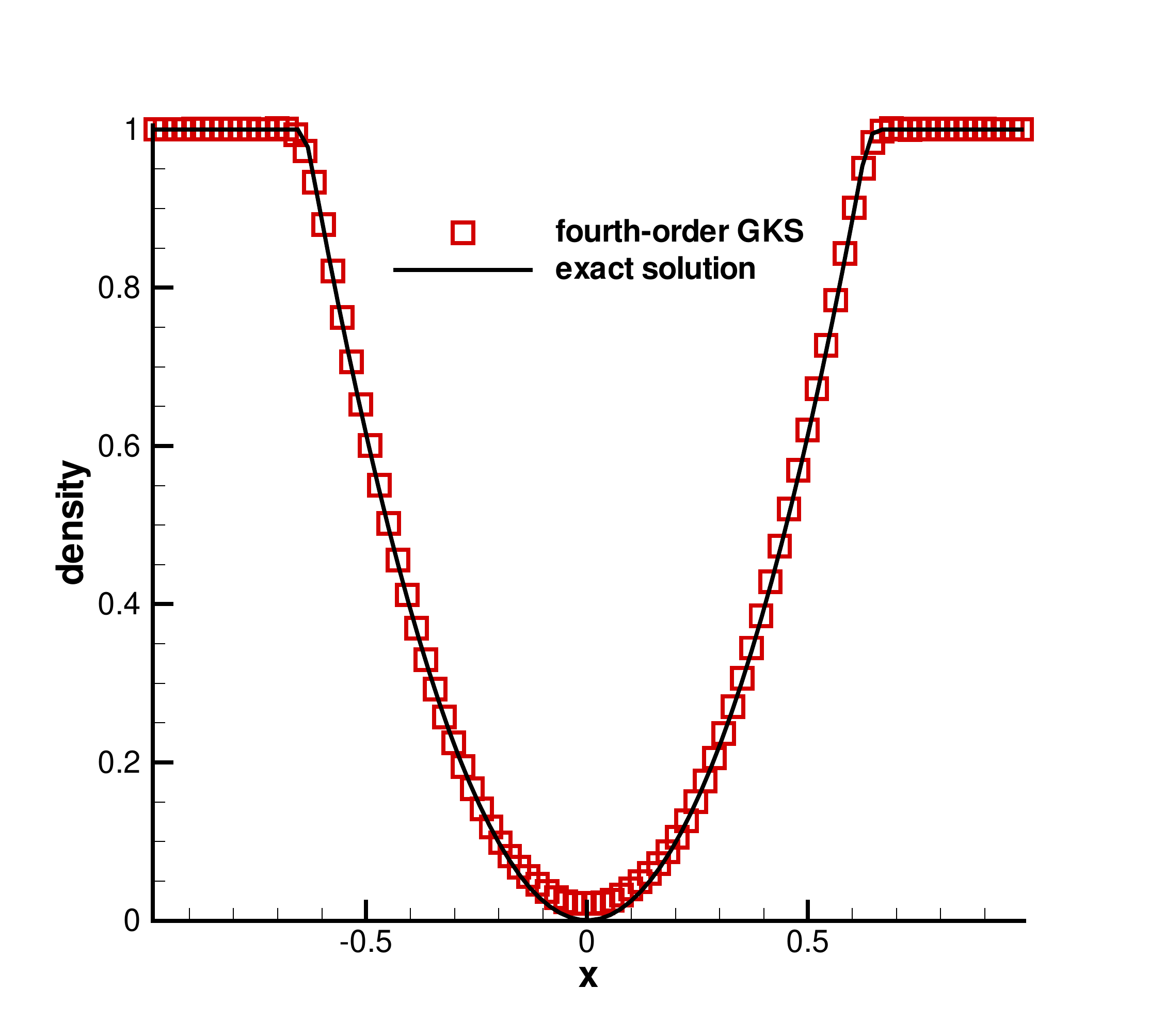}
\includegraphics[width=0.475\textwidth]{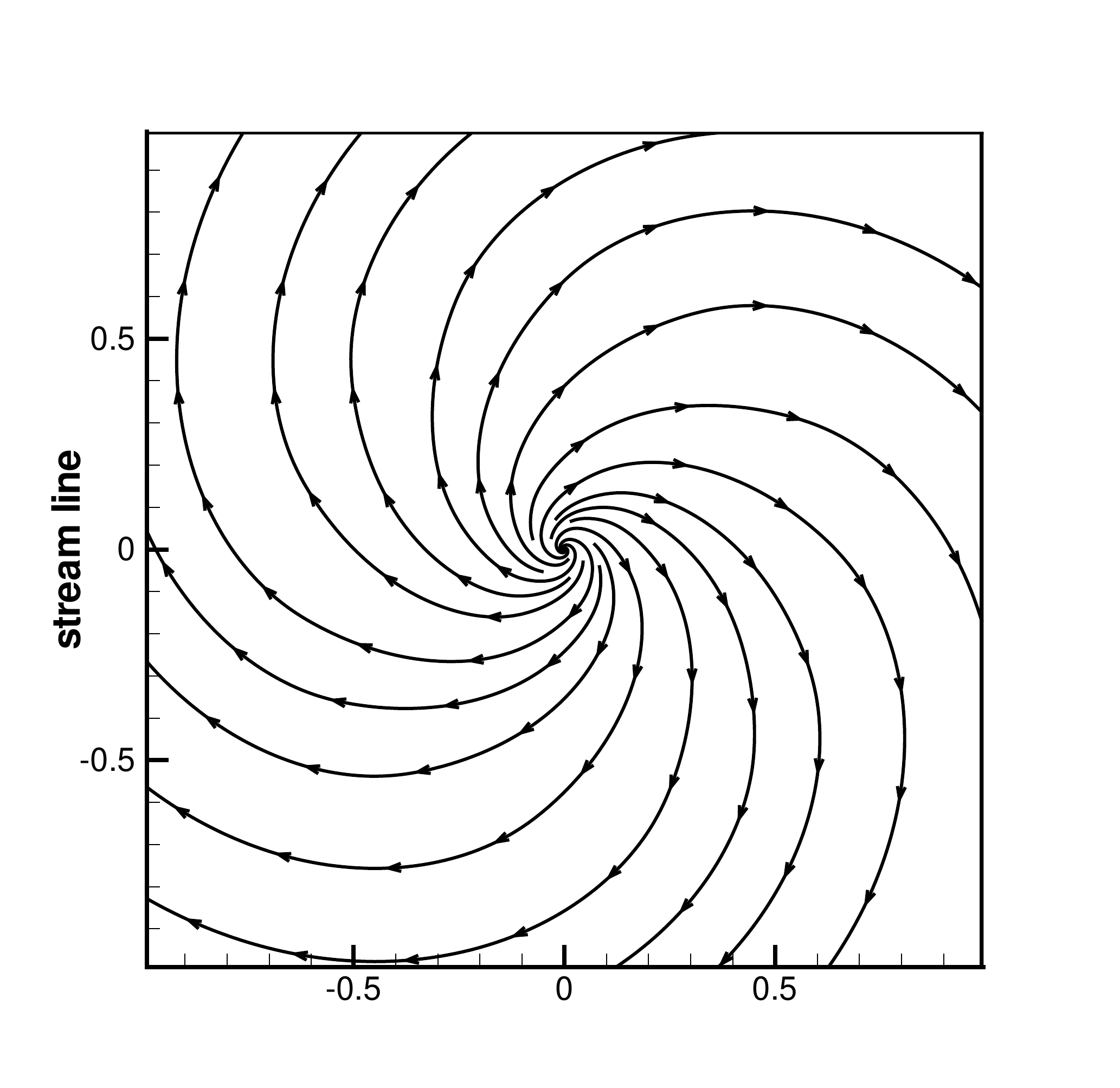}
\caption{\label{rotation-c}Hurricane-like solution: the critical
case.}
\end{figure}

\subsubsection{Hurricane-like solutions}
The first group of two-dimensional time-dependent solutions are the
hurricane-like flow field, whose solution has one-point vacuum in
the center with rotational velocity field. The initial condition are
given as
\begin{align*}
(\rho, U, V, p)=(\rho_0, v_0\sin\theta, -v_0\cos\theta,
A\rho^{\gamma}_0),
\end{align*}
where $\theta =\arctan(y/x)$, $A$ is the initial entropy,
$\gamma=2$. In this case, the initial velocity distribution has a
nontrivial transversal component, which makes the flow rotational.
The solutions are classified into three types \cite{Zhang-Zheng}
according to the initial Mach number $M_0=|v_0|/c_0$, where $c_0$ is
the sound speed.

{\em 1. Critical rotation with $M_0=\sqrt 2$.} For this case, we
have an exact solution with explicit formula. This solution consists
of two parts: a far field solution and a near-field solution. The
former far field solution is defined for $r\geq2t\sqrt{p'(\rho_0)}$,
$r=\sqrt{x^2+y^2}$,
\begin{align*}
\begin{cases}
U(x,y,t)=(2tp_0'\cos\theta+\sqrt{2p_0'}\sqrt{r^2-2t^2p_0'}\sin\theta)/r,\\
V(x,y,t)=(2tp_0'\sin\theta-\sqrt{2p_0'}\sqrt{r^2-2t^2p_0'}\cos\theta)/r,\\
\rho(x,y,t)=\rho_0,
\end{cases}
\end{align*}
and the near-field solution is defined for $r<2t\sqrt{p'(\rho_0)}$
\begin{align*}
U(x,y,t)=\frac{x+y}{2t}, ~V(x,y,t)=\frac{-x+y}{2t},
~\rho(x,y,t)=\frac{r^2}{8At^2}.
\end{align*}
The curl of the velocity in the near-field is
\begin{align*}
curl(U,V)=V_x-U_y=-\frac{1}{2t}\neq 0,
\end{align*}
and the solution has one-point vacuum at the  origin $r=0$. This is
a typical hurricane-like solution that behaves highly singular,
particularly near the origin $r=0$. \vspace{0.2cm}

There are two issues here challenging the numerical schemes: One is
the presence of the vacuum state which examines whether a high order
scheme can keep the positivity preserving property; the other is the
rotational velocity field that tests whether a numerical scheme can
preserve the symmetry. In our computation we take the data $A=25,
v_0=10, \rho_0=1$ and $\Delta x=\Delta y=1/100$. Numerical results
are presented in Fig.\ref{rotation-c} with the distribution of
density and velocity field. This case shows the robustness of the
current scheme. The positivity and symmetry are all preserved well.
\vspace{0.2cm}

\begin{figure}[!htb]
\centering
\includegraphics[width=0.475\textwidth]{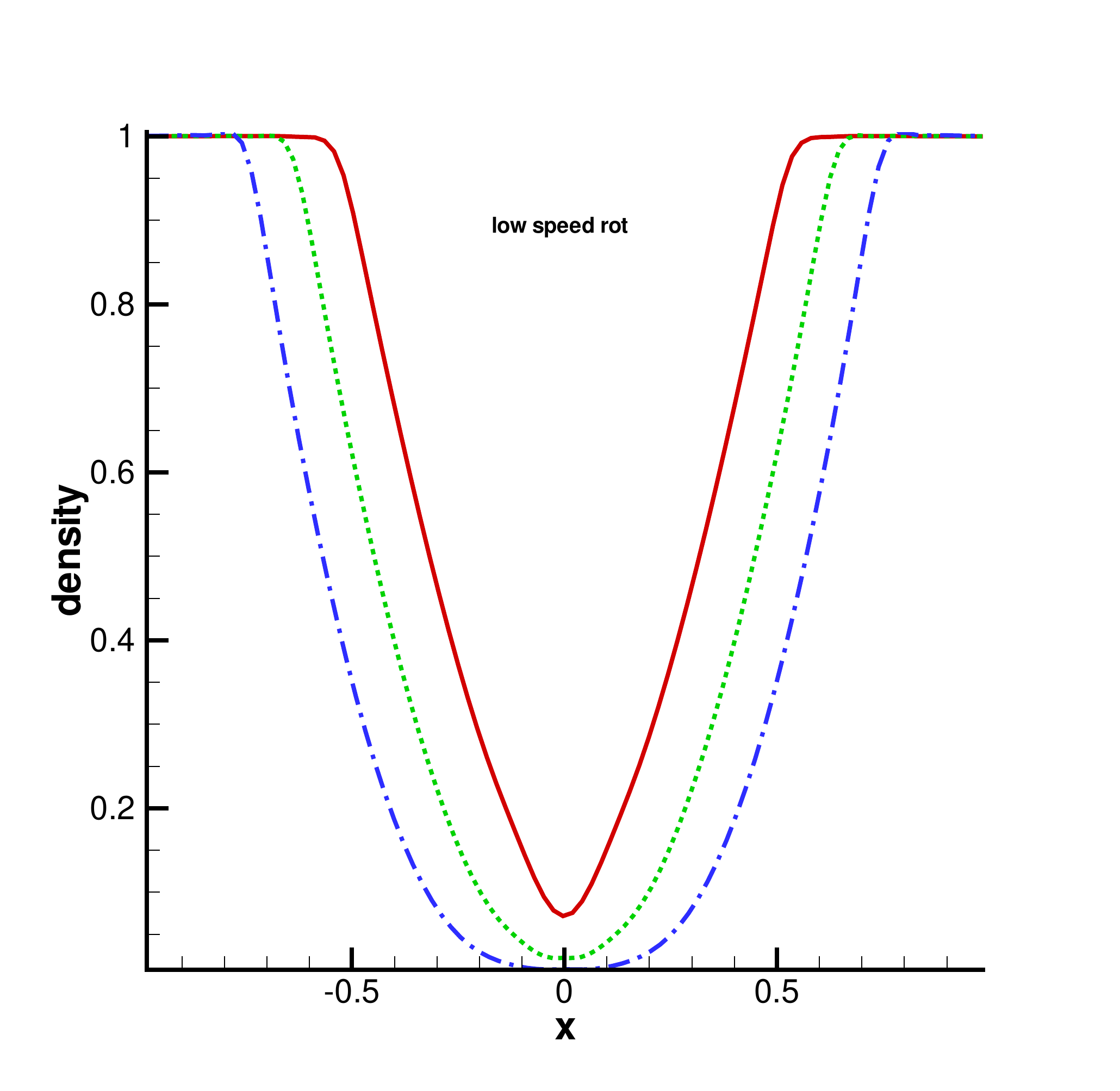}
\includegraphics[width=0.475\textwidth]{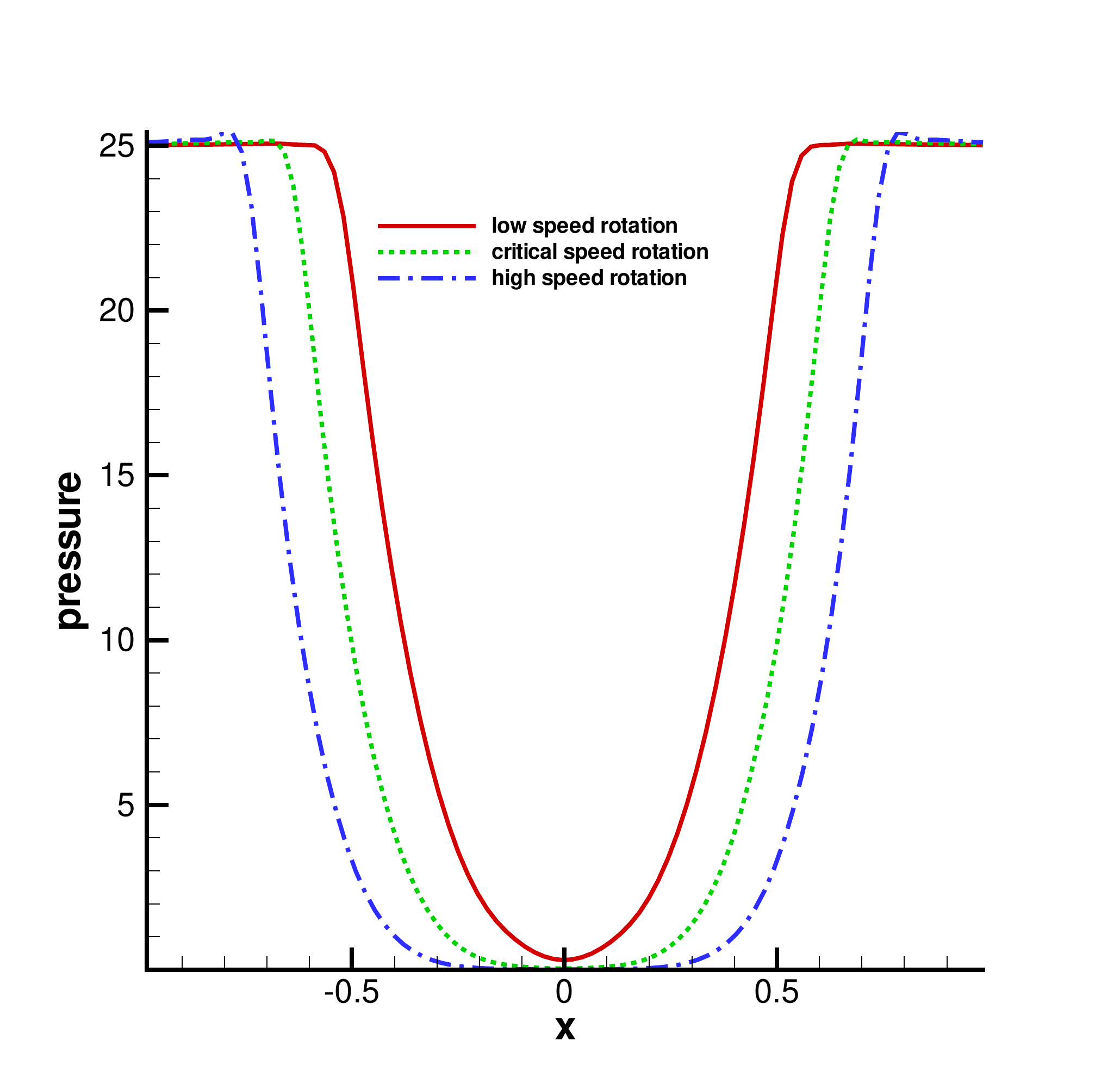}\\
\includegraphics[width=0.475\textwidth]{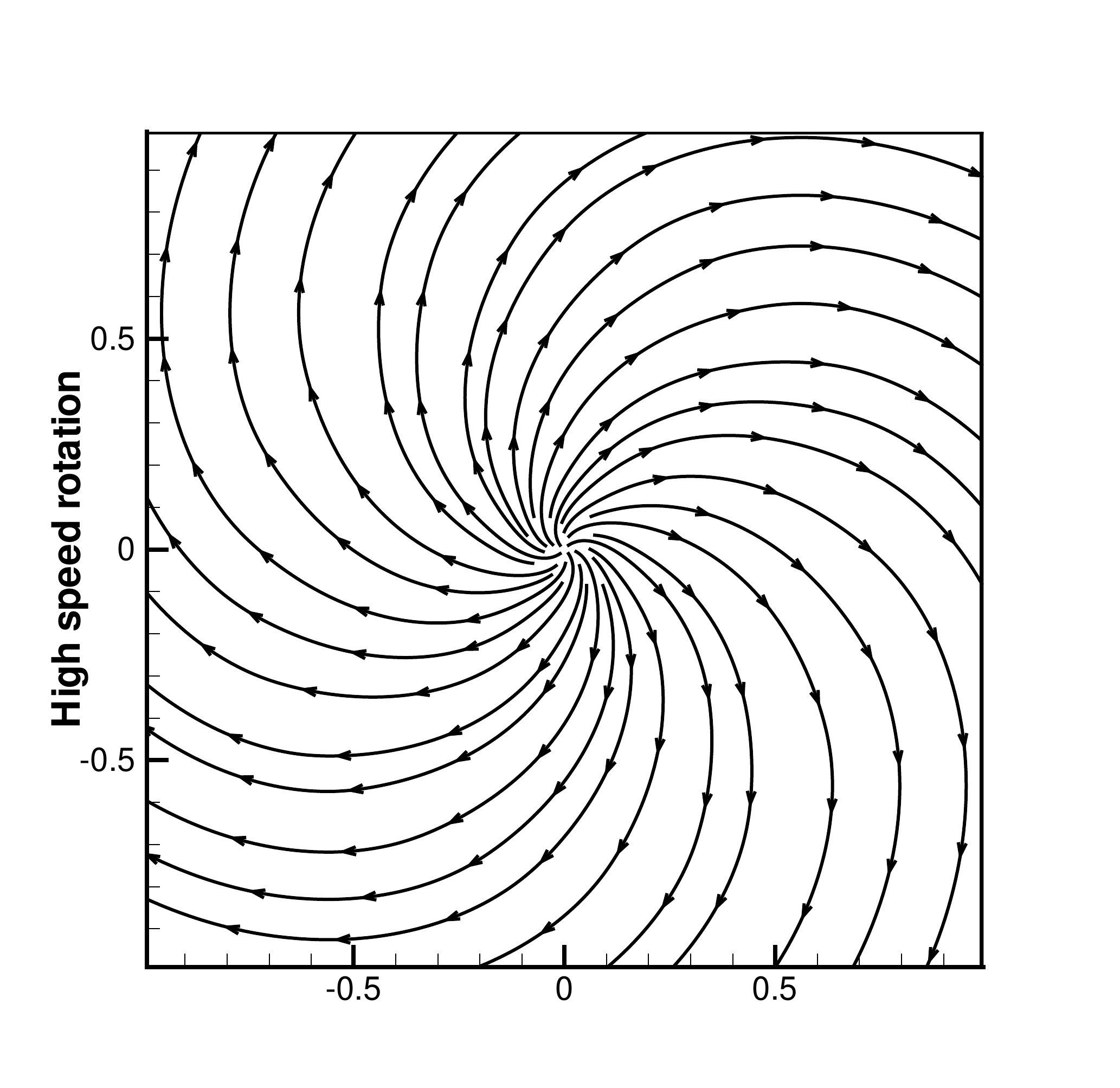}
\includegraphics[width=0.475\textwidth]{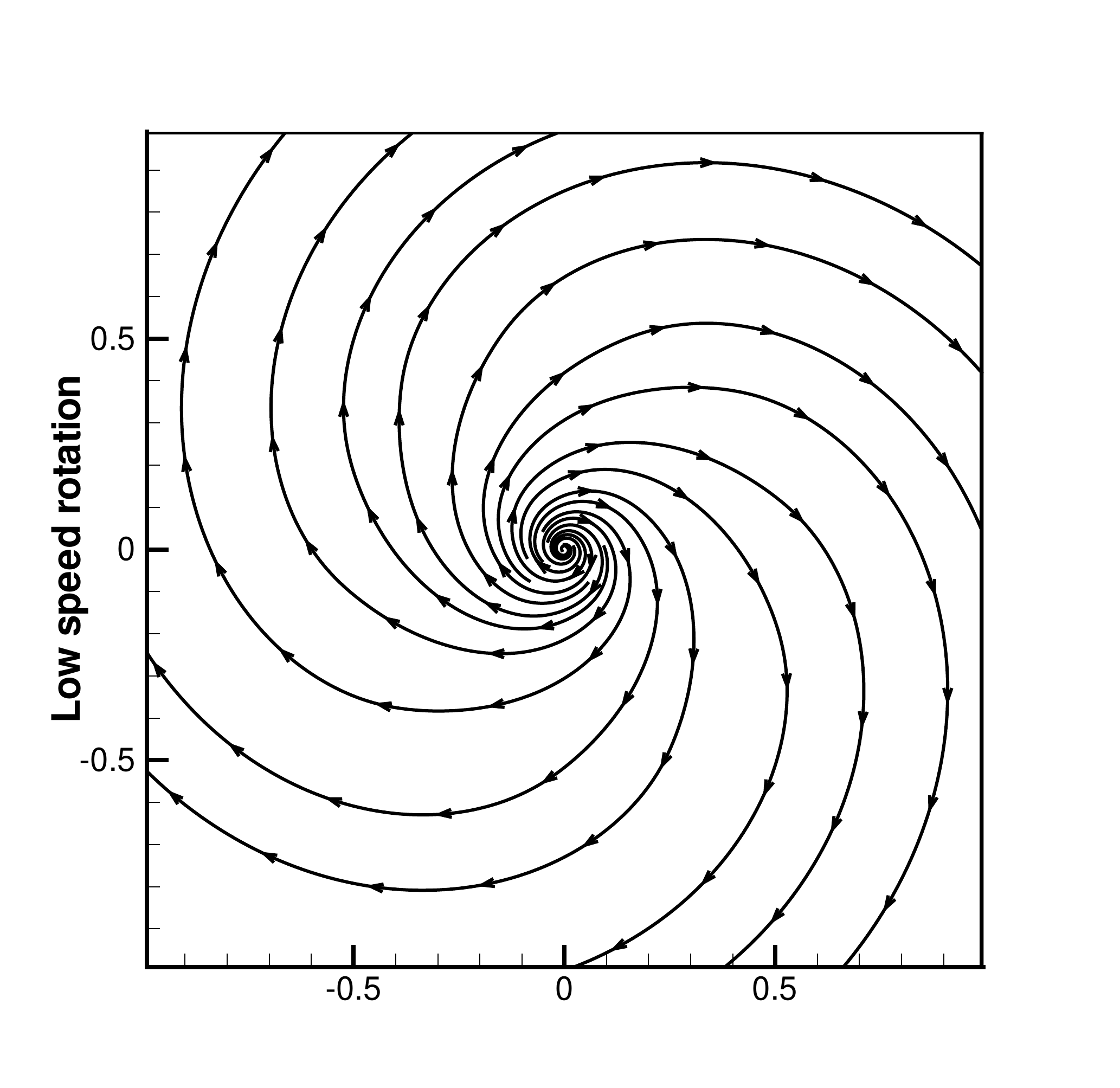}
\caption{\label{high-low-rotation}Hurricane-like solution: the cases
of high speed rotation and low speed rotation.}
\end{figure}

{\em 2. High-speed rotation with $M_0>\sqrt 2$.} For this case, the
density goes faster to the vacuum and the fluid rotates more severe.
In the computation $A=25, v_0=12.5, \rho_0=1$ and $\Delta x=\Delta
y=1/100$. The numerical result is shown in Fig.
\ref{high-low-rotation}. Because of its high rotation speed, this
case is more tough than the first one, and it also validates the
robustness of the scheme. \vspace{0.2cm}

{\em 3. Low-speed rotation with $M_0<\sqrt 2$.}  This case is milder
than the other cases above. There is no vacuum in the solution, but
it is still rotational with a low speed. In the computation $A=25,
v_0=7.5, \rho_0=1$ and $\Delta x=\Delta y=1/100$. The numerical
solution is presented in Fig. \ref{high-low-rotation}. We observe
the symmetry of the flow structure is preserved well.

\subsubsection{Large Mach number limit}
The second group is about the interactions of planar contact
discontinuities, whose solutions in the large Mach number limit are
singular and contain either vacuum or singular shocks or
delta-shocks \cite{2d-riemann3, Li-pressureless}. Two typical cases
related to the large Mach number limit  are provided in this group.
As is well-known, the compressible fluid flows become incompressible
for the low Mach number limit \cite{Majda}. On the other hand, the
large Mach number limit leads to the pressureless model
\cite{Kreiss, Li-pressureless}
\begin{align}\label{eq:pressureless}
\frac{\partial}{\partial{t}}
\begin{pmatrix}
   \rho   \\
   \rho U \\
   \rho V \\
   \rho E \\
 \end{pmatrix}+
\frac{\partial}{\partial{x}}
\begin{pmatrix}
   \rho U  \\
   \rho U^2\\
   \rho UV \\
   \rho EU \\
 \end{pmatrix}+
\frac{\partial}{\partial{y}}
\begin{pmatrix}
   \rho V  \\
   \rho UV \\
   \rho V^2\\
   \rho EV \\
 \end{pmatrix}=0.
\end{align}
This system can be also regarded as the zero moment closure of the
Boltzmann-type equations \cite{Sinai} or describe the single
transport effect of mass, momentum. The last equation of  is
decoupled from the first three equations and thus it is sufficient
to consider the solution $(\rho, U,V)$ of the first three equations
in Eq.\eqref{eq:pressureless}. \vspace{0.2cm}

The Riemann-type initial conditions for Eq.\eqref{eq:Euler} are
given as follows
\begin{align}\label{data:2driemann}
(\rho,U,V,p)&= \left\{\begin{aligned}
&(\rho_1,U_1,V_1,p_1), &x>0.5,y>0.5,\\
&(\rho_2,U_2,V_2,p_2), &x<0.5,y>0.5,\\
&(\rho_3,U_3,V_3,p_3), &x<0.5,y<0.5,\\
&(\rho_4,U_4,V_4,p_4), &x>0.5,y<0.5.
\end{aligned}
\right.
\end{align}
The negative contact discontinuity and positive contact
discontinuity, which connect the $l$ and $r$ areas, are denoted as
$J^-_{lr}$ and $J^+_{lr}$ respectively
\begin{align*}
J_{lr}^-: w_l=w_r, p_l=p_r, w_l'\geq w_r',\\
J_{lr}^+: w_l=w_r, p_l=p_r, w_l'\leq w_r'.
\end{align*}
where $w_l, w_r$ are the normal velocity and $w_l', w_r'$ are the
tangential velocity. Two types of interaction of planar contact
discontinuities are considered as follows. \vspace{0.2cm}

\begin{figure}[!htb]
\centering
\includegraphics[clip, width=0.45\textwidth]{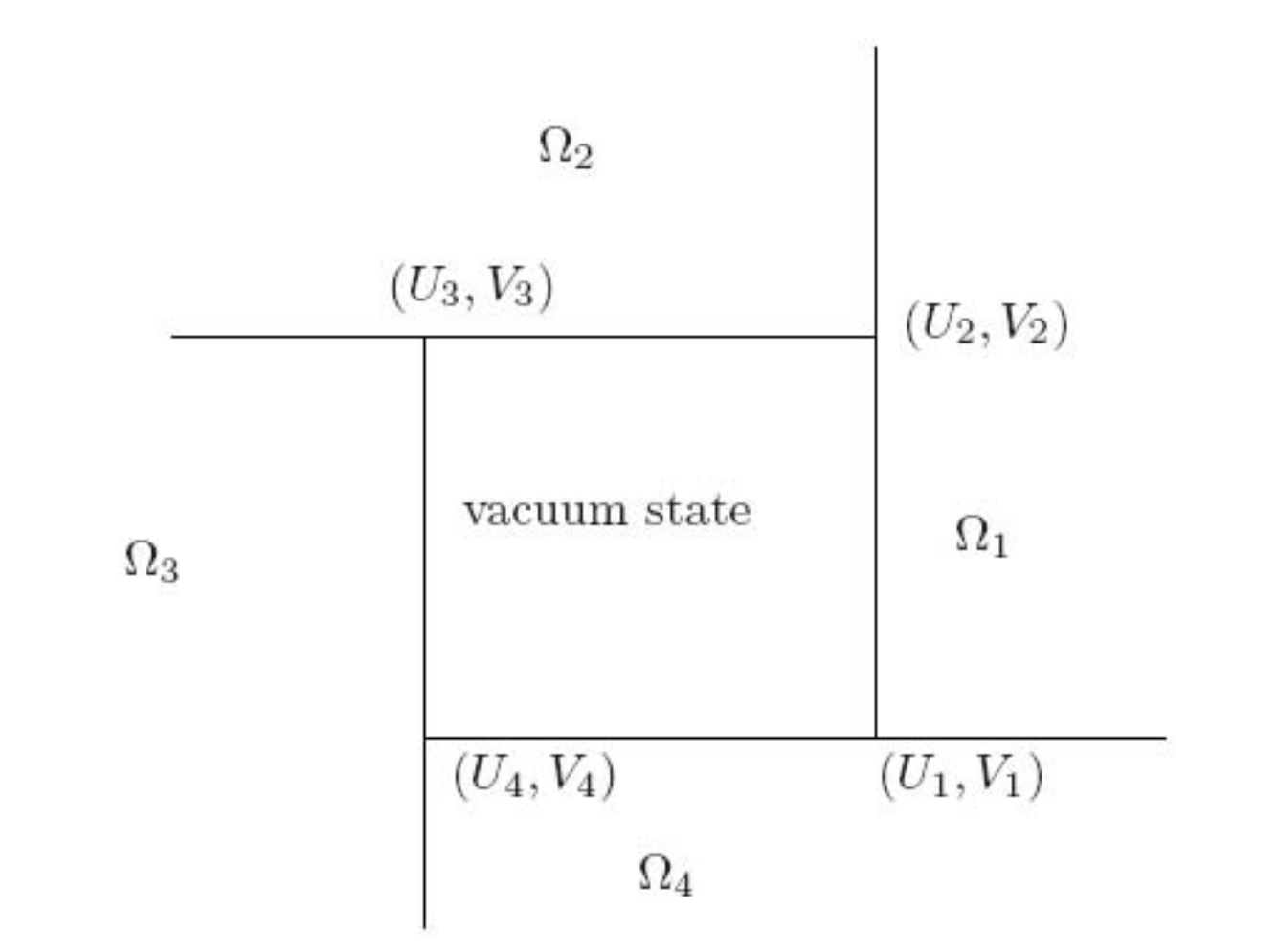}
\caption{\label{Fig:pless} The solution of the pressureless
equations Eq.\eqref{eq:pressureless} with the data
\eqref{data:same-2}. This figure is displayed in the self-similarity
$(x/t,y/t)$--plane. The notation $(U_i,V_i)$ denotes the coordinate
$(x/t,y/t)=(U_i,V_i)$, $i=1,2,3,4$. }
\end{figure}

\begin{figure}[!htb]
\centering
\includegraphics[width=0.45\textwidth]{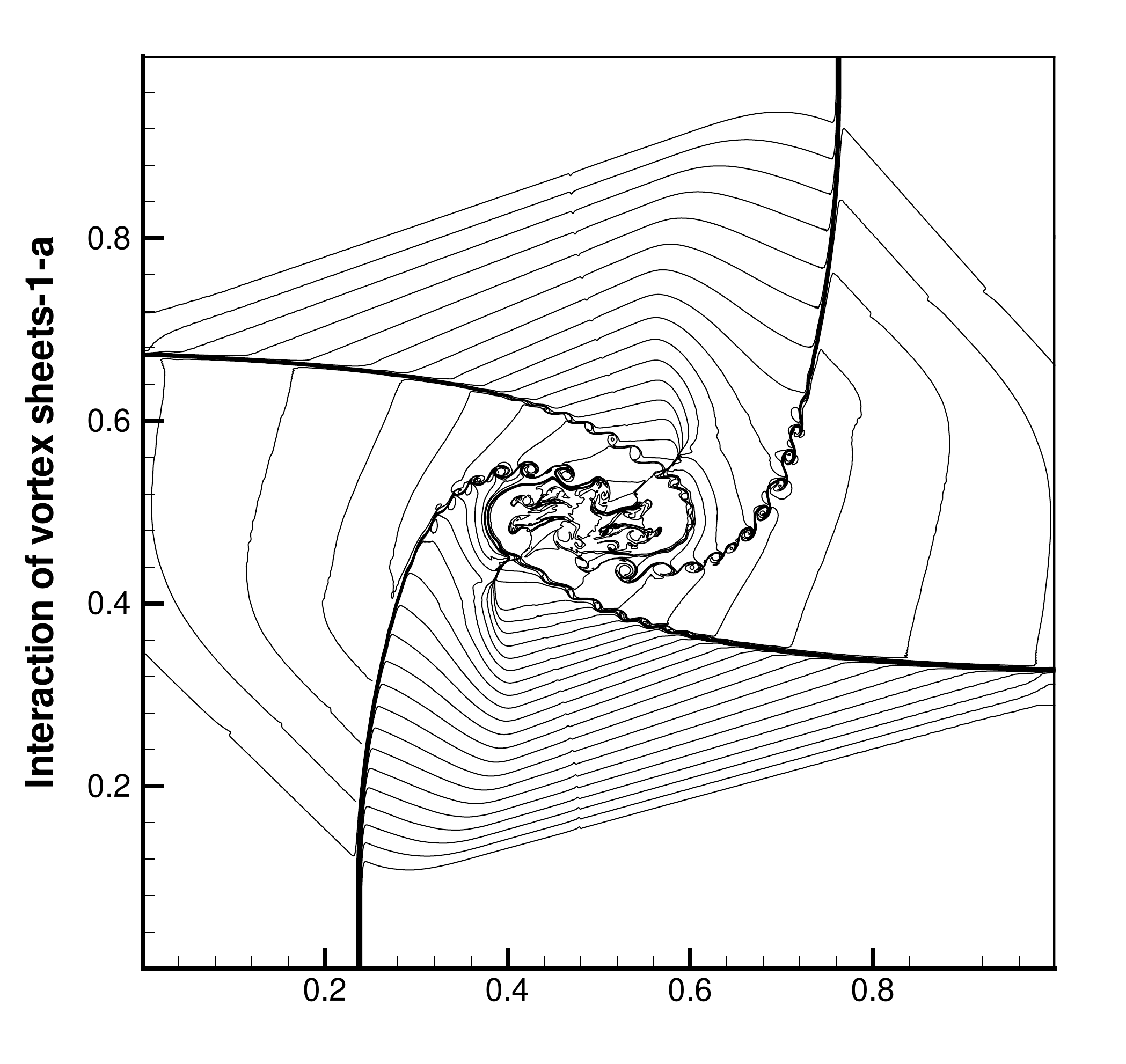}
\includegraphics[width=0.45\textwidth]{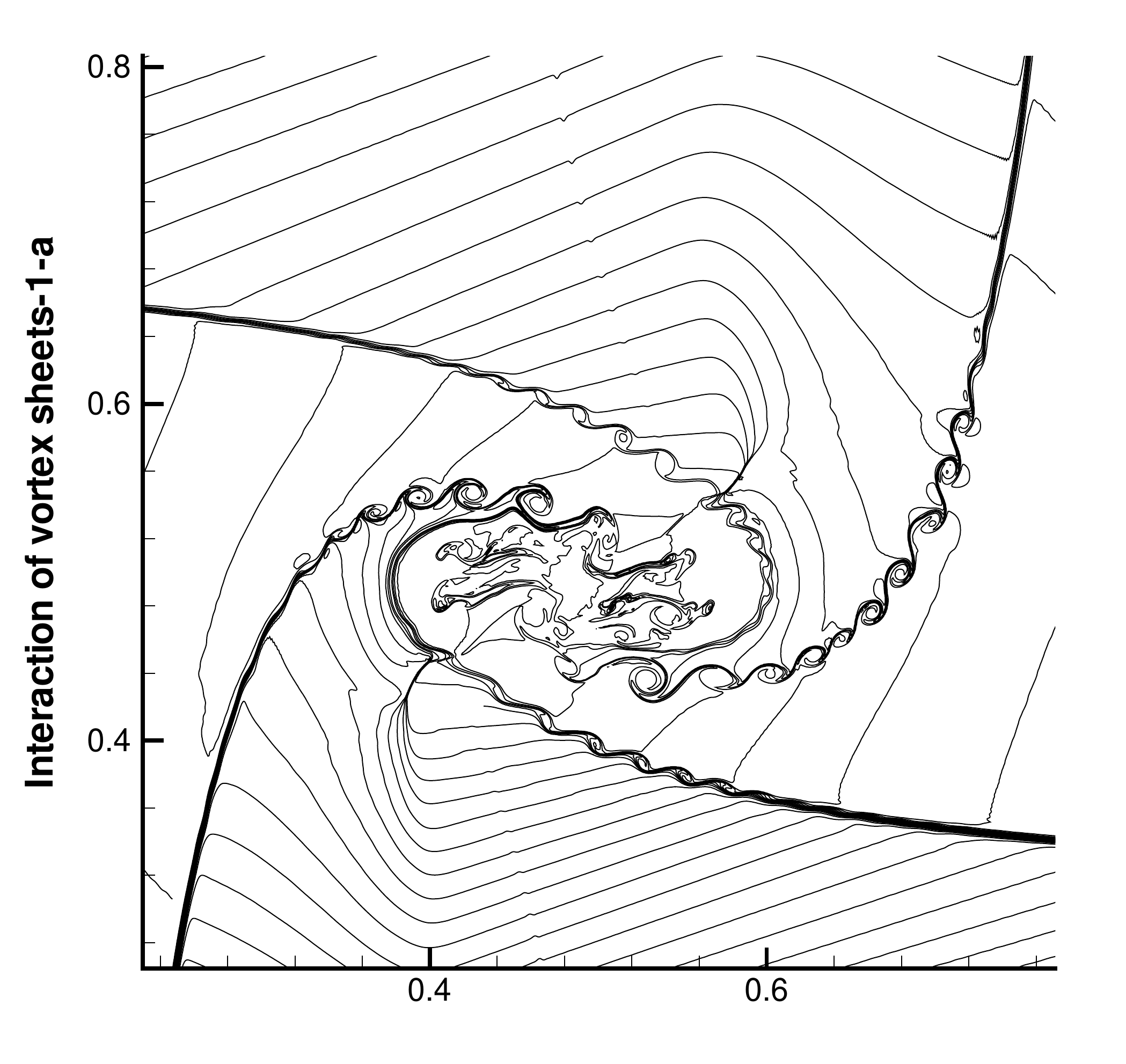}
\caption{\label{4j-1}Density distribution and the local enlargement
for the interaction of vortex sheets with same signs, where
$p_0=1$.}
\end{figure}

\begin{figure}[!htb]
\centering
\includegraphics[width=0.45\textwidth]{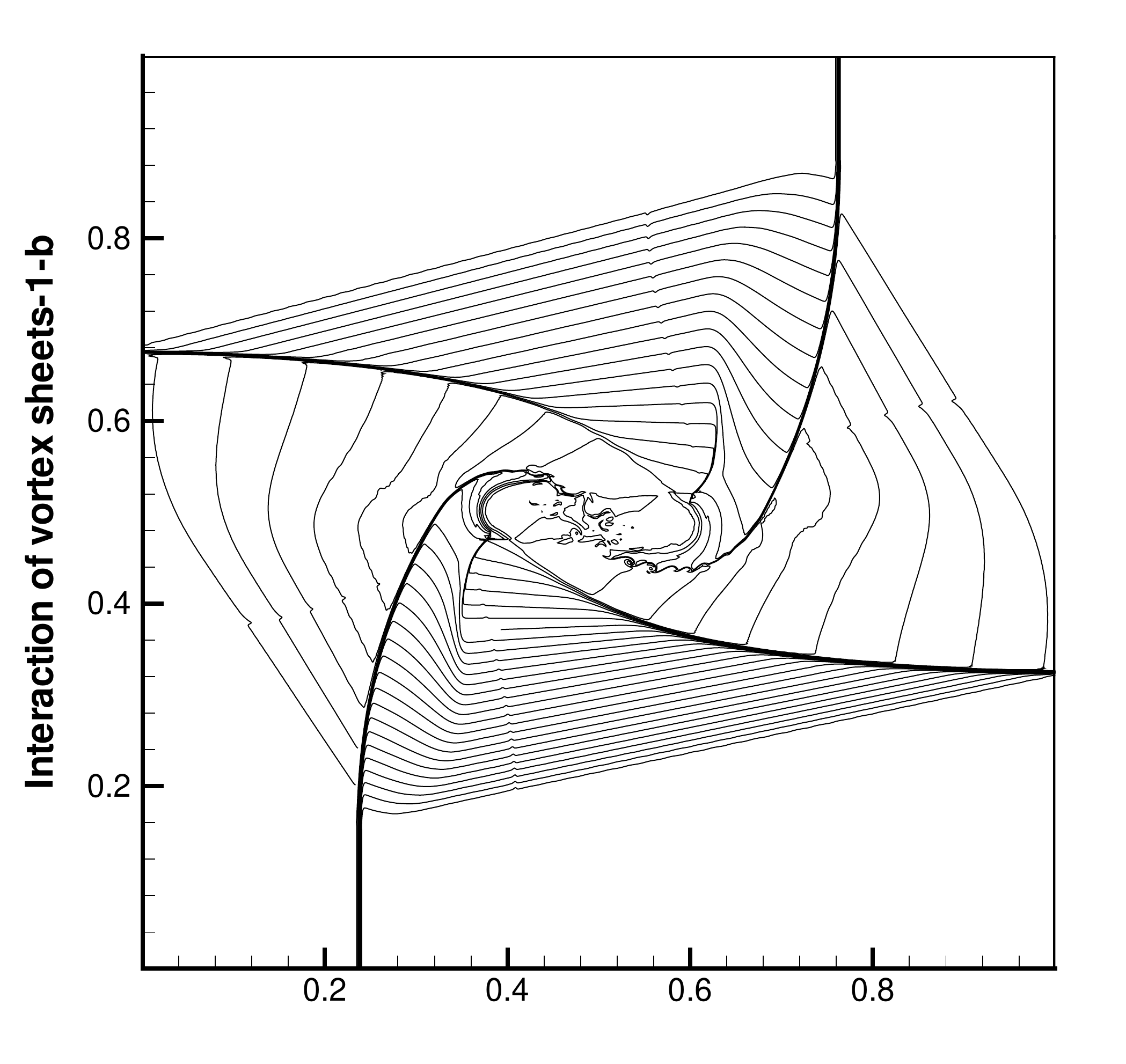}
\includegraphics[width=0.45\textwidth]{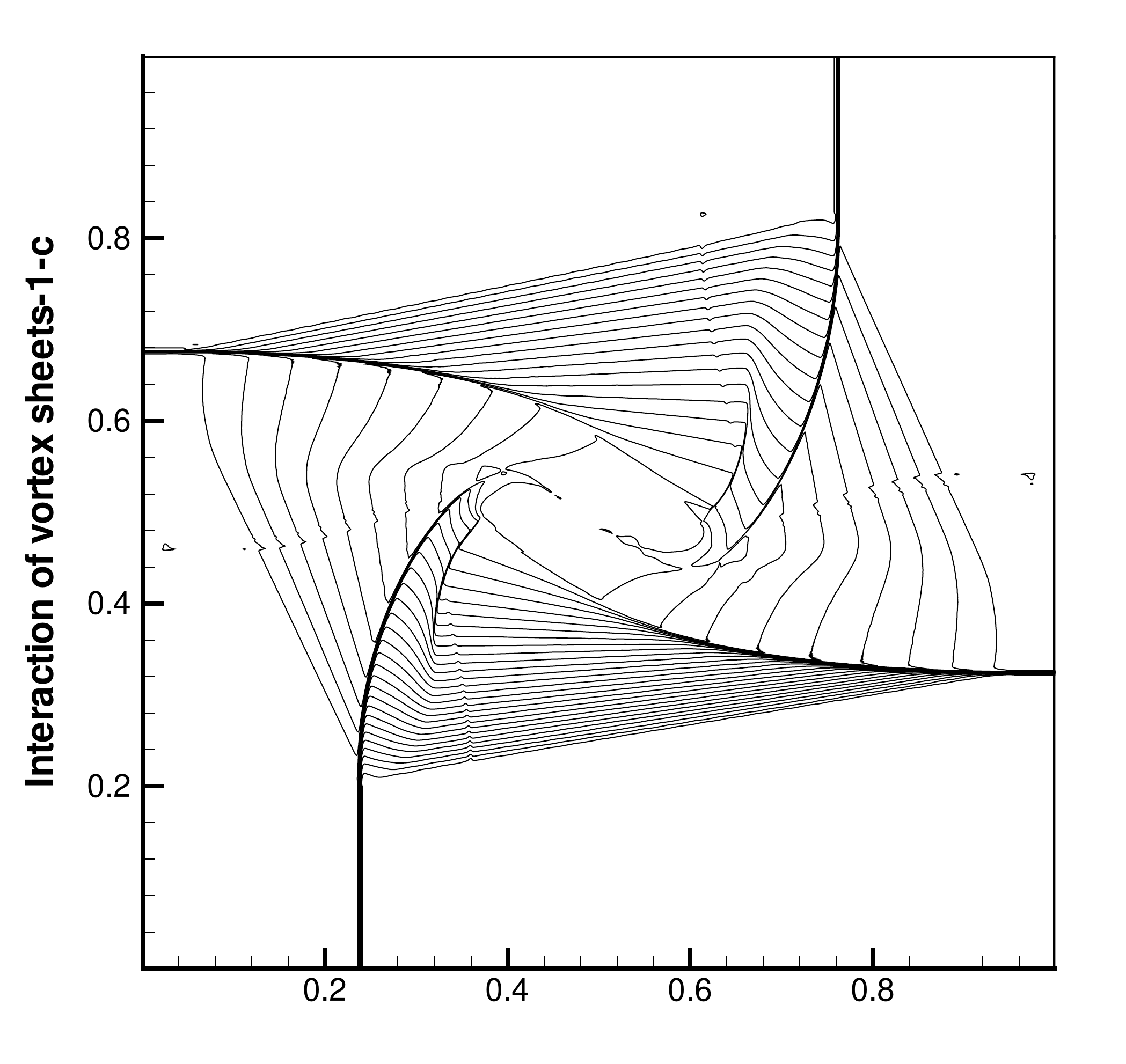}
\includegraphics[width=0.45\textwidth]{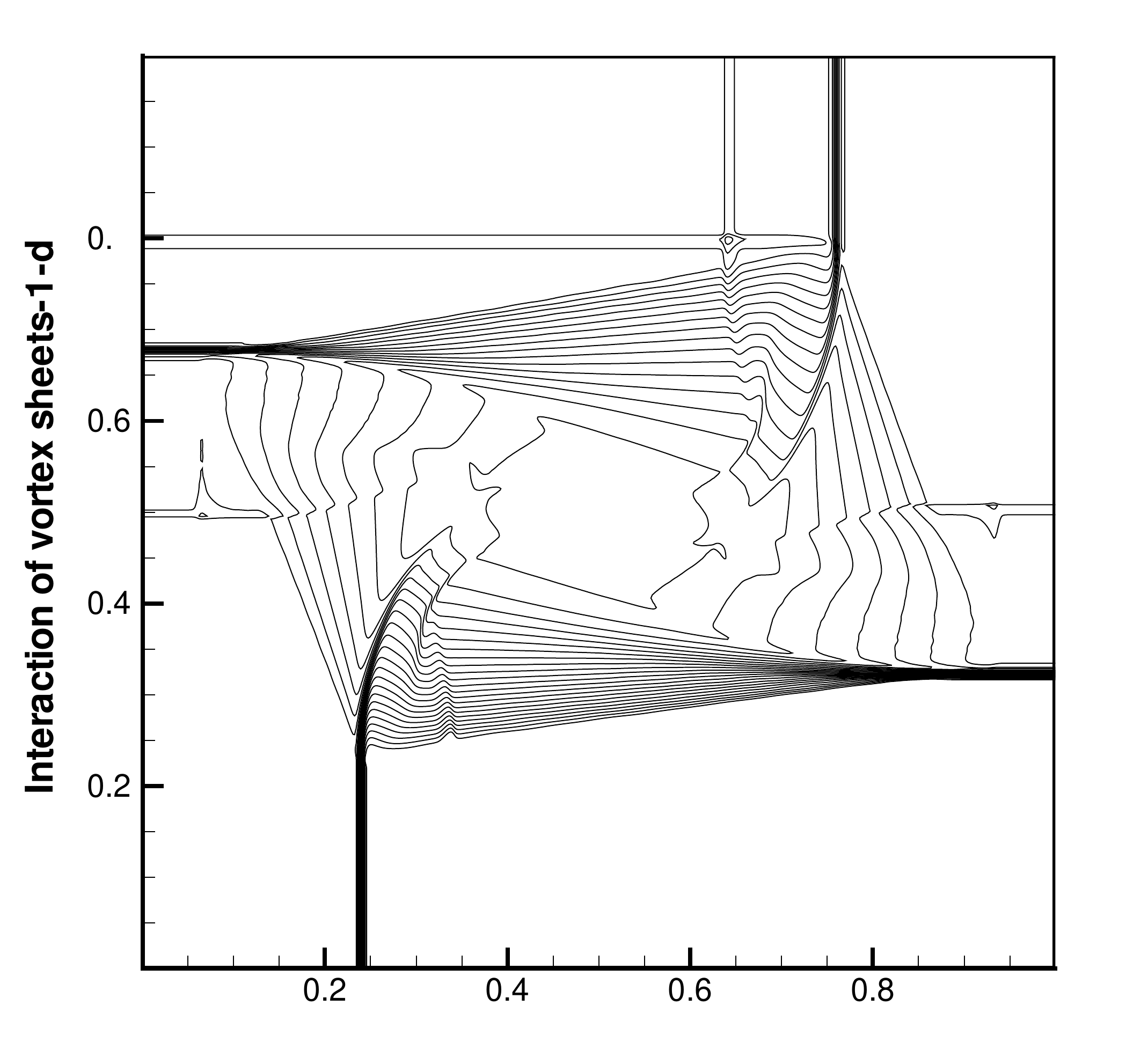}
\includegraphics[width=0.45\textwidth]{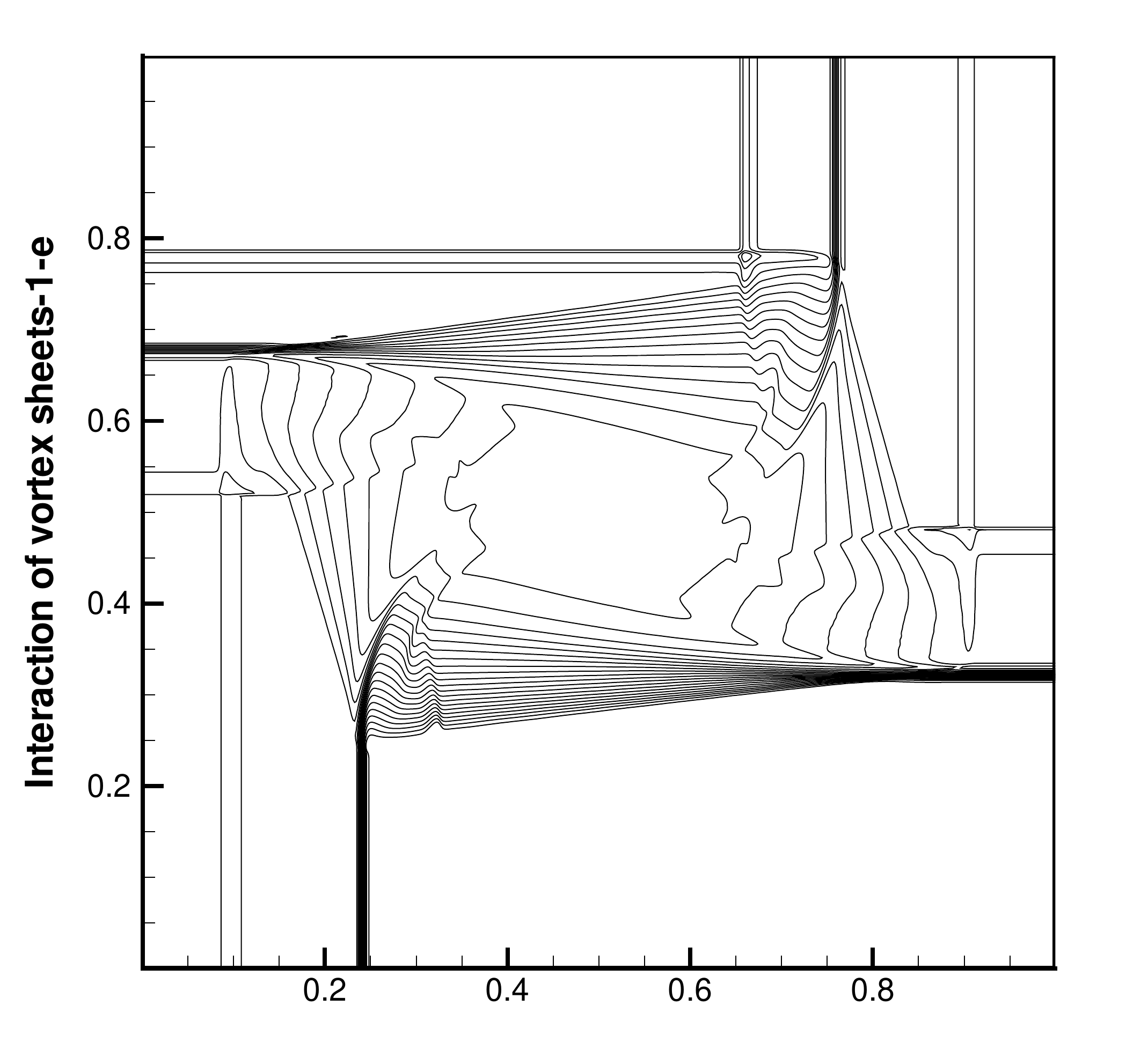}
\caption{\label{4j-2}Density distributions for the interaction of
vortex sheets with same signs, where $p_0=0.5, 0.25, 0.15$ and
$0.1$.}
\end{figure}

{\em 1. Interaction of planar contact discontinuity with same sign
vortex sheets} $J_{12}^- J_{32}^- J_{41}^- J_{34}^-$. The vacuum
solutions will be provided by the large Mach number
limit\cite{Sheng-Zhang, 2d-riemann3}. For such a case, the initial
data satisfies
\begin{equation}
U_1=U_2>U_3=U_4, \ \ \ V_2=V_3>V_1=V_4, \label{data:spiral}
\end{equation}
and in the computation, the initial  data takes
\begin{align}\label{data:same-2}
\begin{cases}
(\rho_1,U_1,V_1,p_1)=(1,-0.75,-0.5,p_0),\\
(\rho_2,U_2,V_2,p_2)=(2,-0.75,0.5,p_0),\\
(\rho_3,U_3,V_3,p_3)=(1,0.75,0.5,p_0),\\
(\rho_4,U_4,V_4,p_4)=(3,0.75,-0.5,p_0).
\end{cases}
\end{align}
This initial pressure distribution $p_0$ is uniform and the density
distribution could be arbitrary. Four planar contact discontinuities
$J_{ij}^-$ separate neighboring states and support the same sign
vortex sheets. On the contact discontinuity $J_{21}^-$, the density
undergoes a jump, and the vorticity is a singular measure
\begin{equation*}
curl(U, V)=(V_1-V_2)\delta(x-U_1t, y,t), V_1-V_2<0,
\end{equation*}
where $\delta(x,y,t)$ is the standard Dirac function (measure).
Therefore this contact discontinuity is the composite of an entropy
wave and a vortex sheet. The same applies for $J_{23}^-$, $J_{34}^-$
and $J_{41}^-$. Their instantaneous interaction results in a complex
wave pattern. The limiting system Eq.\eqref{eq:pressureless} has an
explicit formula that consists of four constant states in $\Omega_i,
i=1,2,3,4$ and a vacuum state inside a pyramid with edges
$(x/t,y/t)=(U_i,V_i)$, i.e.
\begin{equation*}
(\rho,U,V,p)(x,y,t) =\left\{
\begin{array}{ll}(\rho_i,U_i,V_i,p_0), \ \ \ \ & \mbox{$(x,y,t) \in \Omega_i$,}\\
\mbox{vacuum}, & \mbox{$(x,y,t)$ in the pyramid.}
\end{array}
\right.
\end{equation*}
The solution is schematically described in Fig.\ref{Fig:pless}.
\vspace{0.2cm}

We use the fourth order gas-kinetic scheme to look into the
asymptotic process for this problem by letting the pressure to
become smaller and smaller. The initial pressure is taken to be
$p_0=1, 0.5, 0.25, 0.15$ and $0.1$, respectively. The density
distributions are displayed in Fig.\ref{4j-1} and Fig.\ref{4j-2} at
$t=0.35$. The uniform mesh with $\Delta x=\Delta y=1/1500$ are used
for $p_0=1, 0.5$ and $0.25$. It is observed that as $p_0=1$ (the
Mach number $M_0$ is relatively large), the numerical solution
displays more small scale structures. With the increase of the
initial Mach number $M_0$, the complicated flow structure disappears
and the uniform mesh with $\Delta x=\Delta y=1/400$ are used for
$p_0=0.15, 0.1$. The solution becomes much closer to that given in
Fig.\ref{Fig:pless}. \vspace{0.2cm}

\begin{figure}[!htb]
\centering
\includegraphics[clip,width=8cm]{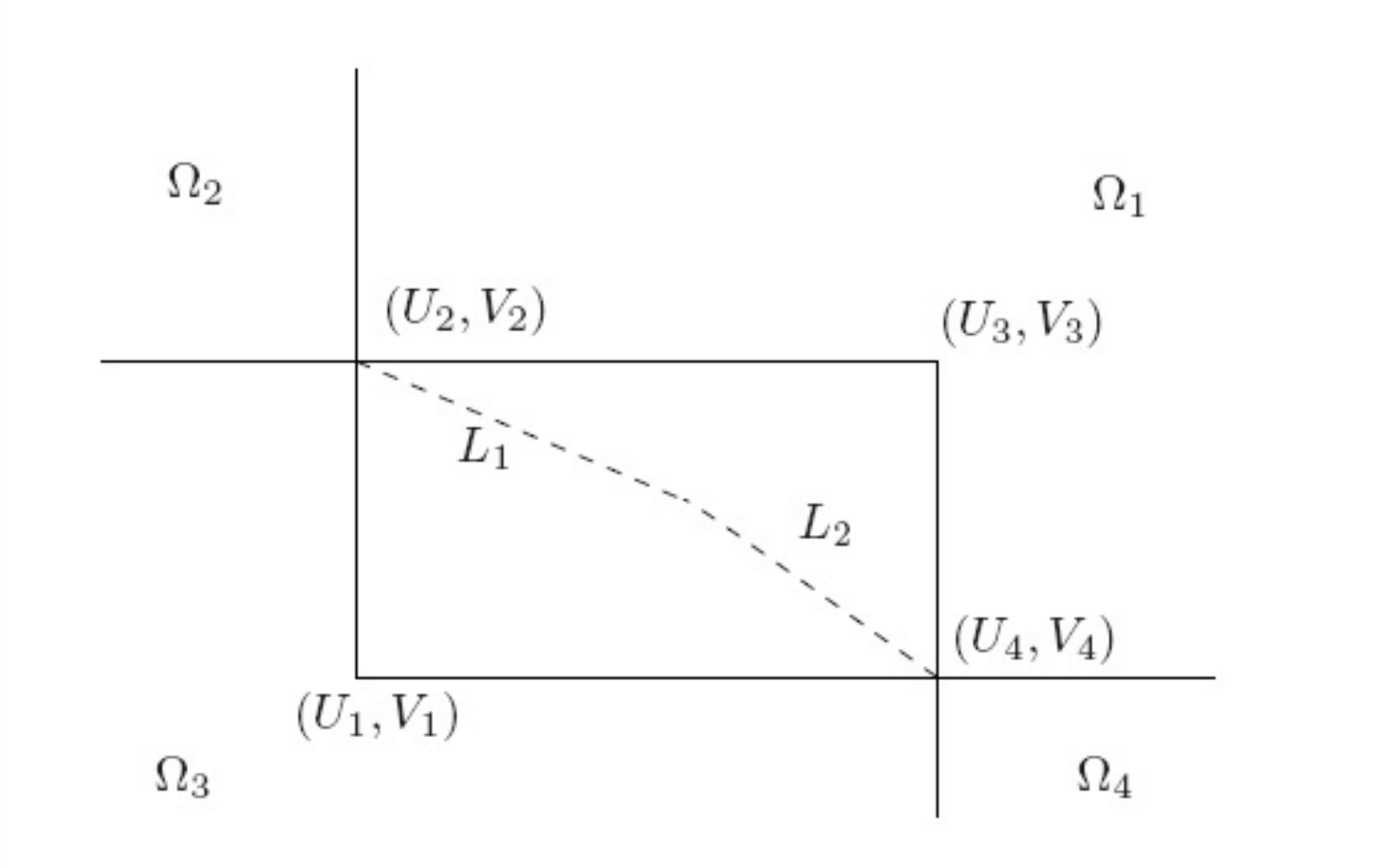}
\caption{\label{Fig:delta-shock}   The delta-shock solution of the
pressureless equations due to the interaction of vortex sheets. This
figure is displayed in the self-similarity $(x/t,y/t)$--plane. The
notation $(U_i,V_i)$ denotes the coordinate $(x/t,y/t)=(U_i,V_i)$,
$i=1,2,3,4$. The dashed line denotes the support of delta-shocks.}
\end{figure}

{\em 2. Interaction of planar contact discontinuities with vortex
sheets with different signs} $J_{21}^-J_{32}^+J_{41}^+J_{34}^-$. For
such a case,  the delta-shock solutions will be provided by the
large Mach number limit \cite{Sheng-Zhang, 2d-riemann3}. The initial
data is so designed as to satisfy
\begin{equation*}
U_3=U_4>U_1=U_2, \ \ \ V_2=V_3> V_1=V_4.
\end{equation*}
The initial pressure distribution is also uniform and the density
distribution is arbitrary. This case is different from the first
case in this group. The four planar contact discontinuities support
vortex sheets of different signs and their interaction produces
totally different flow patterns. With such initial data, the
solution of Eq.\eqref{eq:pressureless} has a singular solution
containing so-called delta-shocks, as shown in
Fig.\ref{Fig:delta-shock}. The solution formula is
\begin{equation}\label{sol:measure1} (\rho, U,V)(x,y,t)=(\rho_i,
U_i, V_i), \mbox{$(x,y,t)\in\Omega_i$.}
\end{equation}
However, the solution becomes singular, particularly, the density
takes a singular measure on the support $L_1\cup L_2$ in Fig.
\ref{Fig:delta-shock}
\begin{equation}\label{sol:measure2}
\rho(x,y,t) =\sqrt{\rho_1\rho_3} \delta(x-x(t,s), y-y(t,s),t),
\end{equation}
where $(x,y) =(x(t,s),y(t,s))$ represents the support $L_1\cup L_2$
of the Dirac measure in Fig.\ref{Fig:delta-shock}. \vspace{0.2cm}

\begin{figure}[!htb]
\centering
\includegraphics[width=0.45\textwidth]{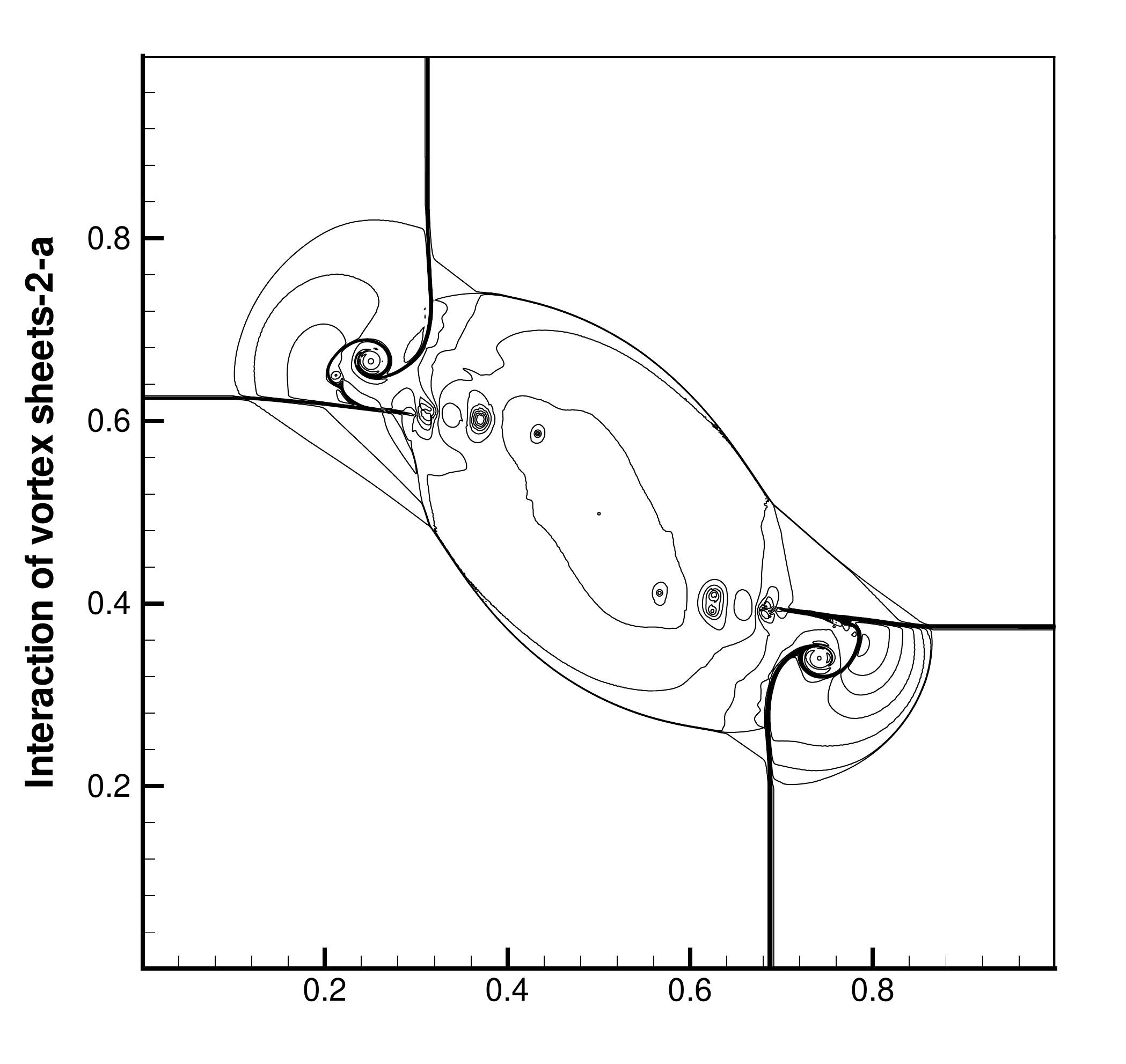}
\includegraphics[width=0.45\textwidth]{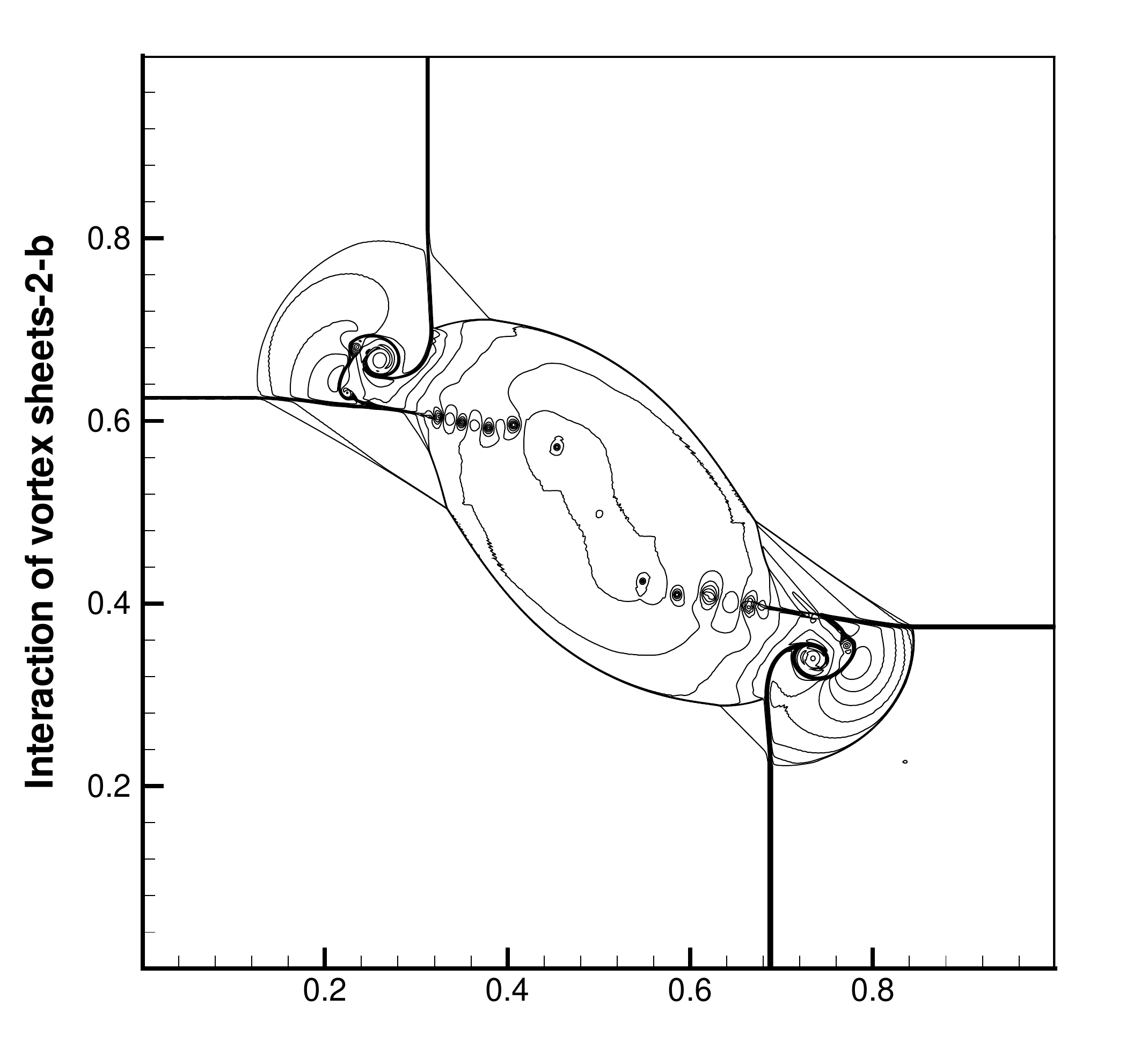}
\includegraphics[width=0.45\textwidth]{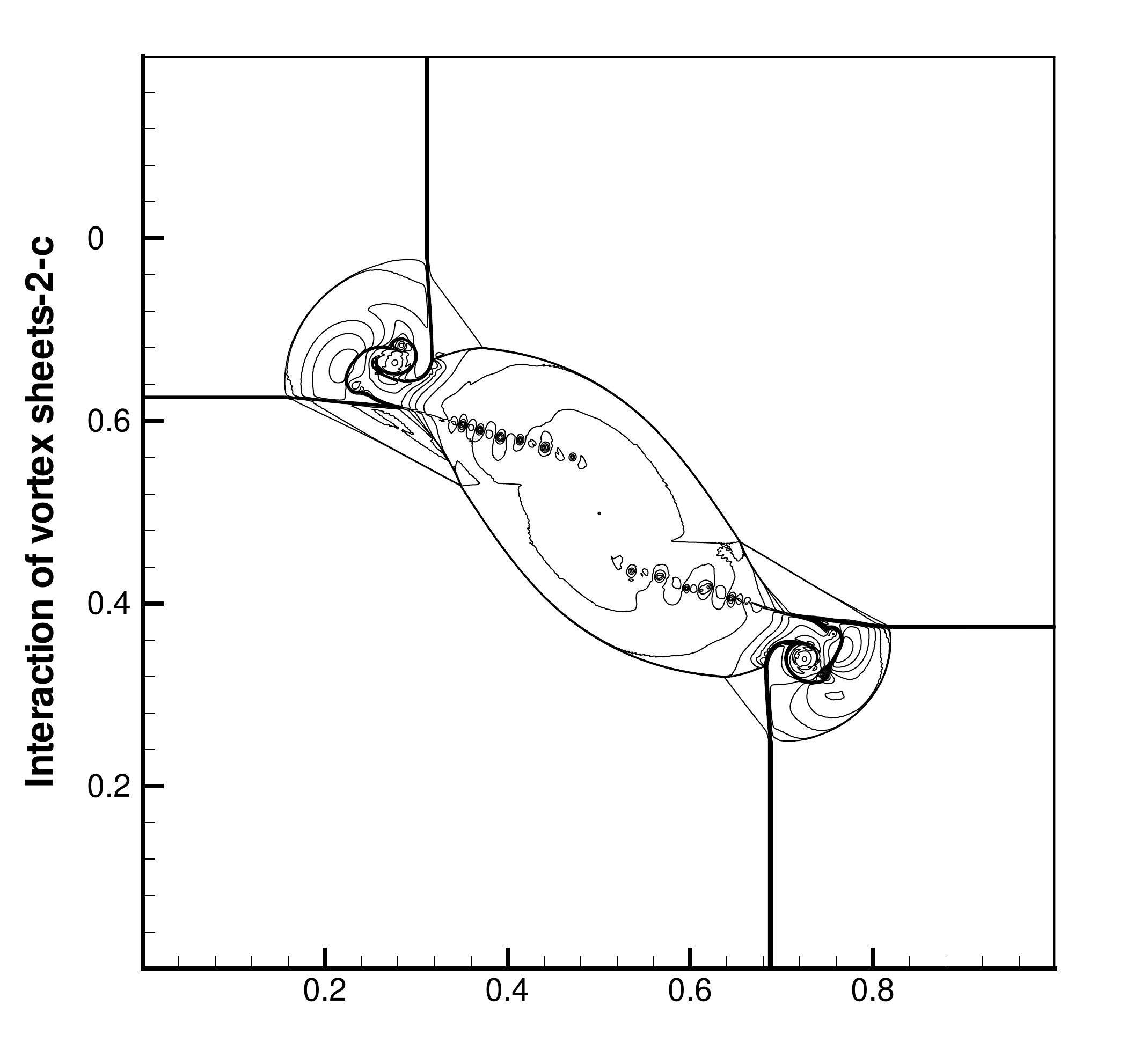}
\includegraphics[width=0.45\textwidth]{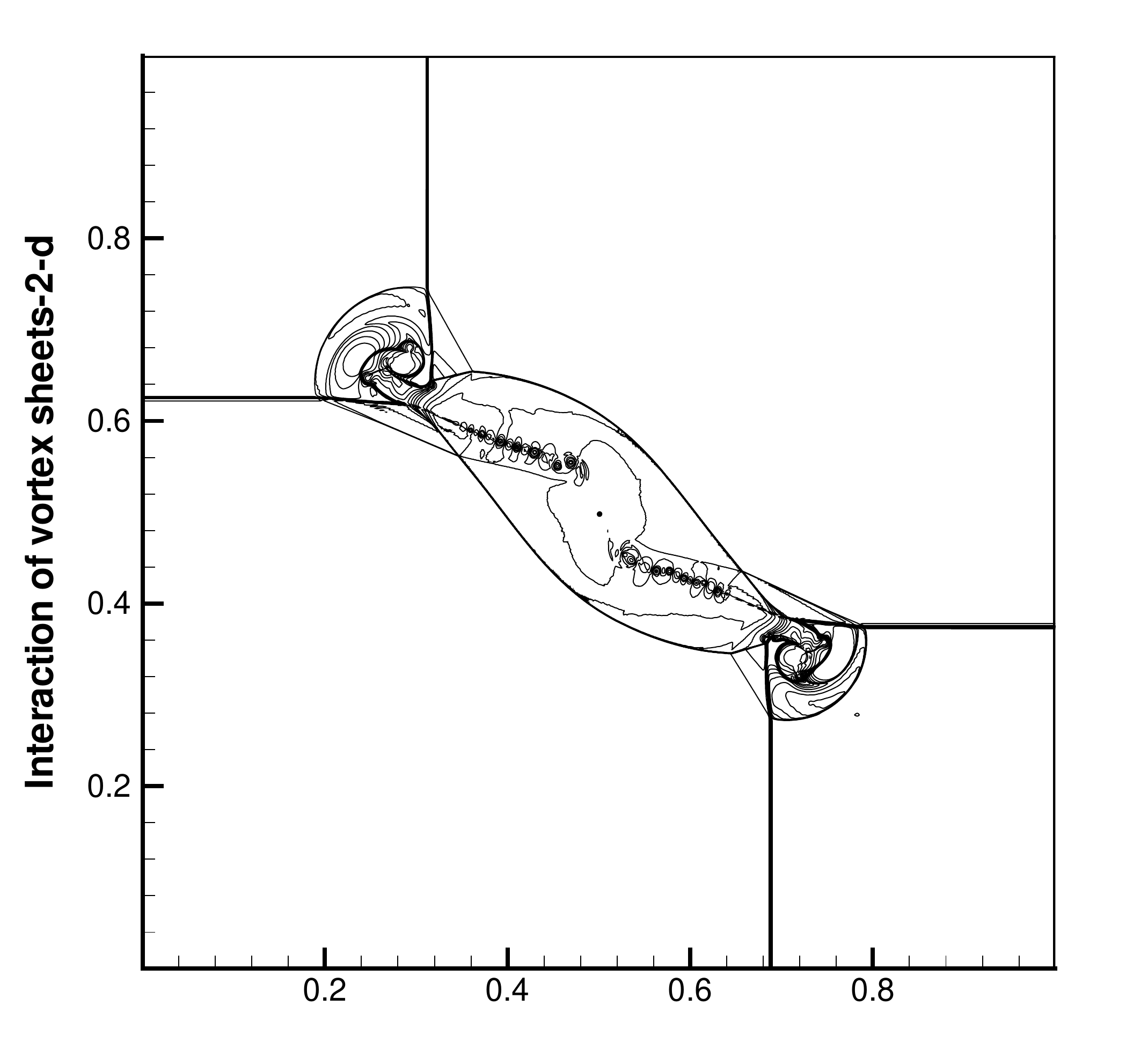}
\caption{\label{4j-3} Density distributions for the interaction of
vortex sheets with different signs, with $p_0=1, 0.75, 0.5$ and
$0.3$.}
\end{figure}

\begin{figure}[!htb]
\centering
\includegraphics[width=0.45\textwidth]{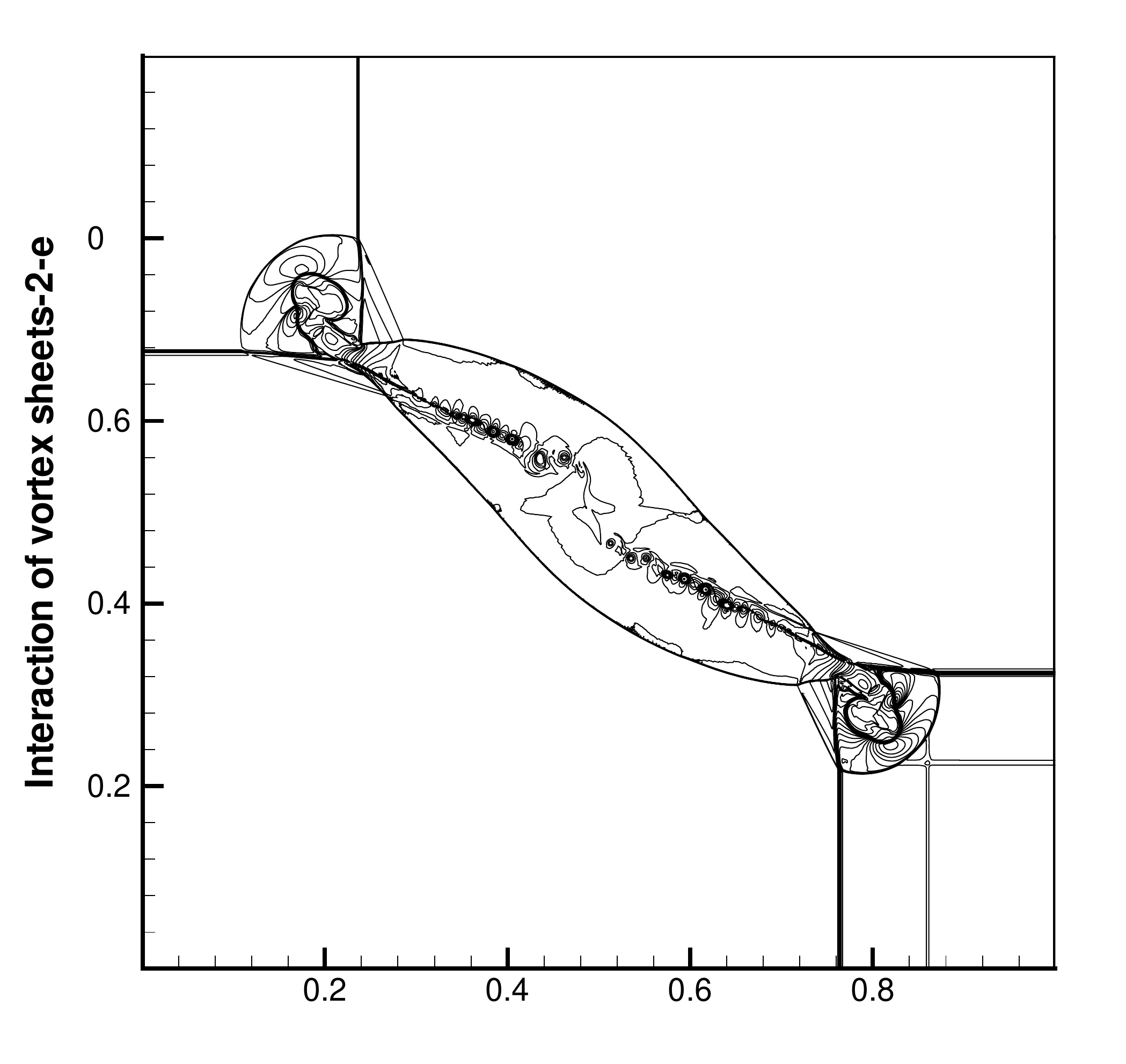}
\includegraphics[width=0.45\textwidth]{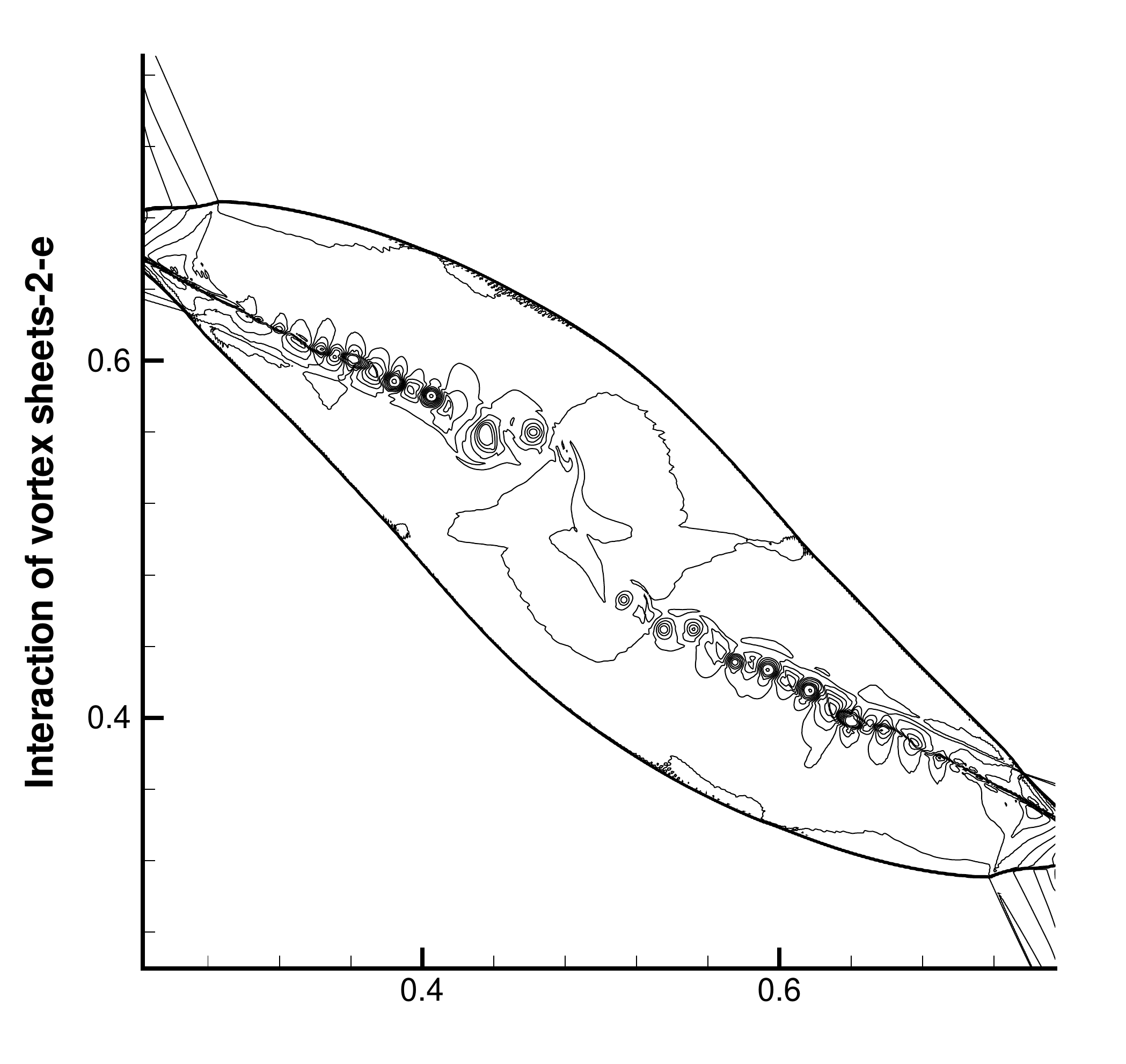}
\caption{\label{4j-4} Density distribution and the local enlargement
for the interaction of vortex sheets with different signs, with
$p_0=0.2$}
\end{figure}

The numerical simulations are designed for the cases of different
Mach numbers. The initial condition is given as follows
\begin{align*}
\begin{cases}
(\rho_1,U_1,V_1,p_1)=(1 ,0.75,-0.5,p_0),\\
(\rho_2,U_2,V_2,p_2)=(2,0.75,0.5,p_0),\\
(\rho_3,U_3,V_3,p_3)=(1,-0.75,0.5,p_0),\\
(\rho_4,U_4,V_4,p_4)=(3,-0.75,-0.5,p_0),
\end{cases}
\end{align*}
The uniform mesh with $\Delta x=\Delta y=1/1500$ is used. The
density distributions are presented in Fig.\ref{4j-3} with the
initial pressure taken as $p_0=1, 0.75, 0.5$ and $0.3$ at $t=0.25$,
respectively. With the decrease of the pressure, or equivalently
with the increase of initial Mach number, it is observed that more
and more small scaled vortices are present in the solutions. For the
case with $p_0=0.2$, the density distribution is presented in
Fig.\ref{4j-4} at $t=0.28$. The solution tends to be very close to
the solution Eq.\eqref{sol:measure1} and \eqref{sol:measure2}. This
shows the capability of the current fourth order scheme to preserve
the asymptotical stability for large Mach number flow simulations.

\begin{figure}[!h]
\centering
\includegraphics[height=0.42\textwidth]{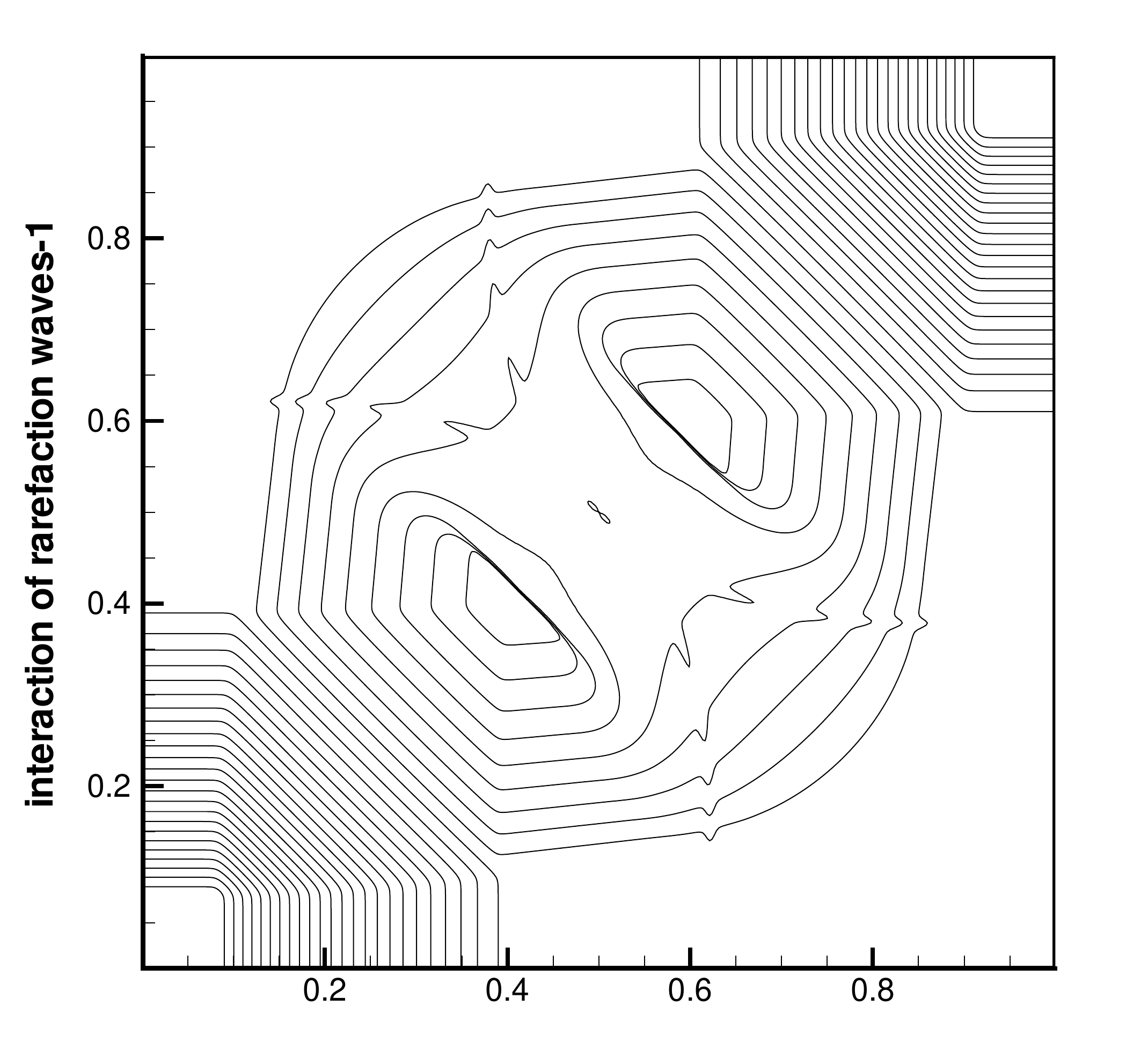}
\includegraphics[height=0.42\textwidth]{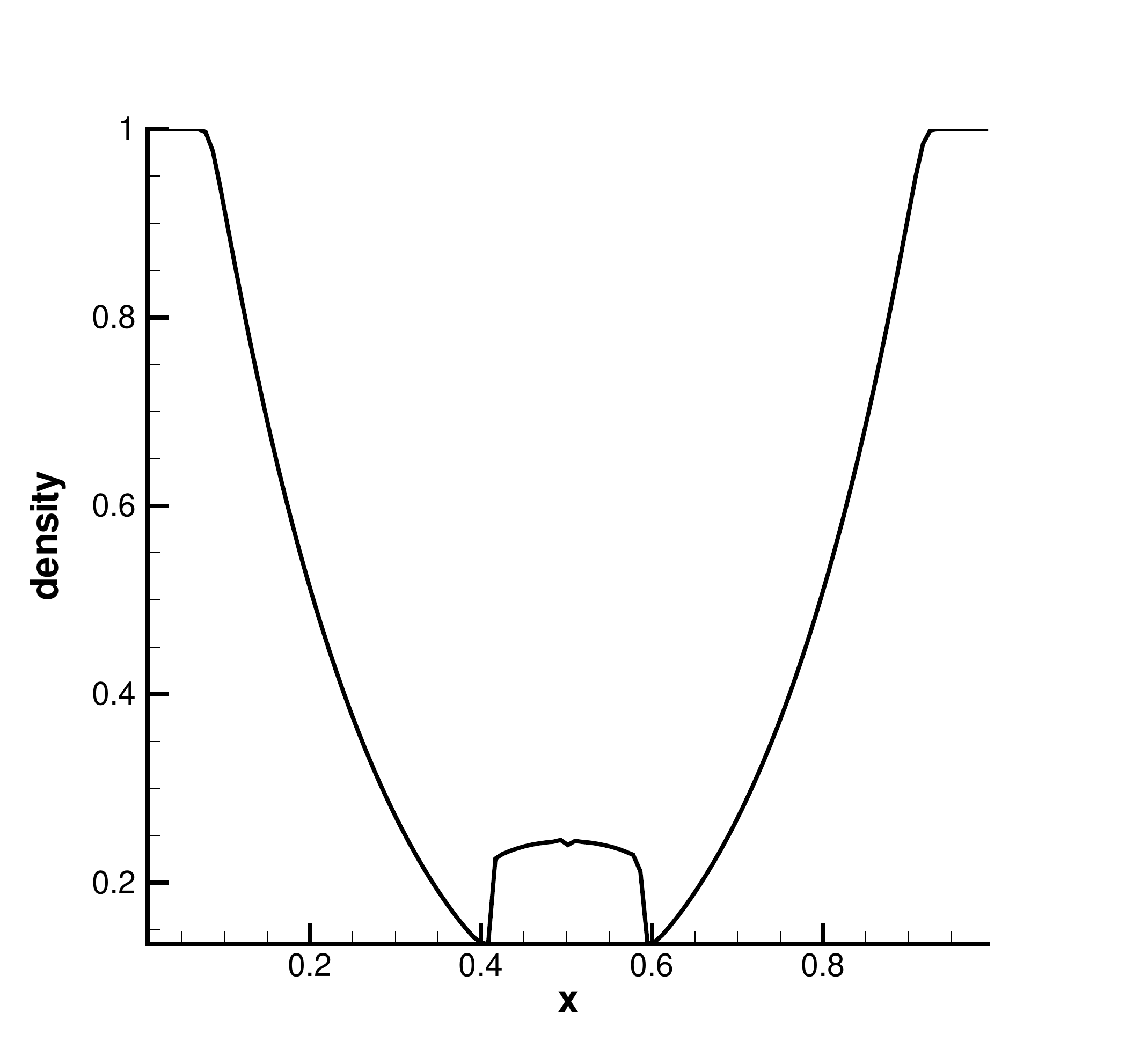}
\caption{\label{Fig:Interaction-Rare1}Interaction of planar
rarefaction waves with initial condition Eq.\eqref{data:rare1}.}
\includegraphics[height=0.42\textwidth]{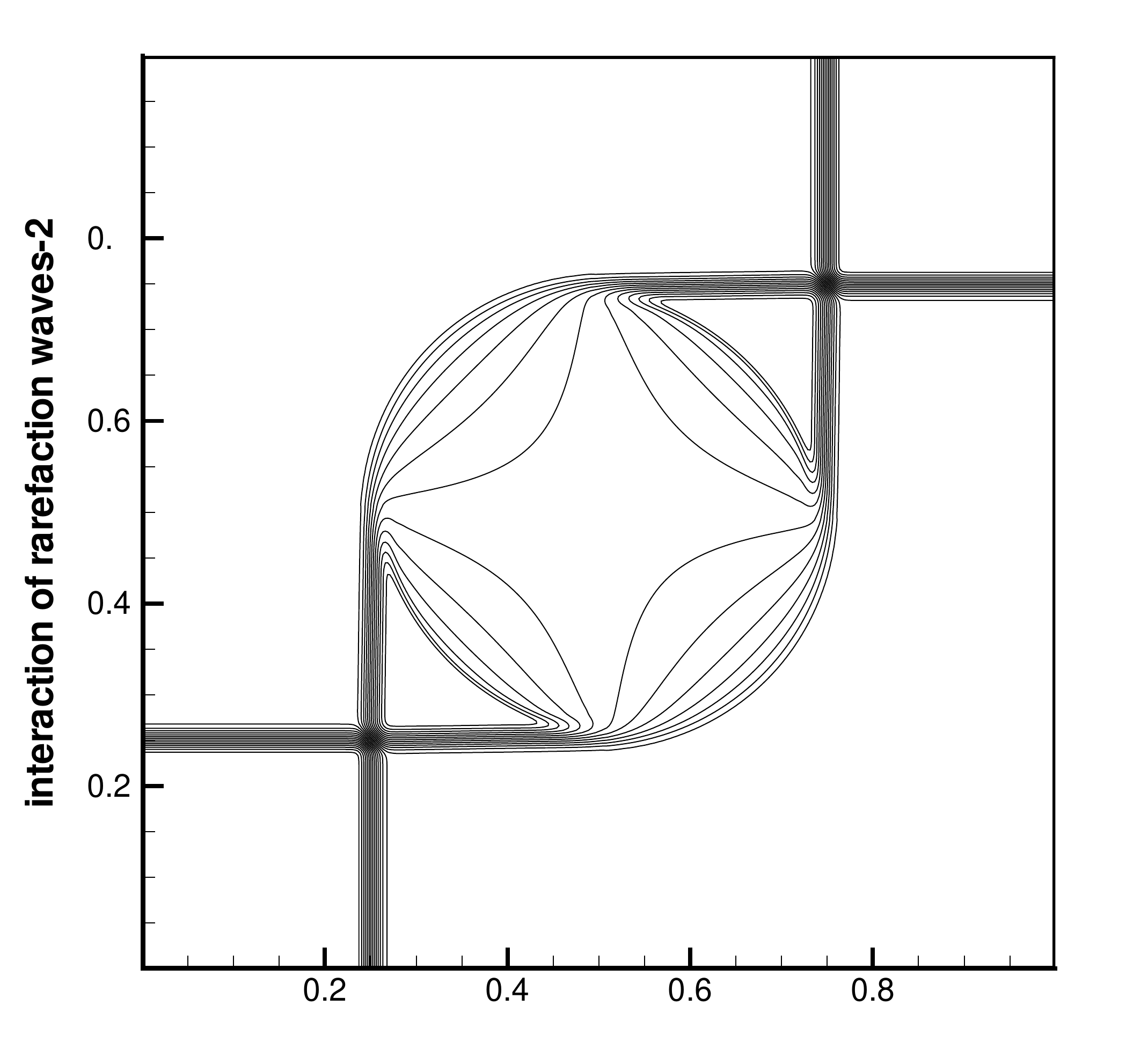}
\includegraphics[height=0.42\textwidth]{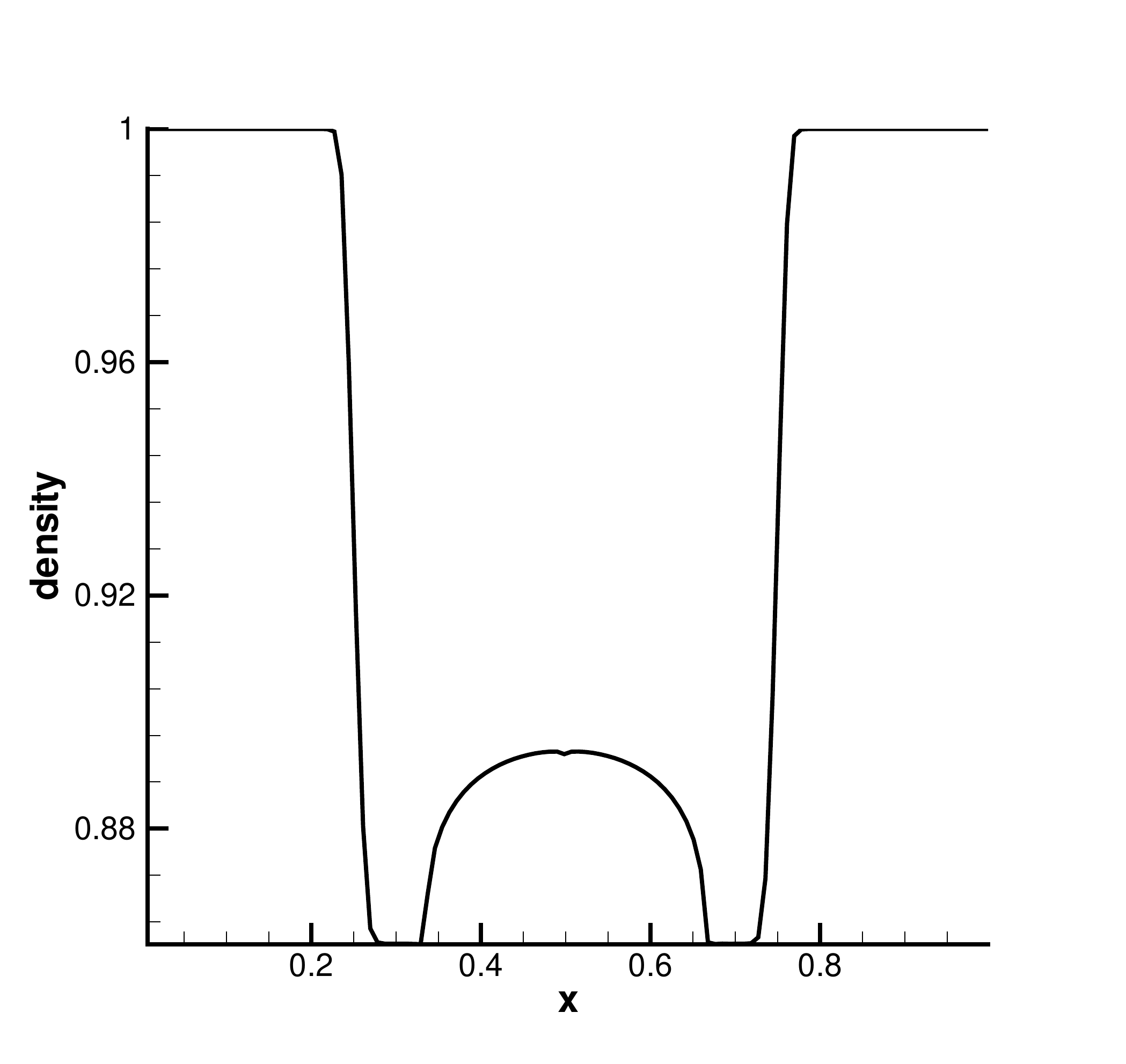}
\caption{\label{Fig:Interaction-Rare2}Interaction of planar
rarefaction waves with initial condition Eq.\eqref{data:rare2}.}
\end{figure}

\subsubsection{Transition from continuous flows to the present of shocks}
In this group, we simulate the interaction of two-dimensional planar
rarefaction waves, from which we observe the continuous transition
from continuous flows to the presence of transonic shocks. A global
continuous solution is constructed in \cite{Li-Zheng} with a clear
picture for the bi-symmetrical interaction only when the rarefaction
waves involved are weak. In general, such a bi-symmetric interaction
may result in the presence of shocks, which never occurs in
one-dimension since the interaction of one-dimensional rarefaction
waves produces only continuous solutions \cite{Courant}. The
backward rarefaction wave and forward rarefaction wave, which
connect the $l$ and $r$ areas, are denoted as
$\overleftarrow{R_{lr}}$ and $\overrightarrow{R_{lr}}$. The four
planar rarefaction waves
$\overrightarrow{R_{21}}\overleftarrow{R_{32}}
\overrightarrow{R_{41}}\overleftarrow{R_{34}}$ with the following
initial conditions are considered
\begin{align}\label{data:rare1}
\begin{cases}
(\rho_1,U_1,V_1,p_1)=(1,0.6233,0.6233,1.5),\\
(\rho_2,U_2,V_2,p_2)=(0.389,-0.6233,0.6233,0.4),\\
(\rho_3,U_3,V_3,p_3)=(1,-0.6233,-0.6233,1.5),\\
(\rho_4,U_4,V_4,p_4)=(0.389,0.6233,-0.6233,0.4),
\end{cases}
\end{align}
and
\begin{align}\label{data:rare2}
\begin{cases}
(\rho_1,U_1,V_1,p_1)=(1,0.0312,0.0312,0.5),\\
(\rho_2,U_2,V_2,p_2)=(0.927,-0.0312,0.0312, 0.45),\\
(\rho_3,U_3,V_3,p_3)=(1,-0.0312,-0.0312,0.5),\\
(\rho_4,U_4,V_4,p_4)=(0.927,0.0312,-0.0312,0.45).
\end{cases}
\end{align}
The numerical results are presented in Fig.
\ref{Fig:Interaction-Rare1} - \ref{Fig:Interaction-Rare2}, in which
the uniform mesh with $\Delta x=\Delta y=1/400$ is used. For the
solution with the initial condition Eq.\eqref{data:rare1}, the
planar rarefaction waves are strong and shocks are in the interior
domain, for which there was a thorough study on how those shocks
occur and on what the criterion is \cite{Glimm-Li}. On the contrary,
if the planar rarefaction waves are relatively weak, e.g. with the
initial data Eq.\eqref{data:rare2}, the global solution is
continuous \cite{Li-Zheng}, but the density in the interior domain
is quite low, as shown in Fig.\ref{Fig:Interaction-Rare2}. So the
simulation on the interaction of planar rarefaction waves is also
multi-scaled and the numerical results in Fig.
\ref{Fig:Interaction-Rare1} - \ref{Fig:Interaction-Rare2} are fully
consistent with the theoretical analysis in \cite{Glimm-Li,
Li-Zheng}.

\subsubsection{Multiscale wave structures resulting from shock wave interactions}

The fourth group deals with the interaction of shocks, in which we
observe the ability of the current scheme capturing the solution
structures of small scales. The study of shock wave interactions is
always one of central topics in gas dynamics and related fields. The
benchmark problems of step-facing shock interaction and the oblique
shock reflection presented in \cite{Woodward-Colella} were carefully
tested with three categories of numerical schemes. Since then,
almost all new designed schemes take those problems to display their
numerical performance, which can be seen from the huge number of the
Google citations. One example from a case of the two-dimensional
Riemann problems in \cite{2d-riemann0} with  four planar shocks
$\overleftarrow{S_{21}}\overleftarrow{S_{32}} \overleftarrow{S_{41}}
\overleftarrow{S_{34}}$ are tested, where the backward rarefaction
wave and forward rarefaction wave, which connect the $l$ and $r$
areas, are denoted as $\overleftarrow{S_{lr}}$ and
$\overrightarrow{S_{lr}}$. To obtain the detailed flow structure
with less computational mesh points, the initial condition are given
as follows
\begin{align}
\left\{\begin{aligned}
&(\rho_1,U_1,V_1,p_1)=(1.5,0,0,1.5), ~~~x>0.8,y>0.8,\\
&(\rho_2,U_2,V_2,p_2)=(0.5323,1.206,0,0.3), ~~~x<0.8,y>0.8\\
&(\rho_3,U_3,V_3,p_3)=(0.138,1.206,1.206,0.029), ~~~x<0.8,y<0.8\\
&(\rho_4,U_4,V_4,p_4)=(0.5323,0,1.206,0.3), ~~~x>0.8,y<0.8.
\end{aligned}
\right.
\end{align}

\begin{figure}[!htb]\centering
\includegraphics[width=0.45\textwidth]{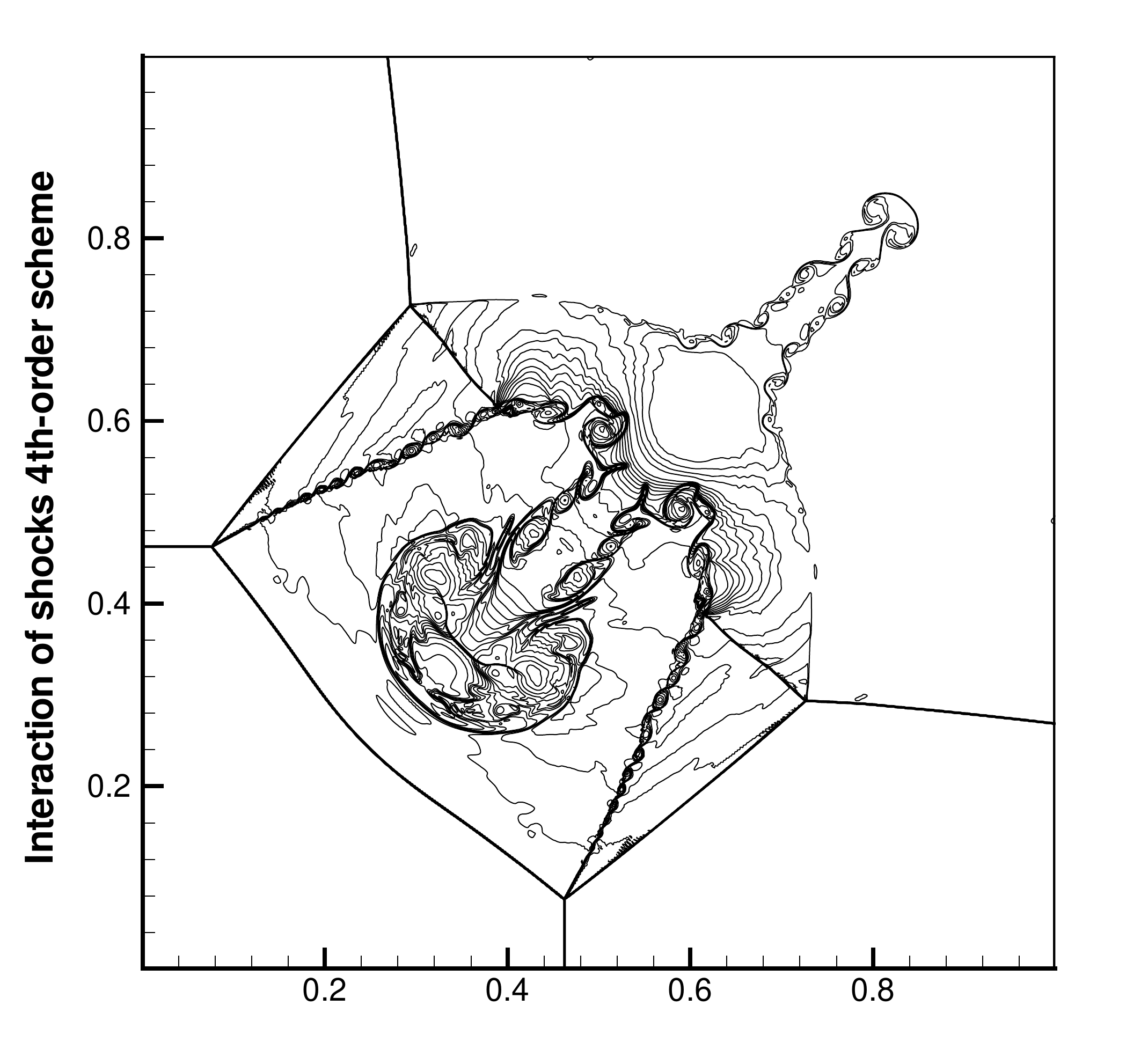}
\includegraphics[width=0.45\textwidth]{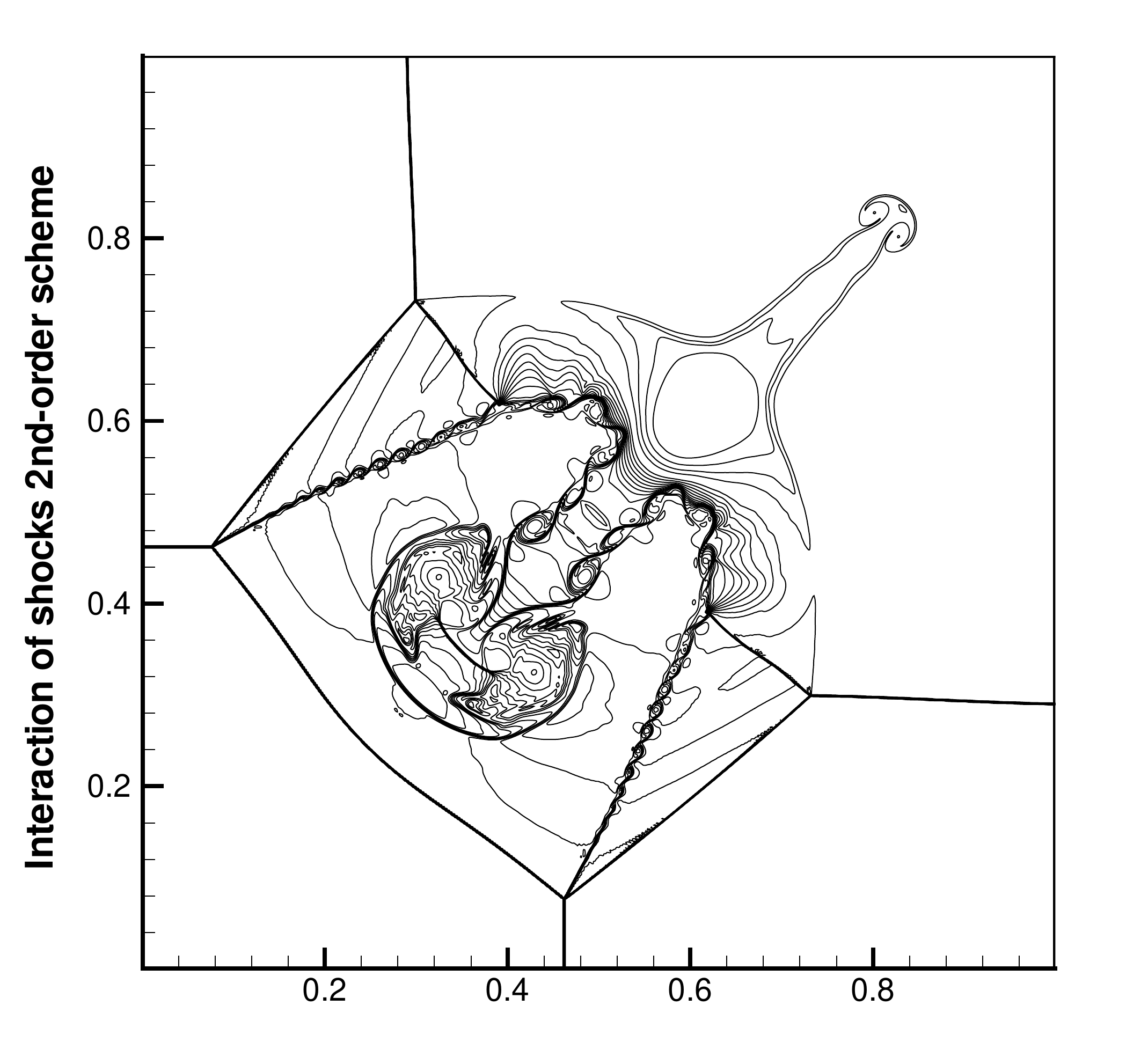}
\caption[small]{The density distribution for the fourth-order scheme
and second-order scheme with $\Delta x=\Delta y=1/1000$ for the
shock wave interactions.} \label{Fig:shock-interaction}
\end{figure}

The numerical solution is displayed in Figure
\ref{Fig:shock-interaction}, where the uniform mesh with $\Delta
x=\Delta y=1/1000$ is used. This case is just the mathematical
formation of the reflection problem of oblique shocks in
\cite{Woodward-Colella} and the symmetric line $x=y$ can be regarded
as the rigid wall. Such a formulation can avoid the complexity of
numerical boundary conditions and make the simulation simple so that
any numerical scheme can perform it easily without worry about
numerical treatment of boundary conditions. The gas-kinetic scheme
with second-order and fourth-order temporal accuracy are tested
using the same fifth-order accurate WENO spatial data reconstruction
at each time level to simulate the wave patters resulting from the
interaction of shocks. It is evident that the second order scheme is
much more dissipative than the four order accurate one. The small
scaled vortices are resolved  sharply using the fourth-order
gas-kinetic scheme, which hints that it may be suitable for the flow
with small structure, such as turbulent flows \cite{Lee-1, Lee-2,
Lee-3}.

\begin{figure}[!h]
\centering
\includegraphics[width=0.23\textwidth]{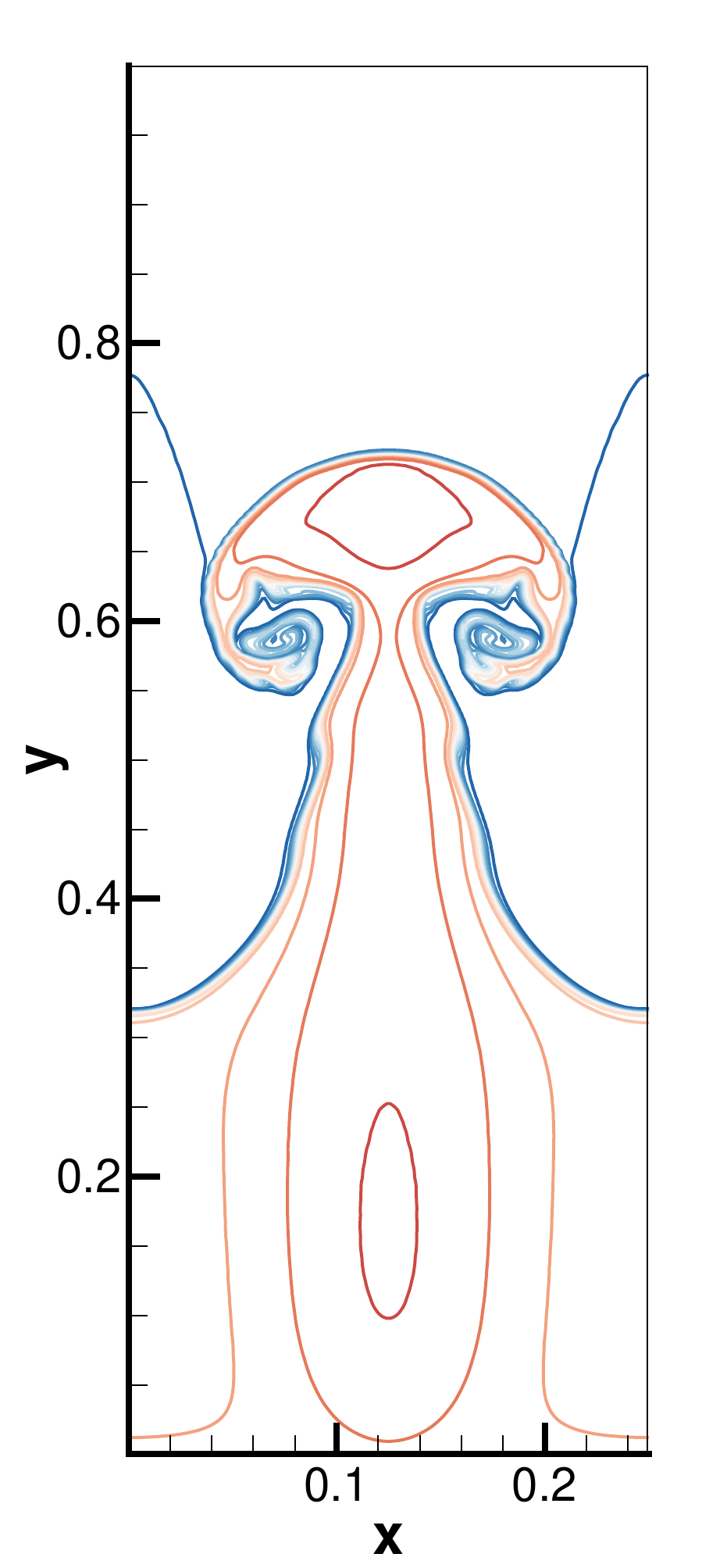}
\includegraphics[width=0.23\textwidth]{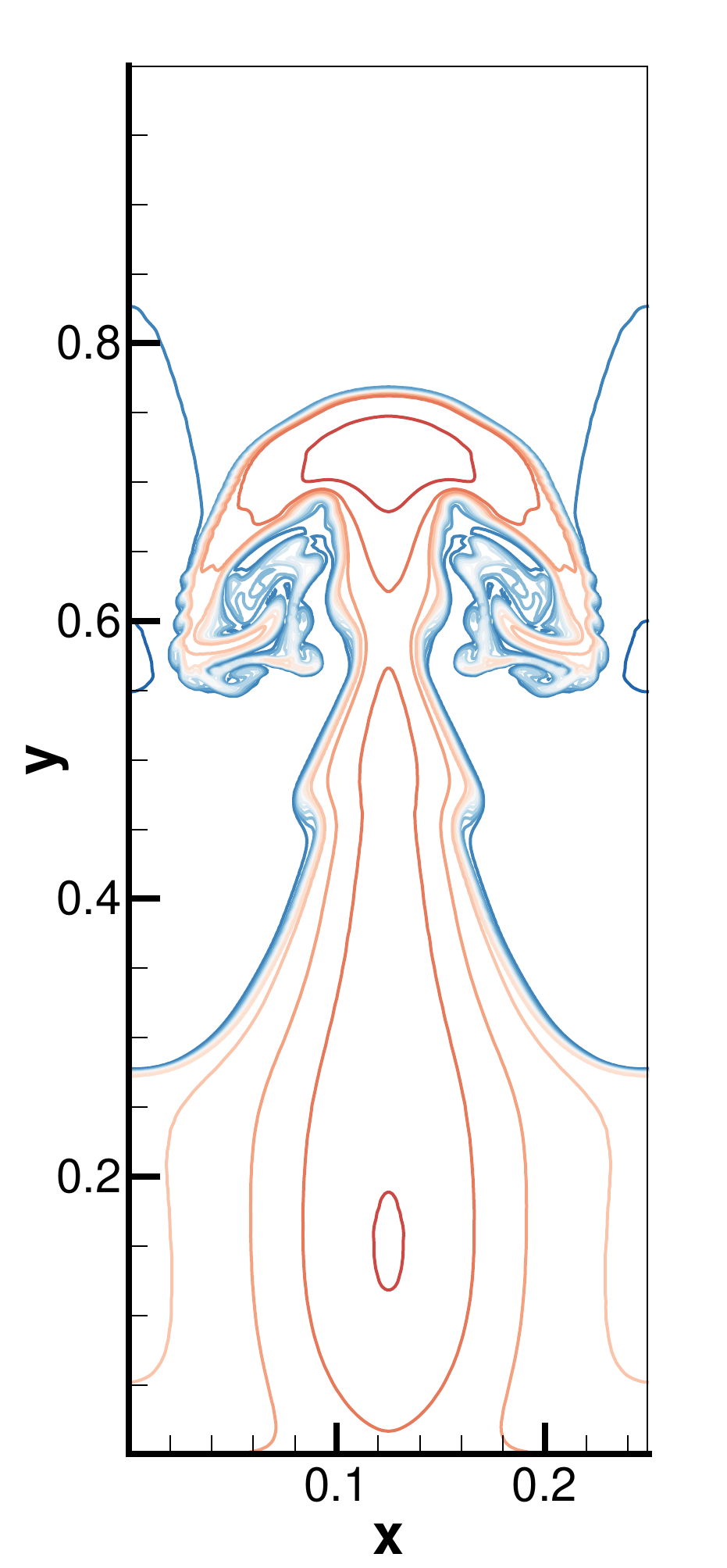}
\includegraphics[width=0.23\textwidth]{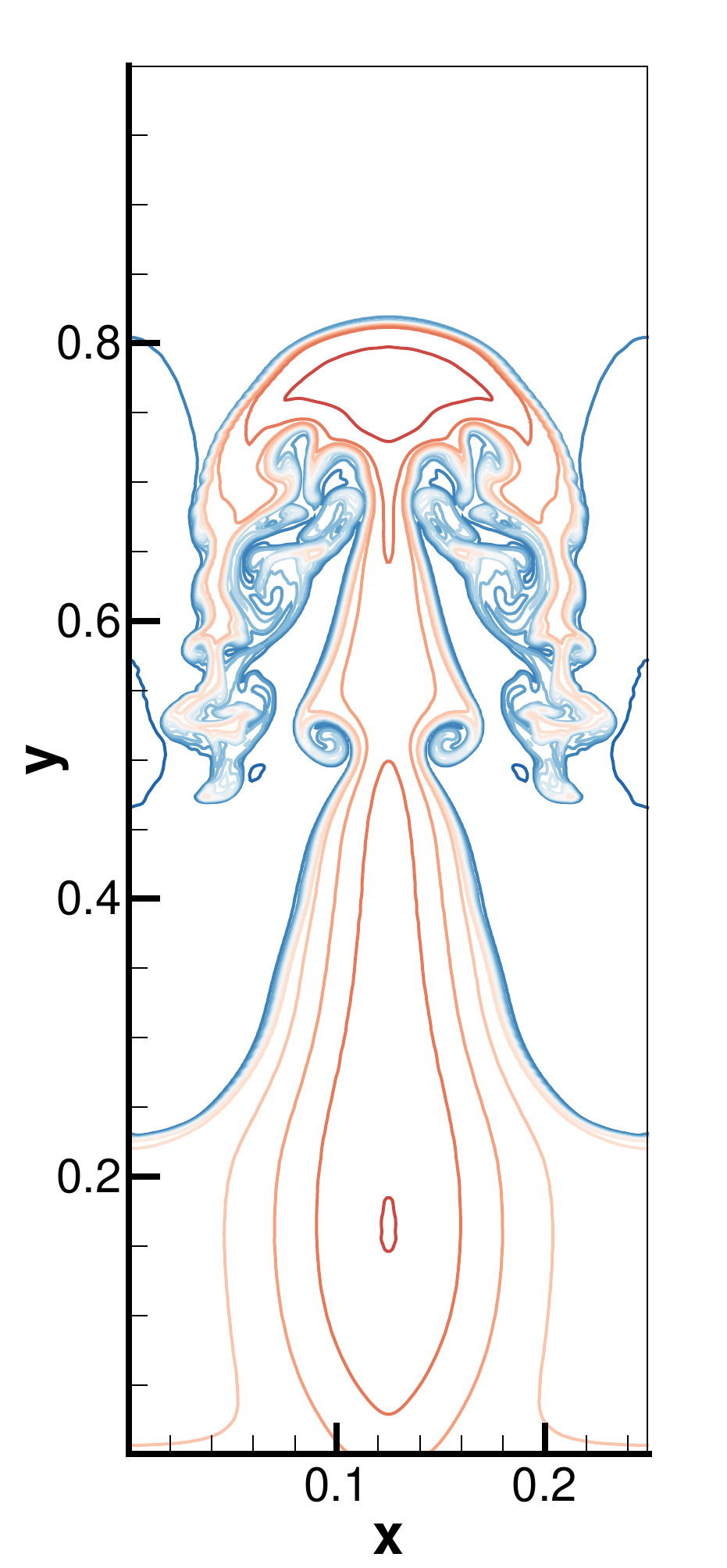}
\includegraphics[width=0.23\textwidth]{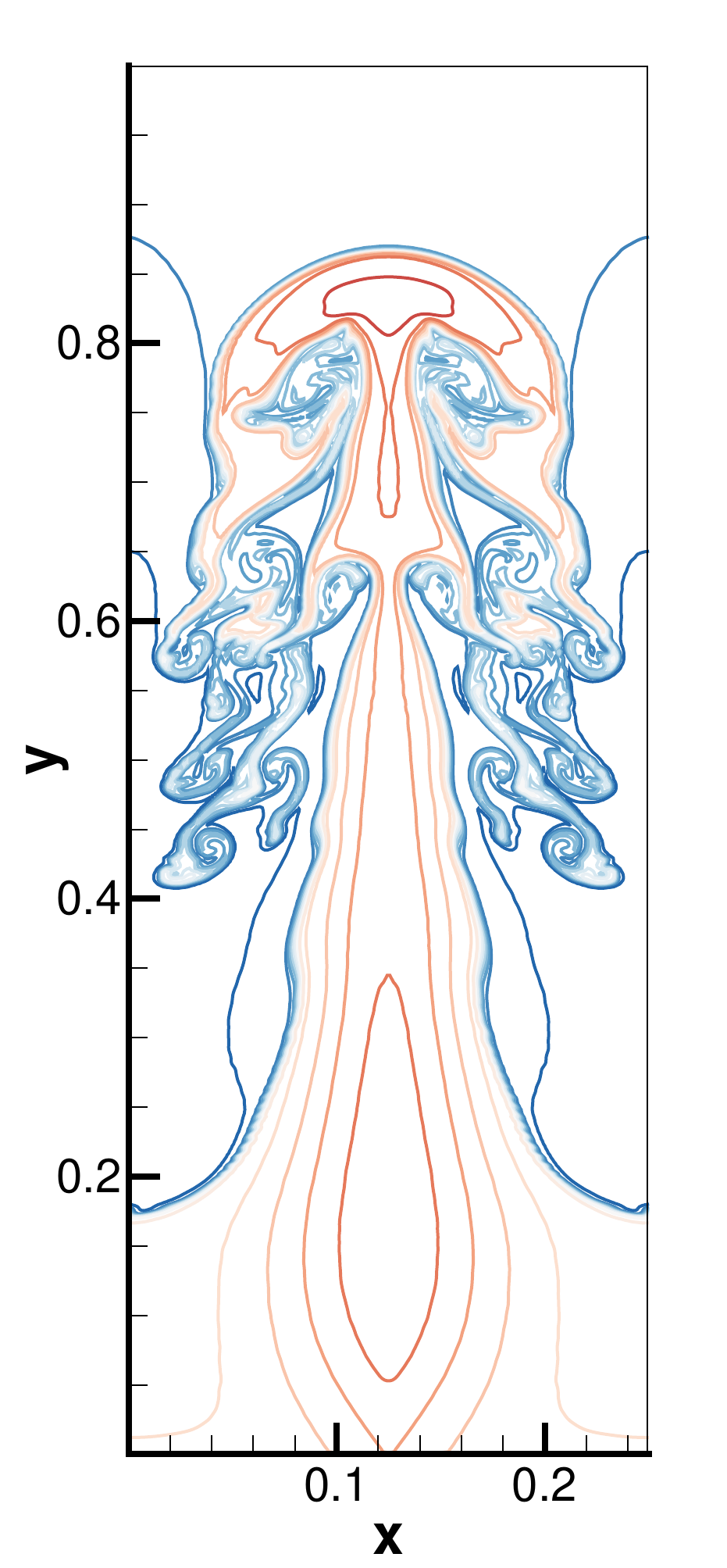}
\caption{\label{Rayleigh-Taylor1} Rayleigh-Taylor instability:
density distribution with the mesh size $\Delta x=\Delta y=1/800$ at
$t=1.75, 2, 2.25$ and $2.5$.}
\end{figure}
\begin{figure}[!h]
\centering
\includegraphics[width=0.23\textwidth]{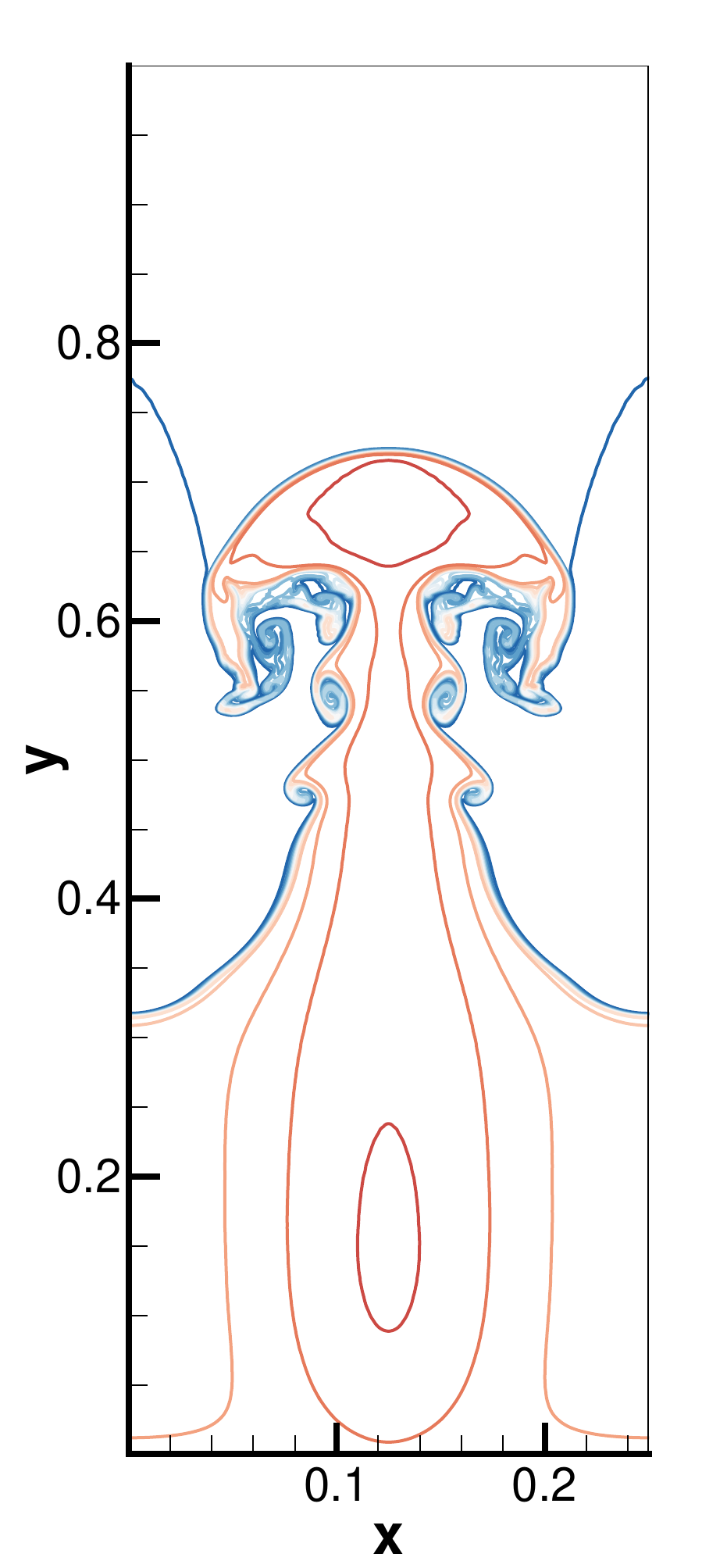}
\includegraphics[width=0.23\textwidth]{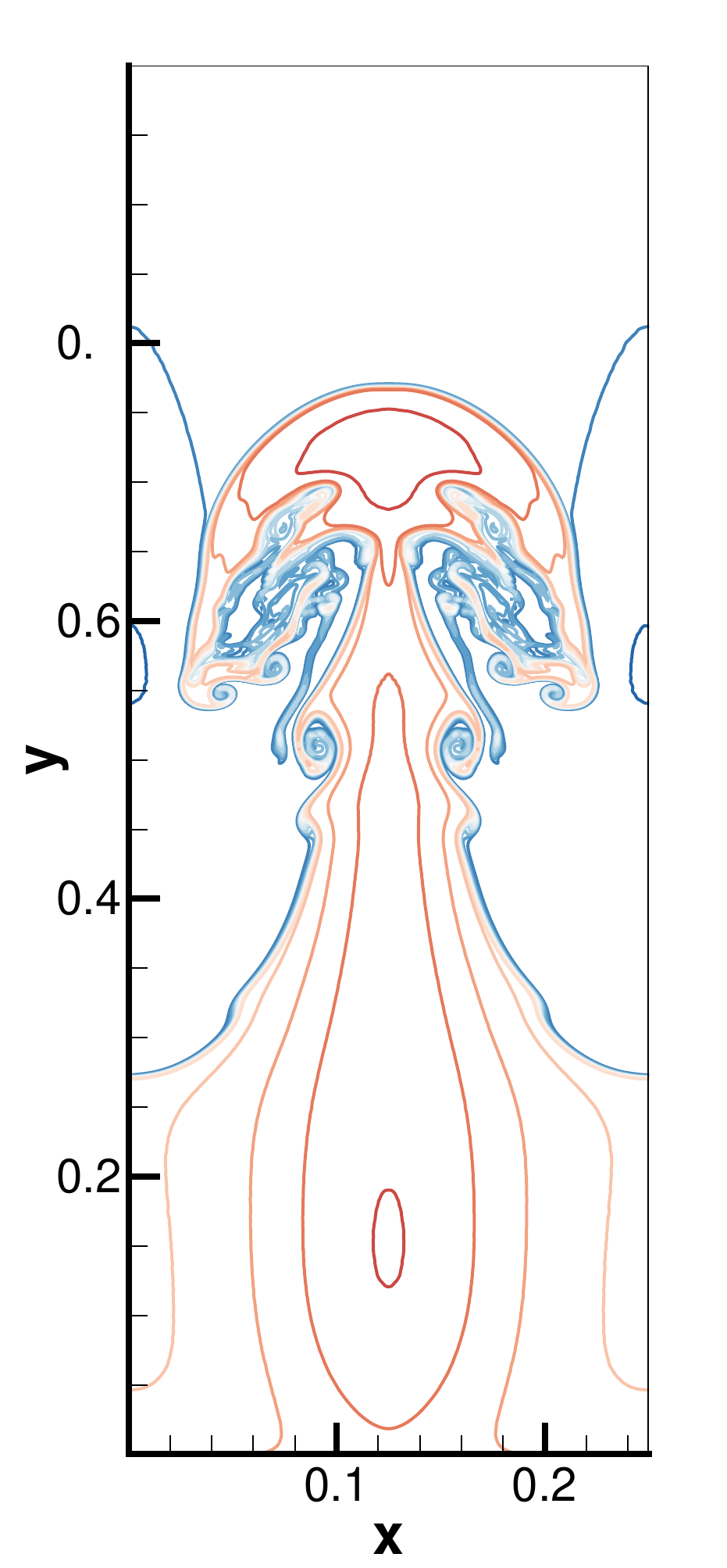}
\includegraphics[width=0.23\textwidth]{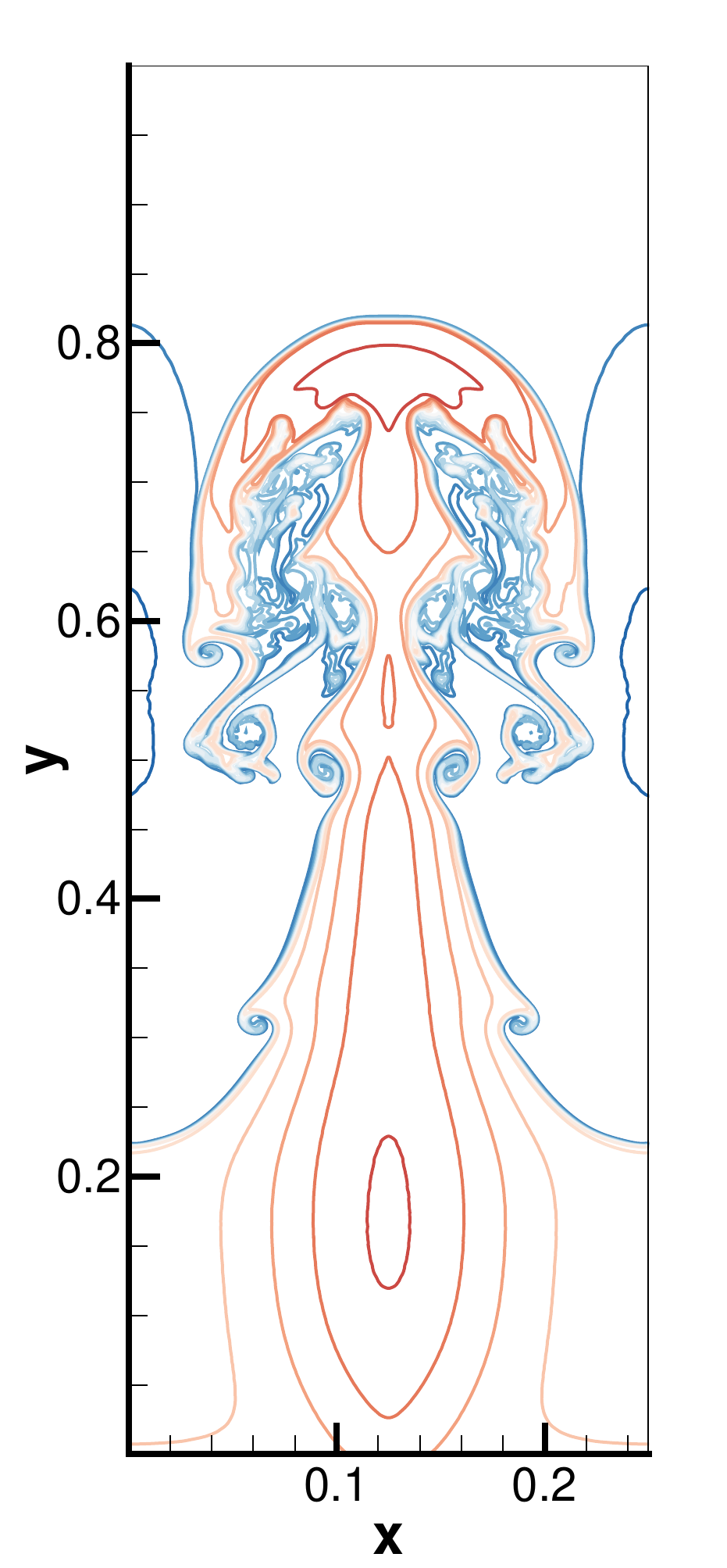}
\includegraphics[width=0.23\textwidth]{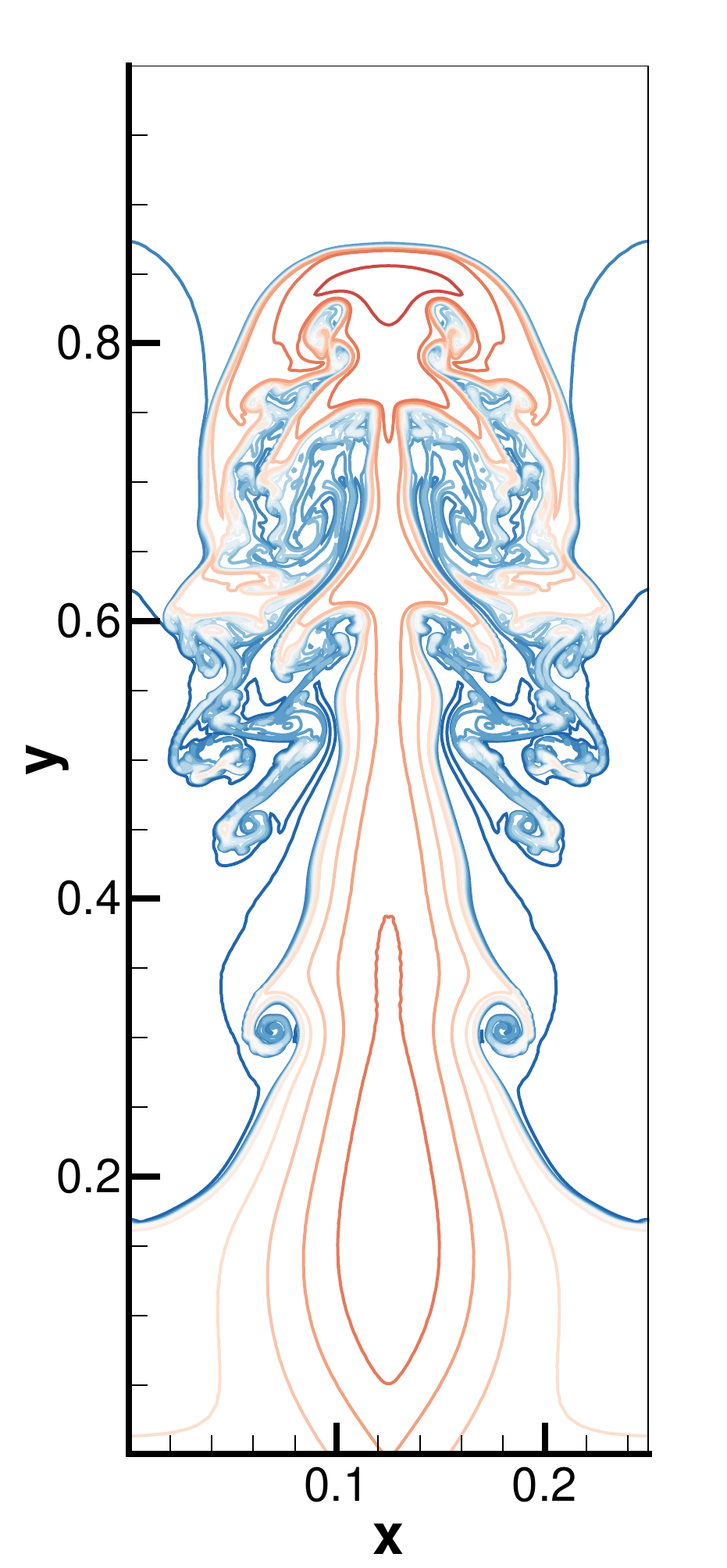}
\caption{\label{Rayleigh-Taylor2} Rayleigh-Taylor instability:
density distribution with the mesh size $\Delta x=\Delta y=1/1600$
at $t=1.75, 2, 2.25$ and $2.5$.}
\end{figure}

\subsection{Conservation laws with source terms}
The last group is the Rayleigh-Taylor instability to test the
performance of the two-stage temporal discretization for the
conservation laws with source terms. Rayleigh-Taylor instability
happens on an interface between fluids with different densities when
an acceleration is directed from the heavy fluid to the light fluid.
The instability has a fingering nature bubbles of light fluid rising
into the ambient heavy fluid and spikes of heavy fluid falling into
the light fluid. The initial condition of this problem
\cite{tr-instability} is given as follows
\begin{align*}
\begin{cases}
(\rho, U, V, p)=(2, 0, -0.025c\cos(8\pi x, 2y+1), x\leq 0.5,\\
(\rho, U, V, p)=(1, 0, -0.025c\cos(8\pi x, y+3/2), x> 0.5,
\end{cases}
\end{align*}
where $c=\sqrt{\displaystyle\frac{\gamma p}{\rho}}$ is the sound
speed and $\gamma=5/3$. The computational domain is
$[0,0.25]\times[0,1]$. The reflective boundary conditions are
imposed for the left and right boundaries; at the top boundary, the
flow values are set as $(\rho, U, V, p)=(1, 0, 0, 2.5)$, and at the
bottom boundary, they are $(\rho, U, V, p)=(2, 0, 0, 1)$. The source
terms on the right side of the governing equations are
$S(\textbf{w})=(0,0,\rho,\rho V)$. To achieve the temporal accuracy,
$S(\textbf{w})$ is given by the cell averaged value $\textbf{w}$ and
$\displaystyle
\partial_tS(\textbf{w})=(0,0,\partial_t\rho,\partial_t(\rho V))$ is
given by the governing equation, respectively.  The uniform meshes
with $\Delta x=\Delta y=1/800$ and $1/1600$ are used in the
computation. The density distributions at $t=1.75, 2, 2.25$ and
$2.5$ are presented in Fig.\ref{Rayleigh-Taylor1} and
Fig.\ref{Rayleigh-Taylor2}. With the mesh refinement, the flow
structures for the complicated flows are observed. It hints that
current scheme may be suitable for the flow with interface
instabilities as well.

\section{Conclusion}
In this paper, we select several  one-dimensional and
two-dimensional challenging problems which can be used to check the
performance of higher order numerical schemes. The scheme for
providing reference solutions is our recently developed two-stage
four order accurate GKS scheme. These cases can be used to test the
performance of other  higher-order schemes as well, which have been
intensively constructed in recent years. The test problems are

\begin{itemize}
\item One-dimensional problems
\begin{enumerate}
\item Titarev-Toro's highly oscillatory shock-entropy wave
interaction;
\item Large density ratio problem with a very strong rarefaction wave;
\end{enumerate}
\item Two-dimensional Riemann problems
\begin{enumerate}
\item  Hurrican-like solutions with one-point vacuum and rotational velocity field;
\item Interaction of planar contact discontinuities, with the involvement of
entropy wave and vortex sheets;
\item Interaction of planar rarefaction waves with the transition from  continuous
flows to the presence of shocks;
\item Interaction of planar (oblique) shocks.
\end{enumerate}
\item Conservation law with source terms
\begin{enumerate}
\item Rayleigh-Taylor instability
\end{enumerate}
\end{itemize}
The construction of accurate and robust higher-order numerical
schemes for the Euler equations is related to many numerical and
physical modeling issues. Based on the exact Riemann solution or
other simplified approximate Riemann solvers,  the traditional
higher-order approaches  mostly concentrate on the different
frameworks, such as DG, finite different, and finite volume, with
all kind of underlying data reconstructions or limiters. The
higher-order dynamics under the flux function has no attracted much
attention. Maybe this is one of the reasons for the stagnation on
the development of higher-order schemes in recent years. Along the
framework of multiple derivatives for the improvement of time
accuracy of the scheme, it has a higher requirement on the accuracy
of the flux modeling than the traditional higher-order methods,
because the scheme depends not only on the flux, but also on the
time derivative of the flux function. Most test cases in this paper
are for the time accurate solutions, which require the close
coupling of the space and time evolution. These test cases can fully
differentiate the performance of different kinds of higher-order
schemes. Their tests and the analysis of the numerical solutions can
guide the further development of higher-order schemes on a
physically and mathematically consistent way.

\section*{Acknowledgements}
 The work of J. Li is supported by NSFC (11371063, 91130021),  the doctoral program from the Education Ministry of China (20130003110004) and the Science Challenge Program in China Academy of Engineering Physics. The research of K.  Xu is supported by Hong Kong Research Grant Council (620813, 16211014, 16207715) and HKUST research fund
(PROVOST13SC01, IRS15SC29, SBI14SC11).

\end{document}